\documentclass[12pt,leqno]{amsart}  
\usepackage{amscd,amsfonts,amssymb,amsthm,amsmath}
\usepackage[mathscr]{eucal}
\usepackage[english]{babel}
\usepackage{color}

\usepackage{amscd}

\usepackage[small,nohug,heads=littlevee]{diagrams} 
\diagramstyle[labelstyle=\scriptstyle,LaTeXeqno]

\theoremstyle{change}

\newtheorem{theorem}{Theorem}[section]
\newtheorem{proposition}[theorem]{Proposition}
\newtheorem{definition}[theorem]{Definition}
\newtheorem{lemma}[theorem]{Lemma}
\newtheorem{corollary}[theorem]{Corollary}
\newtheorem{remark}[theorem]{Remark}
\newtheorem{remarks}[theorem]{Remarks}
\newtheorem{example}[theorem]{Example}
\newtheorem{examples}[theorem]{Examples}
\newtheorem{conjecture}[theorem]{Conjecture}
\newtheorem{notation}[theorem]{Notation}

\newenvironment{demo}{\noindent{\bf Proof.}}{
{{\phantom{a}\hfill\rule{4pt}{6pt}}}}

\def\endproof{
{{\phantom{a}\hfill\rule{4pt}{6pt}}}}

\def\CC{\mathbb{C}}
\def\NN{\mathbb{N}}
\def\QQ{\mathbb{Q}}
\def\RR{\mathbb{R}}
\def\ZZ{\mathbb{Z}}
\def\Aut{\hbox{Aut}}

\def\triv{\{1\}}
\def\adele{{\mathbb A}}
\def\Kprod{{\mathbb O}}

\def\RCom#1#2{
            #1/\kern-3pt / #2
}
\def\Rmap#1#2{
            \tau_{#1, #2}
}

\def\field{{\Bbb{F}}}
\def\symm{\hbox{\rm Symm}}

\title[Commensurated Subgroups of Arithmetic Groups]{Commensurated 
Subgroups of Arithmetic Groups,
 Totally Disconnected Groups
and Adelic Rigidity}

\author{Yehuda Shalom}
\address{Department of Mathematics, UCLA, Los Angeles CA, USA}
\email{yeshalom@math.ucla.edu}

\author{George A. Willis}
\address{School of Mathematical and Physical Sciences, University of 
Newcastle, Building V, Callaghan 2308 AUSTRALIA }
\email{George.Willis@newcastle.edu.au}

\begin{document}

\maketitle


\section{Introduction}
\label{sec:intro}

\noindent {\bf The Margulis-Zimmer conjecture.} The subject of this paper is a well known
question advertised by Gregory Margulis and Robert Zimmer since the late 1970's, 
which seeks refinement of the 
celebrated Normal Subgroup Theorem of Margulis (hereafter abbreviated 
{\bf NST}). 
Although Margulis' NST is stated and proved in the context of (higher rank) 
irreducible lattices in products of simple algebraic groups over local fields, by Margulis'
arithmeticity theorem we may and shall work solely 
in the framework of
(S-){\it arithmetic groups}. One departure point for the Margulis-Zimmer 
conjecture is the phenomenon
that while all higher rank 
$S$-arithmetic groups are uniformly treated by the NST, 
there is a notable difference in the structure of subgroups which are {\it commensurated}, rather than normalized, by the ambient arithmetic group. 
For example, the group 
$SL_n(\ZZ[{\frac{1}{p}}])$ commensurates its subgroup $SL_n(\ZZ)$, while the
latter commensurates no apparent 
infinite, infinite index subgroup of its own. The obvious generalization of this example, which by Margulis arithmeticity theorem and with the aid of the 
restriction of scalars functor is the most general one, goes
as follows:

\begin{definition}
\label{def:setup}
Let $K$ be a global field, $\mathcal O$ its ring of integers, and let $\bf G $ be an absolutely simple, simply 
connected algebraic group defined over $K$. Let $V$ be the set of all inequivalent valuations  on $K$, and let $V^\infty \subset V$ denote the 
archimedean
ones. For a subset $V^\infty \subseteq S \subseteq V$, let  
$\mathcal O _S \subseteq K$
be, as usual, the ring of $S$-integers in $K$, and let $\Gamma <{\bf G} (K) $
be an $S$-arithmetic group, namely, a subgroup commensurable 
with ${\bf G}(\mathcal O _S)$.
Then any $S'$-arithmetic
subgroup $\Lambda <\Gamma$ ($V^\infty \subseteq S'\subseteq S$) is 
commensurated  by
$\Gamma$, and we call 
such $\Lambda$ a {\bf standard commensurated subgroup}. 
We say that the $S$-arithmetic group 
$\Gamma$ has a {\bf standard
description of commensurated subgroups} if every $\Lambda <\Gamma$ 
commensurated by it, is standard or finite.
\end{definition}

A precise definition of the notions appearing in Definition 
\ref{def:setup} is given in Section~\ref{sec:Proof1.4} below.
The Margulis-Zimmer {\bf commensurated subgroup problem}, which we shall
hereafter abbreviate {\bf CmSP} (to distinguish from the celebrated CSP -- Congruence Subgroup
Problem), can now be stated as:

\medskip
\noindent{\bf The Commensurated Subgroup Problem:} Let $K$, $\bf G$, $S$, $\Gamma$ be as above, and assume, as in the NST, that
$S$-rank$(\Gamma):=\Sigma_{\nu \in S}K_{\nu}$-rank$({\bf G}(K_\nu)) \ge 2$. 
Does
$\Gamma$ have a standard description of commensurated subgroups? 

\medskip

While this was never put in writing, Margulis and Zimmer conjectured that 
the CmSP should have positive answer (certainly in characteristic zero, 
which is studied here), in which case we shall again use the abbreviation
CmSP (where P stands for ``Property''). As we shall see, the notational 
similarity to the CSP is not 
merely formal; the two are in fact intimately linked. Indeed, in 
the last section of this paper we shall
propose a rather sweeping conjecture, which unifies the CmSP and many
other well known deep results and conjectures in the theory of
 arithmetic groups including the
CSP, thus providing additional motivation for its study (we return to this issue towards the end of the introduction).
 We further remark that the assumption that $\bf G$ be simply connected is needed only when $S$ is infinite. For reasons to become clear later on, in all of 
our results in the sequel we shall assume that $K$ is a number field, i.e. 
char$(K)=0$,
 although for illustration we shall  also use arithmetic 
groups in positive characteristic.

An important caveat which must be made at this point is that 
the CmSP is not formally well defined as stated. Indeed, an $S$-arithmetic 
subgroup
of ${\bf G} (K)$ is in fact only defined up to commensurability 
(as it depends on the
$K$-embedding of $\bf G$ in $\mathbf{GL_n}$), while, as we shall see in Theorem~\ref{thm:counter} below, the 
property we are after is in general sensitive to passing to finite index
subgroups. Thus, merely proving that every subgroup of 
$\Gamma=SL_3(\ZZ)$ commensurated by it is finite or co-finite 
(which is indeed the case)
does not, in itself, guarantee the same result for any finite index subgroup of 
$\Gamma$.
This is one motivation for the approach we take in this paper, to which we now turn, where a much stronger property which {\emph is} commensurability-stable  
is introduced.

\medskip

\noindent {\bf Some key notions and the main results.}  Since one can naturally view the CmSP as a refinement
of Margulis' NST, a natural approach to it would be to try to push further, 
or give a ``better'' 
proof of Margulis' NST. Extending Margulis' original analytic proof
(the only one that exists in complete generality) turns out to encounter 
serious difficulties and has not become successful. Expanding on 
algebraic approaches to some special cases of the NST, T. N. Venkataramana in the only prior work \cite{Venkat} around the CmSP, 
established a somewhat ``brute
force'' computational proof of the CmSP for 
$\Gamma = {\mathbf G}(\ZZ)$ when $\mathbf G$ is defined over $\QQ$ and 
$\QQ$-rank$({\mathbf G}) \ge2$ 
(e.g. $\Gamma= SL_{n \ge 3}(\ZZ)$ -- see also 
the beginning of Section~\ref{sec:alternative}).
 The approach we take in this paper is very different conceptually, as it completely separates the 
CmSP and the NST by providing the tool to reduce the former to the latter 
(taken as a ``black box'').
However, when applicable, this approach provides a much stronger 
(fixed point) property, new even for  $SL_{n \ge 3}(\ZZ)$.
To explain it we introduce the following notions, where 
groups are assumed countable.

\begin{definition}
\label{def:comm-normalizer}
{\bf 1.} We say that a group $\Gamma$ has the {\bf inner} 
{\it commensurator-normalizer property}, if every
commensurated subgroup $\Lambda< \Gamma$ is {\it almost normal} in the sense
(used throughout this paper) that such
$\Lambda$ is  commensurable in $\Gamma$ with a normal subgroup of $\Gamma$. \newline
{\bf 2.} Say that $\Gamma$ has the {\bf outer} {\it commensurator-normalizer property} if the following
holds: for {\it any} group $\Delta$ and any 
homomorphism $\varphi:\Gamma \to \Delta$, any subgroup
$\Lambda < \Delta$ which is commensurated by (the conjugation action of) 
$\varphi(\Gamma)$, is almost 
normalized by  $\varphi(\Gamma)$, namely, a subgroup commensurable with 
$\Lambda$ in $\Delta$ is normalized by $\varphi(\Gamma)$.

\end{definition} 

Clearly if an arithmetic group has the inner (let alone the outer) property,
then by Margulis' NST every commensurated subgroup is either
finite or co-finite which, as noted above,  is not the case for 
general $S$-arithmetic groups.
Two crucial advantages of the outer property
are that, first, {\it unlike the inner property} it is stable
under commensurability, and second, its strength enables one to automatically 
answer positively the CmSP for all the $S$-arithmetic groups containing
a fixed arithmetic group possessing it. We remark in passing that the outer
commensurator-normalizer property seems remarkably rare. While it may certainly be that many groups (including ``random'' ones) can have the inner property
(e.g., the ``Tarski monster'', in which every proper subgroup is finite), 
it is not clear at all why the same might hold (let alone be proved) 
for the outer property. As one instance of this phenomenon we note in Remark~\ref{rem:positivechar} that 
infinite linear groups in positive characteristic never have the outer
commensurator-normalizer property (although at least some higher rank lattices do have the inner one).

We can now state the main result of the paper:

\begin{theorem}
\label{thm:mainresult}
 
Retain the notations in Definition \ref{def:setup}, and assume ${\mathbf G}$
is a Chevalley (i.e. split) group over $K$ (as always, char$(K)=0$). 
In case ${\mathbf G} \cong {\mathbf{SL}}_2$ assume further that $K$ is not $\QQ$ or an imaginary
quadratic extension of it. Then:

\noindent{\bf 1.} Any group commensurable with ${\mathbf G}({\mathcal O})$ has 
the outer commensurator-normalizer property.

\noindent{\bf 2.} For any $V^\infty \subseteq S \subseteq V$, any 
$S$-arithmetic subgroup of ${\mathbf G}(K)$ has standard description
of commensurated subgroups.
\end{theorem}

One major ingredient in the proof of {\bf 1} is the rich structure theory of 
arithmetic Chevalley groups, particularly their bounded generation 
property (which rests on deep analytic number theoretic tools).
In fact, since this property is known to hold also in the quasi-split case,
the result continues to hold in this case as well (see Theorem~\ref{thm:LMZ}
below for the precise border of our approach). The second main ingredient of 
the proof relies on elements from the structure theory of
automorphisms of locally compact totally disconnected groups, developed by the
second author in a series of papers since the mid 1990's. Although this introduction
revolves mostly around arithmetic groups and the CmSP, a large part of the
paper is in fact devoted to what we believe are interesting results for
their own sake in that theory. Briefly, we show that every
polycyclic group of automorphisms of a l.c.t.d. group $H$ is virtually {\it flat}, which implies that its commutator subgroup (virtually) 
normalizes a compact open
subgroup of $H$.
The connection between the latter and the CmCP is based on a general strategy
to the problem proposed early on by Margulis-Zimmer, which involves topologization process described in Section~\ref{sec:commsub}, fundamental for
this paper,  and was never made to work. In fact part {\bf 1} of 
Theorem~\ref{thm:mainresult} is actually proved first in its topological   
version, a setting which seems essential for the proof, and it 
is the Margulis-Zimmer
strategy that transfers the topological result into the discrete setting
we need.

\medskip

\noindent{\bf The CmSP revisited.} In the last section of the paper 
we change gear, showing how the CmSP combines naturally with
the congruence subgroup problem, the normal subgroup theorem, superrigidity, 
the Margulis-Platonov conjecture, and results originating from work 
of Deligne (extended by Raghunathan and well known for their use by Toledo) 
concerning
the non-residual finiteness of certain ``lifts'' of arithmetic groups to universal covers of algebraic groups. These deep aspects of the theory of arithmetic groups will be shown to fit together as different 
facets of one unified 
conjecture in adelic framework,
building upon work of Deligne (and supported by 
results of Prasad-Rapinchuk on the metaplectic kernel) concerning central extensions of algebraic 
groups over local and global fields.
While we postpone precise details and references to that section, we state
here a simplified version of an {\it unconditional} result, in the direction of this conjecture, which is deduced
from Theorem \ref{thm:mainresult}.

\begin{theorem}
\label{thm:adelic}
Retain the assumptions and notations of Theorem~\ref{thm:mainresult} and Definition~\ref{def:setup}. Assume that $K$ admits a field embedding into $\RR$. 
Let $\adele _f$ denote the ring of finite adeles of $K$, 
 ${\bf G}(\adele _f)$ denote, as usual, the 
restricted
direct product of ${\bf G}(K_\nu)$ over the finite $\nu \in V$,
and consider an $S$-arithmetic subgroup 
$\Gamma <{\bf G} (K)$ identified with its image in ${\bf G}(\adele _f)$
through the diagonal embedding. 
Let $\varphi:\Gamma \to H$ be {\it any} homomorphism
into an {\it arbitrary}  locally compact totally disconnected group $H$. 
Then one, and exactly one of the following occurs: either

\noindent{\bf 1.} Im $\varphi$ is discrete and Ker $\varphi $ is finite (central), in which case $\varphi$ doesn't extend, or

\noindent{\bf 2.} The homomorphism $\varphi$ extends 
to a continuous homomorphism of the closure\\
 $\bar \Gamma < {\bf G}(\adele _f)$ {\it onto}
 the closure ${\overline {\varphi (\Gamma)}} < H$. 
  
Furthermore, if the normalizer of any compact open subgroup of $H$ is compact
(or even merely amenable), then {\bf 2} necessarily holds.
\end{theorem}

Of course, as one can take $H=\Gamma$ in the theorem (or ``encode'' this 
situation
in a non-discrete $H$), {\bf 1} is an inevitable possibility.
 Note that as before, $S$ may be infinite. 
The assumption on $K$ arises from the intimate 
relation to  the congruence subgroup property. Indeed, the proof of 
Theorem \ref{thm:adelic} makes crucial use of the CSP, NST and our
solution above to the CmSP. 
In return, by making different choices of $H$ it immediately 
implies all of them, 
as well
as superrigidity (noting that the condition on $H$ in the last statement is 
satisfied when $H=GL_n(F)$). In the proof of Theorem \ref{thm:adelic}
 in Section~\ref{sec:Proof1.4} we describe also 
what happens exactly when  $K$ is totally imaginary, which yields 
a modified result. The general conjecture, which builds on deep work
of Deligne~\cite{De1}, anticipates a
result in the spirit of Theorem \ref{thm:adelic} where
all places, finite and infinite, play a symmetric role, and assumptions
are relaxed considerably.

\medskip

\noindent {\bf Acknowledgments.} 
The authors thank Gregory Margulis and
Robert Zimmer for their input to the manuscript.
Deep and special gratitude goes to Gopal Prasad for the wealth of 
crucial information he provided during the work on
Section 7.2, and to Andrei Rapinchuk who has read a previous version of this section
with great care, pointing
out many improvements. The first author acknowledges the support
of the ISF and NSF through grants 500/05 and DMS-0701639 resp. The second author 
acknowledges the
 ARC support through grants
LX0667119 and DP0984342.

\section{Preliminaries}
\label{sec:prelim}

Throughout this section, $G$ will be a totally disconnected locally compact group and the group of bi-continuous automorphisms of $G$ will be denoted by $\hbox{Aut}(G)$. Every such $G$ has a base of neighborhoods of the identity consisting of compact open subgroups, see \cite[\S2.3]{MontZip}. The set of all compact open subgroups of $G$ will be denoted by  ${\mathcal B}(G)$. Of course, 
despite having many compact open subgroups, $G$ need not have such a subgroup that is normal. 
Should $x\in G$ fail to normalize any compact open subgroup of $G$, there will nevertheless be subgroups that  are minimally distorted by the inner automorphism $\alpha_x : y\mapsto xyx^{-1}$ $(y\in G)$. 

\begin{definition}
\label{def:tidy&flat}
Let $G$ be a totally disconnected locally compact group and let $\alpha\in \hbox{Aut}(G)$ and $x\in G$. 
\begin{enumerate}
\item The \emph{scale of $\alpha$} is the positive integer
$$
s(\alpha) = \min\left\{ [\alpha(V) : V\cap \alpha(V)] \mid V \in {\mathcal B}(G)\right\}.
$$
(The \emph{scale of  $x$} is the positive integer $s(x) = s(\alpha_x)$.)
\item A compact open subgroup $V$ such that $s(\alpha) = [\alpha(V) : V\cap \alpha(V)] $ is called \emph{minimizing for $\alpha$} (and \emph{minimizing for $x$} when it is minimizing for $\alpha_x$).
\end{enumerate}
\end{definition}

The scale function for automorphisms is analogous to a certain function on linear transformations, namely, the absolute value of the product of all eigenvalues whose absolute value is greater than $1$. Indeed, $s(\alpha)$ may be calculated in that way from the eigenvalues of $\hbox{ad}(\alpha)$ when  $G$ is a $p$-adic Lie group, see \cite{{Gl:skew},Gl:padic}. An automorphism such that $s(\alpha) = 1 = s(\alpha^{-1})$ may therefore be expected to behave like a unipotent linear transformation. The properties of the scale listed in the next theorem are consistent with this analogy (consider the absolute value of the determinant in place of the modular function in \ref{thm:scale}(\ref{determinant})). All properties may derived directly from Definition~\ref{def:tidy&flat}, except for \ref{thm:scale}(\ref{powers}) and continuity, which require the structure of minimizing subgroups described below in Theorem~\ref{thm:tidiness} and a further argument given in \cite[Theorem~3]{Wi94}. 

\begin{theorem}[Properties of the Scale]
\label{thm:scale}
Let $\alpha\in\Aut(G)$. The scale function\\ 
$s: \Aut(G) \to \ZZ^+$ has the following properties.
\begin{enumerate}
\item $s(\alpha) = 1 = s(\alpha^{-1})$ if and only if there is $V\in {\mathcal B}(G)$ such that $\alpha(V) = V$. \label{uniscalar}
\item $s(\alpha^n) = s(\alpha)^n$ for every $n\in\NN$. \label{powers}
\item $\Delta(\alpha) = s(\alpha)/s(\alpha^{-1})$, where $\Delta : \Aut(G) \to (\RR^+,\times)$ is the modular function on $\Aut(G)$.\label{determinant}
\item $s(\beta\alpha\beta^{-1})$ for every $\beta\in \Aut(G)$. \label{conjclass}
\end{enumerate}
In addition, the scale function $s : G\to \ZZ^+$ induced on $G$ by inner automorphisms is continuous with respect the group topology on $G$ and the discrete topology on $\ZZ^+$.
\end{theorem}

Continuing the linear algebra analogy for the scale, the set of subgroups minimizing for $\alpha$ corresponds to a triangularizing basis for a linear transformation. It is seen in \cite{Gl:skew,Gl:padic} that, when $G$ is a $p$-adic Lie group, minimizing subgroups may indeed be described in terms of a triangularizing basis for $\hbox{ad}(\alpha)$. There is a close association between minimizing subgroups and a canonical form also when $G$ is the automorphism group of a homogeneous tree, for, if $x\in G$ is hyperbolic, then the subgroups of $G$ minimizing for $\alpha_x$ are the stabilizers of strings of vertices on the axis of $x$, see~\cite{Wi94}. The structural characterization of minimizing subgroups given in the next theorem lends the analogy substance. 
\begin{theorem}[The Structure of Minimizing Subgroups]
\label{thm:tidiness}
Let $\alpha\in \hbox{Aut}(G)$. For each $V\in {\mathcal B}(G)$ put
$$
V_+ = \bigcap_{k\geq0} \alpha^k(V) \ \hbox{ and }\ V_- = \bigcap_{k\geq0} \alpha^{-k}(V).
$$
Then $V$ is minimizing for $\alpha$ if and only if
\begin{description}
\item[TA($\alpha$)] $V = V_+V_-$ and 
\item[TB($\alpha$)] $V_{++} := \bigcup_{k\geq0} \alpha^k(V_+)$ is closed.
\end{description}
If $V$ is minimizing for $\alpha$, then $s(\alpha) = [\alpha(V_+):V_+]$. 
\end{theorem}

A compact open subgroup satisfying {\bf TA($\alpha$)} is said to be \emph{tidy above for $\alpha$}, while a subgroup satisfying {\bf TB($\alpha$)} is \emph{tidy below}. A compact open subgroup that is tidy above and below is said to be \emph{tidy for $\alpha$}. Tidy subgroups were first defined in~\cite{Wi94} in relation to the conjecture of K. H. Hofmann and A. Mukherjea concerning concentration functions that was stated in~\cite{HM} and solved in~\cite{JRW}. Equivalence of the tidiness and minimizing properties was established in~\cite{Wi01}.

The proof of the structure theorem relies on the following procedure that, given an arbitrary $V\in {\mathcal B}(G)$, modifies it in three steps to produce a subgroup satisfying {\bf TA($\alpha$)} and {\bf TB($\alpha$)}. 
{\paragraph{\bf Step 1:} Choose $n\in {\mathbf N}$ such that $\bigcap_{k=0}^n \alpha^k(V)$ satisfies {\bf TA($\alpha$)}. (That such $n$ exists is shown in~\cite[Lemma~1]{Wi94}.) Set $V' := \bigcap_{k=0}^n \alpha^k(V)$. Then 
$$
[\alpha(V') : \alpha(V')\cap V'] \leq [\alpha(V) : \alpha(V)\cap V]
$$
with equality if and only if $V$ already satisfies {\bf TA($\alpha$)}.
\paragraph{\bf Step 2:} Find a compact, $\alpha$-stable subgroup $K_\alpha$ such that $V_{++}K_\alpha$ is closed. (Two ways of finding $K_\alpha$ that are relevant in Section~\ref{sec:flatness} are discussed below.)
\paragraph{\bf Step 3:} Set $W = \left\{ x\in V' \mid xK_\alpha \subset K_\alpha V'\right\}$ and $V'' = WK_\alpha$. Then $V''$ is a compact open subgroup and satisfies {\bf TA($\alpha$)} and {\bf TB($\alpha$)}.
Furthermore,
$$
[\alpha(V'') : \alpha(V'')\cap V''] \leq [\alpha(V') : \alpha(V')\cap V']
$$
with equality if and only if $V'$ already satisfies {\bf TB($\alpha$)} (and $K_\alpha \leq V'$).
}
\bigskip

Steps 1 and 3 modify $V$ without increasing the index $[\alpha(V) : \alpha(V)\cap V]$. Since, by \cite[Theorem~2]{Wi94}, this index is the same for all subgroups satisfying {\bf TA($\alpha$)} and {\bf TB($\alpha$)}, it must be the minimum value. Since the index is strictly decreased by these modifications unless $V$ already satisfies {\bf TA($\alpha$)} and {\bf TB($\alpha$)}, such subgroups are the only ones where the minimum is attained.

The subgroup $K_\alpha$ found in {\bf Step 2} is the obstruction to $V'$ being tidy in the sense that a group satisfying {\bf TA($\alpha$)} is tidy if and only if it contains $K_\alpha$. In \cite{Wi01} this subgroup is defined in terms of $V'$ by putting
\begin{equation}
\label{eq:define1K}
 {\mathcal K}_\alpha := \left\{ x\in G \mid \alpha^k(x) \in V' \hbox{ for almost every }k\right\} \hbox{ and }K_\alpha^{(1)} := \overline{\mathcal K}_\alpha. 
\end{equation} 
For an argument in \cite{Wi04} it is important that $K_\alpha$ be defined independently of $V'$. This is done in the following way, which may yield a different group. For each $V\in{\mathcal B}(G)$ define  
\begin{align} 
{\mathcal K}_{\alpha,V}  := {}&  \left\{ x\in G \mid \{\alpha^k(x)\} \hbox{ is bounded and } \alpha^k(x)\in V,  \forall k \hbox{ sufficiently large}\right\} \notag\\
\hbox{and } K_\alpha^{(2)} := {}&\bigcap_{V\in {\mathcal B}(G)} \overline{{\mathcal K}_{\alpha,V}}.
\label{eq:define2K}
\end{align}

A compact open subgroup minimizing for an automorphism $\alpha$ may act in proofs as a substitute for an $\alpha$-stable subgroup when these do not exist. Here are some facts about $\alpha$-minimizing subgroups that can be used in place of corresponding facts about $\alpha$-stable subgroups. 
\begin{proposition}
\label{prop:recall}
Let $V$ be minimizing for $\alpha$ and let $n\in\ZZ$. Then:
\begin{enumerate}
\item \label{eq:poweralpha}
$V$ is minimizing for $\alpha^n$;
\item \label{eq:applyalpha}
$\alpha^n(V)$ is minimizing for $\alpha$;
\item \label{eq:tidyprop}
$V\cap \alpha^n(V) = \bigcap_{k=0}^n \alpha^k(V)$; and 
\item  \label{eq:tidyintersect} if $V_1$ and $V_2$ are minimizing for $\alpha$, then so is $V_1\cap V_2$.
\end{enumerate}
\end{proposition}
\begin{demo}
Parts \ref{eq:poweralpha}--\ref{eq:tidyprop} follow immediately either from the definition of minimizing subgroups or from their tidiness. Part~\ref{eq:tidyintersect} requires more work and is \cite[Lemma~10]{Wi94}. 
\end{demo}

One important difference between $\alpha$-stable and $\alpha$-minimizing subgroups is that the set of automorphisms that have $V$ as a minimizing subgroup might fail to be a group because it fails to be closed under multiplication. For example, if $G$ is the automorphism group of a homogeneous tree and $V$ is the stabilizer of the adjacent edges, $v_1$ and $v_2$, then $V$ is minimizing for any elliptic $x\in G$ that leaves $\{v_1,v_2\}$ invariant and any hyperbolic $x\in G$ whose axis contains $v_1$ and $v_2$, but these automorphisms do not form a set that is closed under multiplication. Should a \emph{subgroup} of $\Aut(G)$ have a common minimizing subgroup, the conclusion of Theorem~\ref{thm:tidiness} can be improved. 
\begin{definition}
\label{def:flat}
Let ${\mathcal H}\leq \hbox{Aut}(G)$.
\begin{enumerate}
\item The subgroup $V\in{\mathcal B}(G)$ satisfies  {\bf TA(${\mathcal H}$)} if it satisfies {\bf TA($\alpha$)} for every $\alpha\in {\mathcal H}$, and is \emph{minimizing for ${\mathcal H}$} if it is minimizing for every $\alpha\in {\mathcal H}$.
\item ${\mathcal H}$ is  \emph{flat} if there is a compact open subgroup $V$ that is minimizing for ${\mathcal H}$. 
\item The \emph{uniscalar subgroup} of  the flat group ${\mathcal H}$ is 
$$
{\mathcal H}(1) := \left\{\alpha\in {\mathcal H} \mid s(\alpha) = 1 = s(\alpha^{-1}) \right\}.
$$
\end{enumerate}
\end{definition}
That ${\mathcal H}(1)$ is a subgroup of ${\mathcal H}$ follows from Theorem~\ref{thm:scale}(\ref{uniscalar}). It is not difficult to show that it is in fact a normal subgroup. 

The analogy between automorphisms and linear transformations suggests that flatness of a group of  automorphisms should be equivalent to commutativity modulo the uniscalar subgroup. It is proved in \cite{Wi04} that this 
is indeed the case. In one direction, finitely generated abelian groups are flat. This criterion for flatness is strengthened below in Section~\ref{sec:flatness}. The next theorem states the converse direction and fleshes out the analogy with linear transformations. 
\begin{theorem}
\label{thm:flat}
Let ${\mathcal H}\leq \Aut(G)$ be a flat group of automorphisms and suppose that ${\mathcal H}/{\mathcal H}(1)$ is finitely generated. Let $V$ be minimizing for ${\mathcal H}$. 
\begin{enumerate}
\item There is $r\in {\NN}$ such that
$$
{\mathcal H}/{\mathcal H}(1) \cong \ZZ^r.
$$ \label{eq:flatrankr}
\item There are $q\in {\NN}$ and closed subgroups, $V_j$ for $j\in\left\{0,1,\dots,q\right\}$, of $V$ such that
$$
V = V_0V_1 \cdots V_q,
$$
where for every $\alpha\in {\mathcal H}$ we have $\alpha(V_0) = V_0$ and, for every $j\in \left\{1,\dots,q\right\}$, either $\alpha(V_j)\geq V_j$ or $\alpha(V_j)\leq V_j$. \label{qnumber}
\item For each $j\in \left\{1,\dots,q\right\}$ the group $\tilde V_j := \bigcup_{\alpha\in {\mathcal H}} \alpha(V_j)$ is  closed and $\alpha$-stable in $G$. 
\item For each $j\in \left\{1,\dots,q\right\}$ let $\Delta_j : {\mathcal H} \to (\QQ^+,\times)$ denote the modular function for the restriction of ${\mathcal H}$ to $\tilde V_j$. Then:
\begin{enumerate}
\item $\Delta_j(\alpha) = \begin{cases}
[\alpha(V_j) : V_j], & \hbox{ if }\alpha(V_j)\geq V_j\\
[V_j : \alpha(V_j)]^{-1}, & \hbox{ if } \alpha(V_j)\leq V_j
\end{cases}$; and \label{module}
\item there are $s_j\in \ZZ^+$ and a homomorphism $\rho_j : {\mathcal H} \to \ZZ$ such that 
$$
\Delta_j(\alpha) = s_j^{\rho_j(\alpha)}\hbox{ for every }\alpha\in{\mathcal H}.
$$
\end{enumerate} \label{module}
\item $s(\alpha) = \prod\left\{ \Delta_j(\alpha) \mid \rho_j(\alpha) >0\right\}$ for every $\alpha\in {\mathcal H}$. \label{factorscale}
\end{enumerate}
\end{theorem}

The groups $\tilde V_j$ are analogous of common eigenspaces for elements of ${\mathcal H}$ with the numbers $\Delta_j(\alpha)$ being the corresponding eigenvalues. Obtaining the factoring of $V$ in part~\ref{qnumber} is the main part of the proof and it is achieved by applying Theorem~\ref{thm:tidiness} to a sequence $\alpha_1$, $\alpha_2$, \dots of elements of ${\mathcal H}$: $\alpha_1$ gives a first factoring of $V$, then $\alpha_2$ is chosen to refine it, and so on. This argument has a geometric flavor, where the homomorphisms $\rho_j$ are viewed as `roots' on ${\mathcal H}/{\mathcal H}(1)$ and the $\alpha_k$'s chosen determine hyperplanes that separate the roots. Sufficiently many hyperplanes must be chosen to separate all the roots and that can be many more than the number of generators of ${\mathcal H}/{\mathcal H}(1)$, as is shown by the example where $G = \QQ_p^q$ and ${\mathcal H}$ is generated by the automorphisms $\alpha((x_j)) := (px_j)$ and $\beta((x_j)) := (p^jx_j)$. That there are finitely many roots and the process terminates may be deduced from the fact that ${\mathcal H}$ is finitely generated in the following way. Let $\{\beta_k\}_{k\in{\sf K}}$ is a finite generating set for ${\mathcal H}$. Each factoring of $V$ determines a factoring of the scales $s(\beta_k)$, and when the factoring of $V$ is properly refined by a new hyperplane, so too is the factoring of at least one $s(\beta_k)$. (This is a precursor to the factoring of scales in part~\ref{factorscale}.) Since we are factoring integers from a finite set, there comes a point when all roots have been separated. 

The exponent $r$ appearing in~part \ref{eq:flatrankr} is the \emph{flat-rank} of ${\mathcal H}$ and is always less than or equal to the number $q$ of `eigenspaces' or `roots'. In general, there is no inequality in the reverse direction as the example ${\mathcal H} = \langle \alpha,\beta\rangle \leq \Aut(\QQ_p^q)$ mentioned in the previous paragraph shows. However: when the flat-rank is~$0$, then $q=0$ and ${\mathcal H}$ is uniscalar; and when the flat-rank is~$1$, then $q =1$ or $q=2$. Theorem~\ref{thm:tidiness} may, with hindsight, be thought of as the flat-rank~1 case of Theorem~\ref{thm:flat}. Part~\ref{factorscale} of Theorem~\ref{thm:flat} has the following consequence for the scale function on flat groups.
\begin{corollary}
\label{cor:submultiplicative}
Let ${\mathcal H}\leq \Aut(G)$ be a flat group of automorphisms. Then
\begin{enumerate}
\item the scale function is submultiplicative on ${\mathcal H}$, that is, $s(\alpha\beta) \leq s(\alpha)s(\beta)$ for every $\alpha,\,\beta\in {\mathcal H}$; and
\item  the function
$$
\alpha{\mathcal H}_1 \mapsto \log s(\alpha) + \log s(\alpha^{-1}) : {\mathcal H}/{\mathcal H}(1) \to \RR^+
$$
is a norm on ${\mathcal H}/{\mathcal H}(1)$. 
\end{enumerate}
\end{corollary}
 It is not true that the scale is submultiplicative on groups that are not flat, for the product of two automorphisms with scale~1 can have arbitrary scale. A couple of examples are: the product of two elliptic elements in the automorphism group of a tree can be hyperbolic; and the product of two unipotent elements in a $p$-adic matrix group can be non-unipotent. Note that 
 $$
 \log s(\alpha) + \log s(\alpha^{-1})= \sum \left\{ \left|\log \Delta_j(\alpha) \right|\mid j\in\{1,\dots,q\}\right\}
 $$
 
The following is  \cite[Theorem~3.3]{Wi04}.
\begin{lemma}
\label{lem:joinK}
Let ${\mathcal H}\leq \hbox{Aut}(G)$ be flat and let $V$ be minimizing for  ${\mathcal H}$. Suppose that $K$ is a compact subgroup of $G$ such that $\alpha(K) = K$ for all $\alpha\in {\mathcal H}$. Let
\begin{equation}
\label{eq:joinK}
V' := \left\{ x\in V \mid xK \subset KV\right\}.
\end{equation}
Then $V'$ is a compact open subgroup of $G$ and $V'K$ is a compact open subgroup that is minimizing for ${\mathcal H}$.
\end{lemma}

\section{Commensurated Subgroups and the Topological Completion}
\label{sec:commsub}

Given a group $\Gamma$ and a commensurated subgroup $\Lambda < \Gamma$, there is a natural `completion' of $\Gamma$ that has the closure of $\Lambda$ as a compact open subgroup. This totally disconnected locally compact `completion' has been used to study group actions and representations in \cite{Sch1,Tz,Tz1,GlWi}. Here, it allows the techniques of the previous section to be applied to the study of commensurators and its relevant properties are collected below. 
\begin{definition}
Let $\Gamma$ be a group with subgroup  $\Lambda$ that is commensurated by~$\Gamma$. The left action of $\Gamma$ on $\Gamma/\Lambda$ yields a homomorphism $\Gamma \to \symm(\Gamma/\Lambda)$, denoted by $\Rmap{\Gamma}{\Lambda}$. \\
The \emph{relative profinite completion of $\Gamma$ with respect to $\Lambda$}, denoted $\RCom{\Gamma}{\Lambda}$, is the closure of $\Rmap{\Gamma}{\Lambda}(\Gamma)$ in the topology of pointwise convergence on $\symm(\Gamma/\Lambda)$.
\end{definition}

\begin{example}
\label{ex:Commensurate}
Let $\Gamma = A \rtimes {\mathbb Z}$, where $A :=  \bigoplus_{k\in {\mathbb Z}} C_2$, $C_2 = \{\bar0,\bar1\}$ is the cyclic group of order~$2$ and ${\mathbb Z}$ acts on  $\bigoplus_{k\in {\mathbb Z}} C_2$ by translation. ($\Gamma$ is sometimes called the \emph{lamplighter group}.) Let $\Lambda = \left\{ f\in A \mid f(0) = \bar0\right\}$. Then $\Lambda$ is commensurated by $\Gamma$ and $\RCom{\Gamma}{\Lambda} \cong  \left(\prod_{k\in {\mathbb Z}} C_2\right) \rtimes {\mathbb Z}$. 
\end{example}

\begin{remarks}
\noindent {\bf 1.} The subgroups $\left\{\Rmap{\Gamma}{\Lambda}(\gamma\Lambda\gamma^{-1}) \mid \gamma\in \Gamma \right\}$ form a subbase of neighborhoods of the identity in the topology of pointwise convergence on $\Rmap{\Gamma}{\Lambda}(\Gamma)$. 

\noindent {\bf 2.} If $\Lambda$ is a normal subgroup of $\Gamma$, then $\ker(\Rmap{\Gamma}{\Lambda}) = \Lambda$ and the topology of pointwise convergence on $\symm(\Gamma/\Lambda)$ is discrete, so that $\RCom{\Gamma}{\Lambda}\cong \Gamma/\Lambda$.  

\noindent {\bf 3.} The relative profinite completion is not, strictly speaking, a completion of $\Gamma$ unless $\ker(\Rmap{\Gamma}{\Lambda})$ is trivial. 
\end{remarks}

The $\Lambda$-orbits of the action of $\Gamma$ on $\Gamma/\Lambda$ are finite and $\overline{\Rmap{\Gamma}{\Lambda}(\Lambda)}$ is open because $\gamma\Lambda \ne \Lambda$  when $\gamma\not\in \Lambda$. These observations lead immediately to the following. 
\begin{proposition}
\label{prop:profcomp}
The set $\overline{\Rmap{\Gamma}{\Lambda}(\Lambda)}$ is a compact totally disconnected group, an open subgroup of $\RCom{\Gamma}{\Lambda}$. It follows that $\RCom{\Gamma}{\Lambda}$ is a totally disconnected, locally compact group. \endproof
\end{proposition}

In the reverse direction, suppose that $G$ is a topological group with compact, open subgroup $L$. Then $L$ is commensurated by $G$. 

\begin{lemma}
\label{lem:topcase}
Let $G$ be a topological group and $L$ be a compact, open subgroup of $G$. Then $\Rmap{G}{L} : G \to \symm(G/L)$ is a continuous and open map with closed range. Hence $\RCom{G}{L}\cong G/\ker(\Rmap{G}{L})$, where $\ker(\Rmap{G}{L})$ is the largest normal subgroup of $G$ contained in~$L$. It follows that $\RCom{G}{L}$ is isomorphic to $G$ when the kernel is trivial.
\end{lemma}
\begin{demo}
Continuity of $\Rmap{G}{L}$ holds because $gLg^{-1}$ is open for every $g\in G$. To see that $\Rmap{G}{L}(G)$ is closed, let $\{g_m\}$ be a sequence in $G$ such that $\Rmap{G}{L}(g_m)$ converges, to $x$ say, in $\RCom{G}{L}$. Then $g_mL = xL$ for all sufficiently large $m$, whence $g_m$ belongs to $xL$ for all large $m$ and $\{g_m\}$ is contained in a compact subset of $G$. Choose any accumulation point, $g$ say, of $\{g_m\}$. Then $x = \Rmap{G}{L}(g)$ and we have shown that $\Rmap{G}{L}(G)$ is closed. Hence $\RCom{G}{L} = \Rmap{G}{L}(G)$ and $\Rmap{G}{L}(L)$ is a compact open subgroup of $\RCom{G}{L}$. Therefore $\Rmap{G}{L}$ is an open map.
\end{demo}

\begin{lemma}
\label{lem:homtotop}
Let $\Gamma$ be a group with subgroup  $\Lambda$ and let $G$ be a topological group with a compact, open subgroup $L$. Suppose that there is a homomorphism $\varphi: \Gamma \to G$ such that
\begin{enumerate}
\item $\varphi(\Gamma)$ is dense in $G$, and 
\item $\varphi^{-1}(L)= \Lambda$. 
\end{enumerate} 
Then $\Lambda$ is commensurated by~$\Gamma$ and  $\RCom{\Gamma}{\Lambda}$ is isomorphic to $\RCom{G}{L}$ (which by Lemma \ref{lem:topcase} is 
isomorphic to $G$ when $L$ contains no non-trivial normal subgroups of $G$).  
\end{lemma}
\begin{demo} The conditions imply that $\varphi(\Lambda)$ is dense in $L$. Hence, for any $\gamma\in \Gamma$,   
$$
\psi : \lambda\left(\Lambda\cap \gamma\Lambda\gamma^{-1}\right) \mapsto \varphi(\lambda)\left(L\cap \varphi(\gamma)L\varphi(\gamma)^{-1}\right) 
$$
maps $\Lambda/\left(\Lambda\cap \gamma\Lambda\gamma^{-1}\right)$ onto $L/\left(L\cap \varphi(\gamma)L\varphi(\gamma)^{-1}\right)$. If $\lambda$ belongs to $\Gamma$ and $\varphi(\lambda)$ to $L\cap \varphi(\gamma)L\varphi(\gamma)^{-1}$, then $\lambda$ is in $\Lambda\cap \gamma\Lambda\gamma^{-1}$, whence $\psi$ is also injective. Therefore 
$$
[\Lambda : \Lambda\cap \gamma\Lambda \gamma^{-1}] = [L : L\cap \varphi(\gamma)L\varphi(\gamma)^{-1}]<\infty,
$$ 
and $\Gamma$ commensurates $\Lambda$.

Define $\hat\varphi : \Gamma/\Lambda \to G/L$ by $\hat\varphi(\gamma\Lambda) = \varphi(\gamma)L$. Then: $\hat\varphi$ is surjective because $\varphi$ has dense range and $L$ is open; $\hat\varphi$ is injective because $\varphi^{-1}(L)= \Lambda$; and $\hat\varphi$ intertwines the $\Gamma$- and $G$-actions, that is, 
$$
\hat\varphi(\gamma.x) = \varphi(\gamma).\hat\varphi(x),\quad (x\in \Gamma/\Lambda,\ \gamma\in \Gamma).
$$ 
Hence the closure of $\Rmap{\Gamma}{\Lambda}(\Gamma)$ in $\symm(\Gamma/\Lambda)$ is isomorphic to the closure of $\Rmap{G}{L}(\varphi(\Gamma))$ in $\symm(G/L)$ and, since $\varphi(\Gamma)$ is dense in $G$, the latter is equal to the closure of $\Rmap{G}{L}(G)$ in $\symm(G/L)$. By Lemma~\ref{lem:homtotop} therefore, $\RCom{\Gamma}{\Lambda}$ is isomorphic to $\RCom{G}{L}$. 
\end{demo}

\begin{corollary}
\label{cor:universality}
Let $\Lambda < \Gamma$ be a commensurated subgroup. Then the group
$\RCom{\Gamma}{ \Lambda}$ has the following universality property: 
Whenever $G$ is a totally disconnected locally compact group with $L<G$ a compact open subgroup
and $\varphi:\Gamma \to G$ is a homomorphism as in Lemma \ref{lem:homtotop},
there exists a continuous epimorphism $\psi: G \to \RCom{\Gamma}{ \Lambda}$
such that 
$\psi \circ \varphi :\Gamma \to G \to \RCom{\Gamma}{ \Lambda}$ is the
natural homomorphism. 
\end{corollary}

\begin{demo}
Take $\psi:G \to \RCom{G}{L} \cong \RCom{\Gamma}{ \Lambda}$ the natural 
homomorphism and use Lemmas~\ref{lem:homtotop} and \ref{lem:topcase}.
\end{demo}

The situation where there are nested subgroups $\Lambda\leq \Upsilon \leq \Gamma$ with $\Lambda$ and $\Upsilon$ commensurated by $\Gamma$ is analogous to that of nested subgroups $L\leq N \leq G$ normalized by $G$. In the case of normal subgroups, the Second Isomorphism Theorem gives an epimorphism $\rho : G/L \to G/N$ with $\ker(\rho) \cong N/L$. The following is the corresponding statement for relative profinite completions. 
\begin{lemma}
\label{lem:nestedcommensurated}
Let $\Lambda \leq \Upsilon$ be subgroups commensurated by $\Gamma$. Then: 
\begin{enumerate}
\item the coset inclusion map $\Gamma/\Lambda \to \Gamma/\Upsilon$ determines a continuous homomorphism 
$$
\rho : \RCom{\Gamma}{\Lambda} \to \RCom{\Gamma}{\Upsilon}
$$ 
that has dense range; \label{lem:2ndIso1}
\item the closure of $\rho(\Rmap{\Gamma}{\Lambda}(\Upsilon))$ is a compact open subgroup, call it $V$, of $\RCom{\Gamma}{\Upsilon}$ and\\ 
$\Rmap{\Gamma}{\Lambda}^{-1}\circ \rho^{-1}(V) = {\Upsilon}$ (more generally, it follows that the last equality can be replaced by commensurability
if $V$ is replaced by any compact open subgroup $V'$);  \label{lem:2ndIsoint}
\item the subgroup $\Rmap{\Gamma}{\Lambda}(\ker(\Rmap{\Gamma}{\Upsilon}))$ is contained in 
$\ker(\rho)$; 
and \label{lem:2ndIso2} 
\item  the restriction map $g \mapsto g|_{\Upsilon/\Lambda}$ is a homomorphism $\ker(\rho) \to \RCom{\Upsilon}{\Lambda}$. \label{lem:2ndIso3}
\end{enumerate}
\end{lemma}
\begin{demo}
\ref{lem:2ndIso1}. Define $\hat\rho : \Rmap{\Gamma}{\Lambda}(\Gamma) \to \Rmap{\Gamma}{\Upsilon}(\Gamma)$ by $\hat\rho(\Rmap{\Gamma}{\Lambda}(\gamma)) = \Rmap{\Gamma}{\Upsilon}(\gamma)$. Then $\hat\rho$ is well-defined and surjective because $\Rmap{\Gamma}{\Lambda}(\gamma)(\delta\Lambda) \mapsto \Rmap{\Gamma}{\Upsilon}(\gamma)(\delta\Upsilon)$ under the inclusion of $\Lambda$-cosets into $\Upsilon$-cosets. To see that $\hat\rho$ is continuous, note that 
$$
\hat\rho^{-1}(\Rmap{\Gamma}{\Upsilon}(\gamma\Upsilon\gamma^{-1})) = \Rmap{\Gamma}{\Lambda}(\gamma\Upsilon\gamma^{-1}),
$$
which is open in $\Rmap{\Gamma}{\Lambda}(\Gamma)$ because it is the union of $\Rmap{\Gamma}{\Upsilon}(\gamma\Lambda\gamma^{-1})$-cosets. Extend $\hat\rho$ by continuity to define $\rho$. Then the range of $\rho$ is dense in $\RCom{\Gamma}{\Upsilon}$ because it contains $\Rmap{\Gamma}{\Upsilon}(\Gamma)$. 

\ref{lem:2ndIsoint}. Since $\rho\circ\Rmap{\Gamma}{\Lambda} = \Rmap{\Gamma}{\Upsilon}$, the closure of $\rho\circ\Rmap{\Gamma}{\Lambda}(\Upsilon)$ is compact and open by Proposition~\ref{prop:profcomp}. If $\gamma\in \Gamma\setminus\Upsilon$, then $\gamma.\Upsilon \ne \Upsilon$. Hence $\Rmap{\Gamma}{\Upsilon}(\gamma)\not\in V$ in this case and the second statement follows.  

\ref{lem:2ndIso2}. and \ref{lem:2ndIso3}. follow immediately from the definition of $\rho$. 
\end{demo}

\begin{remarks}
\label{rem:remarks}
\noindent {\bf 1.} The range of $\rho$ need not be closed. Let $\Lambda$ be the trivial subgroup of $\Gamma$ and $\Upsilon$ be any commensurated subgroup. Then $\RCom{\Gamma}{\Lambda} \cong \Gamma$ and $\rho$ is essentially $\Rmap{\Gamma}{\Upsilon}$, which need not have closed range. 

\noindent {\bf 2.}  The restriction homomorphism $\ker(\rho) \to \RCom{\Upsilon}{\Lambda}$ need not be injective. Let $\Lambda$ and $\Gamma$ be as in Example~\ref{ex:Commensurate} and let $\Upsilon = A = \bigoplus_{k\in {\mathbb Z}} C_2$. Then $\RCom{\Gamma}{\Lambda} \cong \left(\prod_{k\in{\mathbb Z}} C_2\right) \rtimes {\mathbb Z}$, $\RCom{\Gamma}{\Upsilon} \cong {\mathbb Z}$ and $\RCom{\Upsilon}{\Lambda} \cong C_2$. Hence $\ker(\rho) \cong \prod_{k\in{\mathbb Z}} C_2$ and the restriction map sends $f\in  \prod_{k\in{\mathbb Z}} C_2$ to $f(0)\in C_2$. 

\noindent {\bf 3.} The homomorphism $\rho$ may be an isomorphism when $\Lambda\ne \Upsilon$ and the restriction homomorphism $\ker(\rho) \to \RCom{\Upsilon}{\Lambda}$ need not be surjective. Let $\Gamma$ be as in Example~\ref{ex:Commensurate} but now let $\Lambda = \left\{ f\in A \mid f(0) = f(1) = \bar0\right\}$ and $\Upsilon = \left\{ f\in A \mid f(0) = \bar0\right\}$. Then $\RCom{\Gamma}{\Lambda}$ and $\RCom{\Gamma}{\Upsilon}$ are both isomorphic to $\left(\prod_{k\in{\mathbb Z}} C_2\right) \rtimes {\mathbb Z}$ and $\rho$ has trivial kernel. However $\RCom{\Upsilon}{\Lambda}\cong C_2$. 
\end{remarks}

Lemma~\ref{lem:homtotop} may be used to determine $\RCom{\Gamma}{\Lambda}$ in some cases. 
\begin{examples}
\noindent {\bf 1.} The group $SL_n({\mathbb Z}[1/p])$ commensurates the subgroup $SL_n({\mathbb Z})$. Moreover, the natural homomorphism $\varphi : SL_n({\mathbb Z}[1/p]) \to SL_n({\mathbb Q}_p)$ has dense range and satisfies that  $\varphi^{-1}(SL_n({\mathbb Z}_p)) = SL_n({\mathbb Z})$. Hence, by Lemma~\ref{lem:homtotop},
$$
\RCom{SL_n({\mathbb Z}[1/p])}{SL_n({\mathbb Z})} \cong \RCom{SL_n({\mathbb Q}_p)}{SL_n({\mathbb Z}_p)}.
$$
Lemma~\ref{lem:topcase} shows that the latter is isomorphic to $SL_n({\mathbb Q}_p)$ modulo its largest normal subgroup contained in $SL_n({\mathbb Z}_p)$, namely, the center. Therefore
$$
\RCom{SL_n({\mathbb Z}[1/p])}{SL_n({\mathbb Z})} \cong PSL_n({\mathbb Q}_p).
$$ 

 \noindent {\bf 2.} For each positive integer $m$, the  \emph{congruence subgroup}, 
$$
\Lambda_m = \left\{ \left(r_{ij}\right)\in SL_n({\mathbb Z}) \mid r_{ij} \equiv \delta_{ij} \pmod{m}\right\},
$$ 
has finite index in $SL_n({\mathbb Z})$ 
 and hence is also commensurated by $SL_n({\mathbb Z}[1/p])$. When $p$ does not divide $m$, $\Lambda_m$ is the intersection of $\ker(\varphi_m)$ with $SL_n({\mathbb Z})$, where\\ 
 $\varphi_m : SL_n({\mathbb Z}[1/p]) \to SL_n({\mathbb Z}/m{\mathbb Z})$ is the homomorphism
$$
\varphi_m(\left(r_{ij}\right)) = \left( r_{ij} + m{\mathbb Z} \right).
$$

Consider
$$
\varphi\times \varphi_m : SL_n({\mathbb Z}[1/p]) \to SL_n({\mathbb Q}_p)\times SL_n({\mathbb Z}/m{\mathbb Z}).
$$
Then $\varphi\times \varphi_m$ has dense range because $\varphi_m$ is onto and $\varphi$ has dense range, and because $SL_n({\mathbb Z}/m{\mathbb Z})$ is finite and $SL_n({\mathbb Q}_p)$ has no finite index subgroup. Moreover, $L := SL_n({\mathbb Z}_p) \times\triv$ is a compact, open subgroup of $SL_n({\mathbb Q}_p)\times SL_n({\mathbb Z}/m{\mathbb Z})$ and $(\varphi\times\varphi_m)^{-1}(L) = \Lambda_m$.  Hence, by  Lemmas~\ref{lem:topcase} and~\ref{lem:homtotop}, 
$$
\RCom{SL_n({\mathbb Z}[1/p])}{\Lambda_m} \cong PSL_n({\mathbb Q}_p)\times SL_n({\mathbb Z}/m{\mathbb Z}).
$$
\end{examples}

The example is a special case of a calculation that is needed for the proof of Theorem~\ref{thm:mainresult} in Section~\ref{sec:Thm1.3Proof}. Let $K$ be a global field, $\mathcal O$ its ring of integers, and let $\bf G <\bf GL_n $ be an absolutely simple, simply connected algebraic group defined over $K$.
 For a subset $V^\infty \subseteq S \subseteq V$, denote the ring of 
$S$-integers in $K$ by $\mathcal O _S$ (see Section~\ref{sec:Proof1.4} below). 
Then ${\mathbf G}({\mathcal O})$ is a commensurated subgroup of ${\mathbf G}({\mathcal O}_S)$. We shall also use the following.

\smallskip
\noindent {\bf Notation.} Let $\adele$ be the ring of adeles of $K$. For a set of places $S \subseteq V$ we denote 
by ${\mathbf G}(\adele^{S})$ the subgroup of ${\mathbf G}(\adele)$ obtained by
taking a restricted direct product over places {\emph in} $S$ only. 
We denote by $S_f$ 
 the set of finite places in $S$, 
and ${\mathbf G}({\mathbb O}^{S_f})= \Pi_{\nu \in S_f} {\mathbf G}({\mathcal O_\nu}) <{\mathbf G}(\adele^{S_f}) $ denotes the 
natural compact
open subgroup of the latter.
 
\smallskip
Assume now further that $\bf G$ is  $K$-isotropic (or at least isotropic over one infinite place). Then by 
strong approximation the diagonal embedding 
$\varphi : {\mathbf G}({\mathcal O}_S)\to {\mathbf G}(\adele^{S_f})$ has dense range and $\varphi^{-1}({\mathbf G}({\mathbb O}^{S_f})) = {\mathbf G}({\mathcal O})$. Hence for any finite index subgroup
$\Gamma < {\mathbf G}({\mathcal O}_S)$, by Lemma~\ref{lem:homtotop} $\RCom{\Gamma}{(\Gamma\cap {\mathbf G}({\mathcal O}))}$ is isomorphic to 
$\RCom{ {\mathbf G}(\adele^{S_f})}{{\mathbf G}({\mathbb O}^{S_f})}$. As the 
center is the only normal subgroup of ${\mathbf G}(\adele^{S_f})$ which
 is contained in ${{\mathbf G}({\mathbb O}^{S_f})}$, Lemma~\ref{lem:topcase} shows that $\RCom{ {\mathbf G}(\adele^{S_f})}{{\mathbf G}({\mathbb O}^{S_f})}$ is isomorphic to ${\mathbf G}(\adele^{S_f})$ divided by its center. 
The following isomorphism has thus been established. 

\begin{proposition}
\label{prop:G(A)}
Let $K$, $S$ and  $\bf G $ be as in the beginning of the preceding paragraph. 
Suppose that $\Gamma$ is a finite index subgroup of ${\mathbf G}({\mathcal O}_S)$. 
Then $\RCom{\Gamma}{(\Gamma\cap {\mathbf G}({\mathcal O}))}$ is isomorphic to 
${\mathbf G}(\adele^{S_f})$ divided by its center.
\end{proposition}

Finally, the following general observation will be useful in Section~\ref{sec:adeles} below.

\begin{proposition}
\label{prop:chain}
Let $\Lambda <\Gamma$ be a commensurated subgroup of the discrete group 
$\Gamma$. Then there exists a transfinite increasing chain of 
commensurated subgroups, $\{\Lambda_{\alpha}\}$, beginning at $\Lambda$ and 
terminating at some $\Lambda '<\Gamma$, such that the following holds:\begin{enumerate}
\item $\Lambda_{\alpha}<\Lambda_{\alpha +1}$ has finite index for all $\alpha$,
\item $\Lambda _{\beta}=\bigcup _{\alpha<\beta}\Lambda _{\alpha}$ for limit
ordinals $\beta$, and 
\item  The group $\RCom{\Gamma}{\Lambda'}$ has
no non-trivial compact normal subgroup.
\end{enumerate}
\end{proposition}

\begin{demo} For notational convenience let us assume throughout the proof 
that $\Lambda$ does not contain a non-trivial normal subgroup of $\Gamma$ (which is anyway the case, up to a finite center, of interest to us here), so that we may identify $\Lambda$ and 
$\Gamma$ with their images in $G=\RCom{\Gamma}{\Lambda}$. Denote by $V<G$ the (compact open) closure of $\Lambda$ in $G$.
  Our purpose is to factor $G$ by all 
its compact normal subgroups. Had the group $M$ 
generated by all of
those been compact, it would readily follow that $G/M$ 
has no non-trivial compact normal subgroups (a property hereafter 
called  {\it reduced}). Correspondingly, 
the subgroup $V\cdot M$ would be compact, hence $V<M\cdot V$ of finite index, and
 the group $\Lambda '= \Gamma \cap (M\cdot V)$ would be a finite 
extension of $\Lambda$ for which $\RCom{\Gamma}{\Lambda'} \cong G/M$ (by Lemma 
\ref{lem:homtotop})
is reduced, thereby leaving one step in the chain of the Proposition (in the preceding examples
the above procedure applied to $\Lambda_m$ in {\bf 2} yields 
$SL_n(\ZZ)$ of {\bf 1}).
Unfortunately, it may generally happen that the above subgroup $M<G$ is not 
 compact, precompact, or even closed,
and that even dividing $G$ by its closure does not result in a reduced group.
This fact leads to the transfinite process that is described next and appears in the proposition.

\begin{definition}
\label{defn:Wang}
Let $G$ be a locally compact group. Define characteristic closed subgroups of $G$, $t^\alpha(G)$, $\alpha$ an ordinal, and $T(G)$ as follows. 
\begin{itemize}
\item $t(G) = \hbox{closure}\bigcup\left\{ N\triangleleft G \mid N \hbox{ a compact normal subgroup of }G \right\}$
\item $t^0(G) = \triv$ and, supposing that $t^\alpha(G)$ has been defined for some  ordinal $\alpha$ and that $q_\alpha : G\to G/t^\alpha(G)$ denotes the quotient map,  
$$
t^{\alpha+1}(G) = q_\alpha^{-1}(t(G/t^\alpha(G))).
$$ 
\item For $\alpha$ a limit ordinal, $t^\alpha(G)$ is the closure of $\bigcup\left\{ t^\beta(G) \mid \beta<\alpha \right\}$. 
\item $T(G) = \lim_\alpha t^\alpha(G)$. 
\end{itemize}
The subgroup $T(G)$ is the \emph{Wang radical} of $G$.
\end{definition}
It is clear from the definition that $t(G)$ is a closed characteristic subgroup of $G$ and it follows by induction that each $t^\alpha(G)$ is a closed characteristic subgroup. Then $\{t^\alpha(G)\}$ ($\alpha$ an ordinal) is an increasing family of closed subgroups of $G$ and it follows by a cardinality argument that $\lim_\alpha t^\alpha(G)$ exists. 

The Wang radical was defined by S.\,P.\,Wang in \cite{Wang}. As shown there (Theorem~1.5) it does indeed have the radical property that $T(G/T(G))$ is trivial. However, in order to be compatible with 
the formulation of the Proposition one refines the family $\{t^\alpha(G)\}$
by  
``unfolding'' (non-canonically) 
the first step in the Definition to an increasing chain of compact
normal (yet no longer 
characteristic) subgroups 
$$
t^{\alpha}(G)= t^{(\alpha ,0)}(G) < 
t^{(\alpha ,1)}(G) < \dots < t^{(\alpha ,i)}(G)< \dots 
< t^{(\alpha ,\eta_\alpha)}(G)=t^{\alpha+1}(G)
$$ 
obtained by adding compact normal subgroups of $G/t^{\alpha}(G)$ one at a time and taking the closure at limit ordinals. This defines the refined transfinite chain  
$\{t^{(\alpha ,i)}(G)\}$ denoted hereafter
by $ \bar t^{\beta}(G)$. 
Recall now that $V<G=\RCom{\Gamma}{\Lambda}$ denotes the closure of $\Lambda$. 
Define next the chain of subgroups $V^\beta:= V \cdot \bar t^ \beta(G)$.
Being open, it is easy to see that for limit ordinals $\beta$:
$\bigcup _{\delta <\beta}V^\delta =V^\beta$, that is, a closure operation
as in the first and third step of the Definition is not required. Finally,
set $\Lambda^\beta := \Gamma \cap V ^\beta$. Then by construction 
all $\Lambda^\beta$ are commensurated by~$\Gamma$ and, by Lemma \ref{lem:homtotop},
$\RCom{\Gamma}{\Lambda^\beta} \cong \RCom{G}{V ^\beta}$. By Lemma \ref{lem:topcase} the latter
is isomorphic to the group $G/ \bar t ^\beta(G)$, hence when the process
terminates we find, following the construction of the Wang radical,
a commensurated subgroup $\Lambda'<\Gamma$ with reduced $\RCom{\Gamma}{\Lambda} \cong
G/T(G)$.
\end{demo} 


\section{Flatness of Nilpotent and Polycyclic Groups}
\label{sec:flatness}

Flat groups of automorphisms are abelian modulo the uniscalar subgroup, as seen in Theorem~\ref{thm:flat}. It is shown in this section that the converse also holds: any group ${\mathcal H} \leq \Aut(G)$  that is finitely generated and abelian modulo a subgroup that stabilizes a compact open group is flat. This extends \cite[Theorem~5.5]{Wi04}, where it is shown that finitely generated abelian groups of automorphisms are flat. The condition that ${\mathcal H}$ be abelian modulo the stabilizer of a compact group can be weakened, and later results in this section treat finitely generated nilpotent and polycyclic groups of automorphisms. Further weakening is not possible, as a group described at the end of the section that is finitely generated and solvable, but not flat, shows. 


\subsection{The Abelian Case}
\label{sec:AbelianFlat}

Suppose that ${\mathcal H}\leq \Aut(G)$ has a normal subgroup, ${\mathcal N}$, such that:
\begin{description}
\item[\ref{sec:AbelianFlat}(a)] there is $V\in {\mathcal B}(G)$ such that $\alpha(V) = V$ for every $\alpha\in {\mathcal N}$, \phantom{AA} and
\item[\ref{sec:AbelianFlat}(b)]  ${\mathcal H}/{\mathcal N}$ is finitely generated and abelian.
\end{description}
It will be shown that ${\mathcal H}$ is flat by adapting the idea behind the notion of `local tidy subgroups' that is used in \cite[Definition~4.2]{Wi04} to take account of the subgroup ${\mathcal N}$. The following theorem is the adapted version of \cite[Theorem~5.5]{Wi04}, stated in a way that avoids the need for a new notion of  `local tidy subgroups'.
\begin{theorem}
\label{def:localtidy}
Let ${\mathcal N}\triangleleft {\mathcal H} \leq \Aut(G)$ and suppose that for every finite set ${\mathbf h}\subset {\mathcal H}$ there is $V_{\mathbf h}\in{\mathcal B}(G)$ such that $\beta(V_{\mathbf h})$ is tidy for $\langle \alpha, {\mathcal N}\rangle$ and stabilized by ${\mathcal N}$ for every $\alpha\in {\mathbf h}$ and $\beta\in {\mathcal H}$. 
Then ${\mathcal H}$ is flat.
\end{theorem}
\begin{demo}
The proof follows the same steps as that in \cite{Wi04}. A complication in the argument is that a group tidy for a finite set ${\mathbf h}$ is not necessarily tidy for $\langle{\mathbf h}\rangle$. However the hypothesized conditions on finite subsets ${\mathbf h}\subset {\mathcal H}$ allow $V_{\mathbf h}$ to be factored as a product of subgroups on which each element of ${\mathbf h}$ is either expanding or contracting, see \cite[Theorem~4.6]{Wi04}. Adding further elements to ${\mathbf h}$ gives a finer factoring. By choosing ${\mathbf h}$ sufficiently large, it may be shown that $V_{\mathbf h}$ is tidy for ${\mathcal H}$ and that the factoring of $V_{\mathbf h}$ is the one described in Theorem~\ref{thm:flat}(\ref{qnumber}), see \cite[Theorem~5.5]{Wi04}.
\end{demo}
 
To show that a group ${\mathcal H}$ satisfying \ref{sec:AbelianFlat}(a) and \ref{sec:AbelianFlat}(b) is flat, it suffices to show that these conditions imply that the hypothesis on finite subsets ${\mathbf h}$ in Theorem~\ref{def:localtidy} is satisfied. The first two lemmas deal with the case when ${\mathbf h}$ has one element.

\begin{lemma}
\label{lem:oneelement}
Let ${\mathcal N}\triangleleft {\mathcal H} \leq \Aut(G)$ be as in \ref{sec:AbelianFlat}(a) and \ref{sec:AbelianFlat}(b) and let $\alpha\in {\mathcal H}$. Then $\langle \alpha, {\mathcal N}\rangle$ is flat.
\end{lemma}
\begin{demo}
Let $V\in {\mathcal B}(G)$ be such that $\beta(V)=V$ for every $\beta\in {\mathcal N}$. It will be shown that when the tidying procedure for $\alpha$ is  applied to~$V$, the resulting group $V''$ is tidy for $\langle \alpha, {\mathcal N}\rangle$.

Normality of ${\mathcal N}$ in ${\mathcal H}$ implies that $\beta(\alpha^k(V)) = \alpha^k(V)$ for every $\beta\in {\mathcal N}$. Hence the subgroup $V' = \bigcap_{k=0}^n \alpha^k(V)$, defined in the first step of the tidying procedure, is stable under ${\mathcal N}$. If ${\mathcal N}$ has finite index in $\langle \alpha , {\mathcal N}\rangle$, then $\alpha(V') = V'$ and $V'$ is tidy for $\langle \alpha, {\mathcal N}\rangle$. Otherwise, continue the tidying procedure. 
Let $x$ belong to ${\mathcal K}_\alpha$, the group defined in Equation~(\ref{eq:define1K}), and let $\beta$ be in ${\mathcal N}$. Then $\alpha^k(\beta(x)) = \beta^{*}(\alpha^k(x))$, for some $\beta^{*}\in {\mathcal N}$. Hence $\alpha^k(\beta(x))$ belongs to $V'$ whenever $\alpha^k(x)$ does, which implies that ${\mathcal K}_\alpha$ and $K_\alpha^{(1)}$ are stable under ${\mathcal N}$. It follows that $V''$, defined in the third step, is stable under ${\mathcal N}$ in addition to being tidy for $\alpha$. 

To see that $V''$ is tidy for every element of $\langle \alpha , {\mathcal N}\rangle$, let $\alpha^k\beta$ be such an element. Then $\left(\alpha^k\beta\right)^n = \alpha^{kn}\beta^{*}$ for some $\beta^{*} \in {\mathcal N}$. Hence,  for each $k\in\ZZ$, $\bigcap_{n\geq0} \left(\alpha^k\beta\right)^n(V'') = \bigcap_{n\geq0} \alpha^{kn}(V'')$ and $\bigcap_{n\leq0} \left(\alpha^k\beta\right)^n(V'') = \bigcap_{n\leq0} \alpha^{kn}(V'')$. If $k>0$ the first set is $V''_+$ and the second $V''_-$, by (\ref{eq:tidyprop}), and {\it vice versa\/} if $k<0$. Since $V''$ is tidy for $\alpha$, it follows from Theorem~\ref{thm:tidiness} that $V''$ is also tidy for $\alpha^k\beta$ as required.
\end{demo}

\begin{lemma}
\label{lem:stably}
Let $\alpha\in {\mathcal H}$ and suppose that $V$ is tidy for $\langle \alpha , {\mathcal N}\rangle$ and that ${\mathcal N}$ stabilizes~$V$. Then $\beta(V)$ is tidy for $\langle \alpha , {\mathcal N}\rangle$ whenever $\beta\in {\mathcal H}$. 
\end{lemma}
\begin{demo}
Let $\alpha'\in \langle \alpha , {\mathcal N}\rangle$ and $\beta\in {\mathcal H}$. Then $\alpha'\beta = \beta\alpha'\gamma$ for some $\gamma\in {\mathcal N}$ and it follows that
$$
[\alpha'(\beta(V)) : \alpha'(\beta(V))\cap \beta(V)] = [\beta(\alpha'(\gamma(V))) : \beta(\alpha'(\gamma(V))) \cap \beta(V)].
$$ 
Since ${\mathcal N}$ stabilizes $V$, $\gamma(V)$ may be replaced by $V$ on the right hand side of this equation. Hence, since $\beta$ is an automorphism, the right hand side is equal to the scale of $\alpha'$ and it follows $\beta(V)$ is tidy for $\alpha'$.
\end{demo}

\begin{lemma} 
\label{lem:kalpha}
For each $\alpha\in {\mathcal H}$ set $K_{\langle\alpha, {\mathcal N}\rangle}:= \bigcap\left\{ V \mid V\hbox{ is tidy for } \langle \alpha, {\mathcal N}\rangle\right\}$. Then
\begin{enumerate}
\item \label{part:one} $\beta(K_{\langle\alpha, {\mathcal N}\rangle}) = K_{\langle\alpha, {\mathcal N}\rangle}$ for every $\beta\in {\mathcal H}$ and
\item \label{part:two} every $V'\in {\mathcal B}(G)$ that satisfies {\bf TA($\langle \alpha, {\mathcal N}\rangle$)} and such that $K_{\langle\alpha, {\mathcal N}\rangle} \subseteq V'$ is tidy for $\langle \alpha, {\mathcal N}\rangle$.
\end{enumerate}
\end{lemma}
\begin{demo}
(\ref{part:one}) Lemma~\ref{lem:stably} shows that $\left\{ V \mid V\hbox{ is tidy for } \langle \alpha, {\mathcal N}\rangle\right\}$ is unchanged if $\beta$ is applied to each element. Therefore the intersection of this set is stable under $\beta$.

(\ref{part:two}) If $V'$ satisfies {\bf TA($\alpha'$)} for some $\alpha'\in \langle \alpha, {\mathcal N}\rangle$, then $V'$ will be tidy for $\alpha'$ provided that it contains the subgroup $K_{\alpha'}^{(2)}$, defined in (\ref{eq:define2K}). Since $K_{\alpha'}^{(2)}$ is contained in every subgroup tidy for $\alpha'$, it is contained in $K_{\langle\alpha, {\mathcal N}\rangle}$. (This is the same argument as used to prove one direction of Lemma~3.31(3) in \cite{BaW}.)
\end{demo}

\begin{proposition}
\label{prop:localtidy}
Let ${\mathcal N}\triangleleft {\mathcal H} \leq \Aut(G)$ be as in \ref{sec:AbelianFlat}(a) and \ref{sec:AbelianFlat}(b). Then for every finite set ${\mathbf h}\subset {\mathcal H}$ there is $V_{\mathbf h}\in{\mathcal B}(G)$ such that $\beta(V_{\mathbf h})$ is tidy for $\langle \alpha, {\mathcal N}\rangle$ for every $\alpha\in {\mathbf h}$ and $\beta\in {\mathcal H}$. 
\end{proposition}
\begin{demo}
The proof is by induction on the number of elements in ${\mathbf h}$. Lemmas~\ref{lem:oneelement} and \ref{lem:stably} establish the case when ${\mathbf h}$ has one element. 

Assume that the claim has been established for $k$-element subsets of ${\mathcal H}$  and let ${\mathbf h} = \{\alpha_1,\dots,\alpha_k,\alpha_{k+1}\}$. Then there is $V\in {\mathcal B}(G)$ such that for every  $\beta\in {\mathcal H}$,  $\beta(V)$ is tidy for $\langle \alpha_j, {\mathcal F}\rangle$, $j\in\{1,\dots, k\}$. Applying the first step of the tidying procedure for $\alpha_{k+1}$ to $V$, there is an $n\in \NN$ such that $V' := \bigcap_{j=0}^n \alpha_{k+1}^j(V)$ satisfies {\bf TA($\alpha_{k+1}$)}. Then $\beta(V')$ is also tidy for $\langle \alpha_j, {\mathcal F}\rangle$ for each $j$ by Proposition~\ref{prop:recall}(\ref{eq:tidyintersect}). Lemma~\ref{lem:kalpha}(\ref{part:one}) shows that $K_{\langle\alpha_{k+1},{\mathcal N}\rangle}$ is $\beta$-stable for every $\beta\in {\mathcal H}$ and so, by \cite[Theorem~3.3]{Wi04}, $V''$ is tidy for $\langle \alpha_j, {\mathcal N}\rangle$ for each $j\in\{1,\dots,k\}$, where $V'' = K_{\langle\alpha_{k+1},{\mathcal N}\rangle}W$ and $W = \left\{x\in V' \mid xK_{\langle\alpha_{k+1},{\mathcal N}\rangle} \subseteq K_{\langle\alpha_{k+1},{\mathcal N}\rangle}V'\right\}$. By construction, $V''$ is also tidy for $\alpha_{k+1}$ and so is tidy for ${\langle\alpha_{k+1},{\mathcal N}\rangle}$ by the same argument as in the last paragraph of the proof of Lemma~\ref{lem:oneelement}.
\end{demo}

The proposition implies that ${\mathcal H}\leq \Aut(G)$ satisfying \ref{sec:AbelianFlat}(a) and 
\ref{sec:AbelianFlat}(b) also satisfies the hypotheses of Theorem~\ref{def:localtidy}. Therefore ${\mathcal H}$ is flat. 
\begin{theorem}
\label{thm:abelian}
Let ${\mathcal H}\leq \hbox{Aut}(G)$ and suppose that there is ${\mathcal N}\triangleleft {\mathcal H}$ that stabilizes some compact open $V\in {\mathcal B}(G)$ and such that ${\mathcal H}/{\mathcal N}$ is a finitely generated abelian group. Then ${\mathcal H}$ is flat. \endproof
\end{theorem}

\subsection{Nilpotent Groups are Flat}
\label{sec:NilpotentFlat}


Suppose that ${\mathcal H}\leq \Aut(G)$ has a normal subgroup, ${\mathcal N}$, such that:
\begin{description}
\item[\ref{sec:NilpotentFlat}(a)] there is $V\in {\mathcal B}(G)$ such that $\alpha(V) = V$ for every $\alpha\in {\mathcal N}$, \phantom{AA} and
\item[\ref{sec:NilpotentFlat}(b)]  ${\mathcal H}/{\mathcal N}$ is finitely generated and nilpotent.
\end{description}
It will be shown that ${\mathcal H}$ is flat by reducing to the case when ${\mathcal H}/{\mathcal N}$ is abelian. 

Denote by  ${\mathcal Z}_1$ the inverse image in ${\mathcal H}$ of the center of ${\mathcal H}/{\mathcal N}$ under the quotient map. Should ${\mathcal H}/{\mathcal N}$ be non-abelian, there are $\alpha,\,\beta\in {\mathcal H}$ such that $\gamma := [\alpha,\beta]$ belongs to ${\mathcal Z}_1\setminus {\mathcal N}$. The group $\langle \beta, \gamma, {\mathcal N} \rangle$ is abelian and finitely generated modulo ${\mathcal N}$ and so is flat, by Theorem~\ref{thm:abelian}. It may be supposed therefore that the group $V$ in \ref{sec:NilpotentFlat}(a) is tidy for $\langle \beta, \gamma, {\mathcal N} \rangle$. 

The next, easily verified, lemma facilitates the reduction. 
\begin{lemma}
\label{lem:Stabnormal}
The subgroup ${\mathcal S} := \left\{ \alpha\in {\mathcal Z}_1 \mid \alpha(V) = V \right\}$ is normal in ${\mathcal H}$. \endproof
\end{lemma}

The subgroup ${\mathcal S}$ normalizes $V$ and ${\mathcal H}/{\mathcal S}$ is nilpotent and finitely generated because it is isomorphic to $({\mathcal H}/{\mathcal N})/({\mathcal S}/{\mathcal N})$. Therefore, if ${\mathcal S} \ne {\mathcal N}$, we may replace ${\mathcal N}$ in \ref{sec:NilpotentFlat}(a) and \ref{sec:NilpotentFlat}(b) by the larger group ${\mathcal S}$. In doing so, it might happen that $\gamma$ now belongs to ${\mathcal N}$. Should that be the case and ${\mathcal H}/{\mathcal N}$ still be non-abelian, choose new $\alpha$, $\beta$ and $\gamma$. If, redefining ${\mathcal Z}_1$ and ${\mathcal S}$, the new ${\mathcal S}$ is not equal to the new ${\mathcal N}$, then ${\mathcal N}$ may be replaced by a still larger subgroup. Since ascending chains of subgroups of the nilpotent group ${\mathcal H}/{\mathcal N}$ are finite, this process will terminate after a finite number of iterations, at which point  either:\\
(i) ${\mathcal H}/{\mathcal N}$ is abelian; \qquad or \\
(ii) there are $\alpha,\,\beta\in {\mathcal H}$ such that $\gamma := [\alpha,\beta] \in {\mathcal Z}_1\setminus {\mathcal N}$, $V$ is tidy for $\langle \beta, \gamma, {\mathcal N} \rangle$ and ${\mathcal S} = {\mathcal N}$. 
\\
The next lemma will show that (ii) leads to a contradiction. 

\begin{lemma}
\label{lem:centreuniscalar}
Let ${\mathcal H}$, ${\mathcal N}$ and ${\mathcal Z}_1$ be as above. Let $\alpha,\,\beta\in {\mathcal H}$ and suppose that $\gamma := [\alpha,\beta]$ belongs to ${\mathcal Z}_1$.  Then $s(\gamma) = 1$.
\end{lemma}

\begin{demo}
We begin by showing that $s(\beta\gamma^n)$ does not depend on $n$. Since the scale is constant on conjugacy classes, see Theorem~\ref{thm:scale}(\ref{conjclass}), $s(\beta\gamma^n) = s(\alpha\beta\gamma^n\alpha^{-1})$. On the other hand, $\alpha\beta\gamma^n\alpha^{-1} = \beta\gamma^{n+1} \gamma'$, where $\gamma'\in {\mathcal N}$. Hence 
\begin{align*} 
s(\beta\gamma^n) = {}& [\beta\gamma^{n+1} \gamma'(V) : \beta\gamma^{n+1} \gamma'(V) \cap V] \\
 = {}& [\beta\gamma^{n+1}(V) : \beta\gamma^{n+1}(V) \cap V] = s(\beta\gamma^{n+1})
\end{align*}
and $s(\beta\gamma^n)$ does not depend on $n$.

Since $V$ is tidy for $\langle \beta, \gamma, {\mathcal N} \rangle$, there are subgroups $V_j\leq V$ as in Theorem~\ref{thm:flat}(\ref{module}) such that  $V = V_0V_1\cdots V_q$ and $s(\gamma) = \prod_{\rho_j(\gamma)>0} \Delta_j(\gamma)$. Hence, if $s(\gamma) \ne 1$, there is $j\in\{1,2,\dots,q\}$ such that $\Delta_j(\gamma)>1$. Then $\Delta_j(\beta\gamma^n) = \Delta_j(\beta)\Delta_j(\gamma)^n \to\infty$ as $n\to\infty$. Since $s(\beta\gamma^n) = \prod_{\rho_j(\beta\gamma^n)>0} \Delta_j(\beta\gamma^n)$, it follows that $s(\beta\gamma^n)\to \infty$ as $n\to\infty$, in contradiction to $s(\beta\gamma^n)$  not depending on~$n$. 
\end{demo}

Lemma~\ref{lem:centreuniscalar} implies that each commutator $\gamma\in {\mathcal Z}_1$ has scale~1. Since $\gamma^{-1}$ is also a commutator, $s(\gamma^{-1}) = 1$ as well and so $\gamma(V) = V$. Hence $\gamma\in {\mathcal S}$ and is not possible that (ii) holds. Therefore ${\mathcal H}/{\mathcal N}$ is abelian and ${\mathcal H}$ satisfies the hypotheses of Theorem~\ref{thm:abelian}, yielding the main result of this subsection. 

\begin{theorem}
\label{thm:nilpotent}
Let ${\mathcal N}\triangleleft {\mathcal H}\leq \hbox{Aut}(G)$ and suppose that ${\mathcal N}$ stabilizes some $V$ in ${\mathcal B}(G)$ and that ${\mathcal H}/{\mathcal N}$ is a finitely generated nilpotent group. Then ${\mathcal H}$ is flat. 
\endproof
\end{theorem}

A topological version of this assertion may be made for inner automorphisms. 

\begin{theorem}
\label{thm:compactly}
Let $H$ be a closed nilpotent subgroup of $G$ that is topologically generated by a compact set. Then $H$ is flat.
\end{theorem}

\begin{demo}
Since compactly generated, locally compact nilpotent groups are pro-Lie, \cite{HLM}, there is a compact, open, normal subgroup $K\triangleleft {H}$. The quotient group ${H}/K$ is finitely generated, let $\{  x_1K, \dots, x_nK\}$ be a generating set.  Then   $\langle x_1, \dots, x_n \rangle$ is a finitely generated, nilpotent subgroup of $G$ and is therefore flat by Theorem~\ref{thm:nilpotent}. Choose $V$ tidy for $\langle x_1, \dots, x_n \rangle$. Then $VK$ is a compact open subgroup of $H$. Since $K$ is normal, $H = \langle x_1, \dots, x_n\rangle K$ and for each $xk\in H$ we have $xkVK(xk)^{-1} = xVx^{-1}K$. Hence $(VK)_\pm = V_\pm K$ and it may be verified that $VK$ is tidy for $xk$. 
\end{demo}

\subsection{Polycyclic Groups are Virtually Flat}
\label{sec:polycyclic}

Suppose that ${\mathcal H}\leq \Aut(G)$ has a normal subgroup, ${\mathcal N}$, such that:
\begin{description}
\item[\ref{sec:polycyclic}(a)] there is $V\in {\mathcal B}(G)$ such that $\alpha(V) = V$ for every $\alpha\in {\mathcal N}$, \phantom{AA} and
\item[\ref{sec:polycyclic}(b)]  ${\mathcal H}/{\mathcal N}$ is polycyclic.
\end{description}
It will be shown that ${\mathcal H}$ has a finite index subgroup that is flat. 

Since ${\mathcal H}/{\mathcal N}$ is polycyclic, there is a series
\begin{equation}
\label{eq:series}
{\mathcal N} = {\mathcal H}_0 \triangleleft {\mathcal H}_1 \triangleleft \cdots \triangleleft {\mathcal H}_l = {\mathcal H}
\end{equation} such that ${\mathcal H}_{j}/{\mathcal H}_{j-1}$ is a finitely generated abelian group for each $j\in \{1,2,\dots,l\}$. The proof is by induction on the length, $l$, of the series. 

When $l=1$, we have that ${\mathcal H}/{\mathcal N}$ is finitely generated and abelian and ${\mathcal H}$ is consequently flat by Theorem~\ref{thm:abelian}. Assume that it has been established that any group having a series (\ref{eq:series}) of length $l$ has a finite index subgroup that is flat and suppose that ${\mathcal H}$ has a series of length $l+1$,
\begin{equation}
\label{eq:seriesplus}
{\mathcal N} = {\mathcal H}_0 \triangleleft {\mathcal H}_1 \triangleleft \cdots \triangleleft {\mathcal H}_l \triangleleft {\mathcal H}_{l+1} = {\mathcal H}.
\end{equation}
Then, by the induction hypothesis, there is a finite index subgroup, $\hat{\mathcal H}_l$, of ${\mathcal H}_l$ that is flat. 

Let $N_{\mathcal H}(\hat{\mathcal H}_l) := \left\{ \alpha\in {\mathcal H} \mid \alpha\hat{\mathcal H}_l\alpha^{-1} = \hat{\mathcal H}_l\right\}$ be the normalizer in ${\mathcal H}$ of $\hat{\mathcal H}_l$. Since every subquotient of ${\mathcal H}/{\mathcal N}$ is polycyclic with rank at most $l+1$, the next lemma, applied with $C = {\mathcal H}_l$ and $B = \hat{\mathcal H}_l$, will imply that $\bigcap_{\alpha\in{\mathcal H}} \alpha\hat{\mathcal H}_l\alpha^{-1}$ has finite index in ${\mathcal H}_l$, whence it will follow that $N_{\mathcal H}(\hat{\mathcal H}_l)$ has finite index in ${\mathcal H}$. 
\begin{lemma}
\label{lem:infinite}
Let $C$ be a group and $B$ be a finite index subgroup of $C$. Let ${\mathcal A} = \Aut(C)$ and suppose that $\bigcap_{\alpha\in {\mathcal A}} \alpha(B)$ has infinite index in $C$. Then there is a finite index subgroup $E\leq C$ and a finite (possibly abelian) simple group $F$ such that for every $n\in \NN$ there is a surjective homomorphism $E \rightarrow F^n$.
\end{lemma}

\begin{demo}
The kernel of the representation of $C$ on $C/B$ is a finite index subgroup that is normal in $C$. Replacing $B$ by this kernel if necessary, it may be assumed that $B$ is a normal subgroup. Consider a composition series
$$
B = B_0 \triangleleft B_1 \triangleleft \cdots B_{r-1} \triangleleft B_r = C
$$ 
for $C/B$, where $B_{j+1}/B_j$ is simple for each $j \in \{0,1,\dots, r-1\}$. Choose $B_j$ to be the largest subgroup in the series such that $\bigcap_{\alpha\in {\mathcal A}} \alpha(B_j)$ has infinite index in $C$ and put $E = \bigcap_{\alpha\in {\mathcal A}} \alpha(B_{j+1})$. Then $E$ has finite index in $C$ and is ${\mathcal A}$-invariant. Since $\bigcap_{\alpha\in {\mathcal A}} \alpha(B_{j}\cap E)$ has infinite index in $C$,  $B_{j}\cap E$ is a proper subgroup of $E$ and it is a normal subgroup because $E\leq B_{j+1}$. The index of $E\cap B_j$ in $E$ is bounded by that of $B$ in $C$ and is strictly less if $r>1$. Hence, if $E/(E\cap B_j)$ is not simple, the argument may be repeated with $E$ and $E\cap B_j$ in place of $C$ and $B$ until a simple quotient is obtained. 

Assuming now that $E/B$ is simple, put $F = E/B$. Then for any finite $A\subset {\mathcal A}$, the subgroup $\bigcap_{\alpha\in {A}} \alpha(B)$ has finite index in $E$ and $E/\left(\bigcap_{\alpha\in {A}} \alpha(B)\right)$ is isomorphic to $F^n$ for some $n$. Since $\bigcap_{\alpha\in {\mathcal A}} \alpha(B)$ has infinite index in $E$, every exponent $n$ occurs.
\end{demo}

Lemma~\ref{lem:infinite}, the induction hypothesis and (\ref{eq:seriesplus}) between them imply that 
$$
\hat{\mathcal H}_l \triangleleft {\mathcal H}_l\cap N_{\mathcal H}(\hat{\mathcal H}_l) \triangleleft N_{\mathcal H}(\hat{\mathcal H}_l),
$$ 
where $\hat{\mathcal H}_l$ has finite index in ${\mathcal H}_l\cap N_{\mathcal H}(\hat{\mathcal H}_l)$, $N_{\mathcal H}(\hat{\mathcal H}_l)$ has finite index in ${\mathcal H}$ and $N_{\mathcal H}(\hat{\mathcal H}_l)$ is abelian modulo ${\mathcal H}_l\cap N_{\mathcal H}(\hat{\mathcal H}_l)$. Hence $({\mathcal H}_l\cap N_{\mathcal H}(\hat{\mathcal H}_l))/\hat{\mathcal H}_l$ is a finite normal subgroup of $N_{\mathcal H}(\hat{\mathcal H}_l)/\hat{\mathcal H}_l$ and the centralizer of $({\mathcal H}_l\cap N_{\mathcal H}(\hat{\mathcal H}_l))/\hat{\mathcal H}_l$ consequently has finite index in $N_{\mathcal H}(\hat{\mathcal H}_l)/\hat{\mathcal H}_l$. Denote the inverse image of this centralizer under the quotient map by ${\mathcal C}$, so that ${\mathcal C}$ is a finite index subgroup of~$\mathcal H$ containing $\hat{\mathcal H}_l$ and ${\mathcal C}/\hat{\mathcal H}_l$ is abelian. The sought after flat group is a finite finite index subgroup of ${\mathcal C}$ that  will be found with the aid of the next lemma. 

\begin{lemma}
\label{lem:finiteorbits}
Let ${\mathcal J}\leq \Aut(G)$ be flat and finitely generated and suppose that ${\mathcal J} \triangleleft {\mathcal K}$ for some subgroup ${\mathcal K}$ of $\Aut(G)$. Then ${\mathcal J}(1)$ is normalized by ${\mathcal K}$ and the induced action of ${\mathcal K}$ on ${\mathcal J}/{\mathcal J}(1)$ has finite orbits.
\end{lemma}

\begin{demo}
The identity $[(\beta\alpha\beta^{-1})\beta(V) : (\beta\alpha\beta^{-1})\beta(V) \cap \beta(V)] = [\alpha(V) : \alpha(V)\cap V]$, holding for any automorphisms $\alpha$ and $\beta$ and compact open subgroup $V$, implies that $V$ is minimizing for $\alpha$ if and only if $\beta(V)$ is minimizing for $\beta\alpha\beta^{-1}$ and that $s(\beta\alpha\beta^{-1}) = s(\alpha)$. Hence, if $\beta \in {\mathcal K}$ and $V$ is tidy for ${\mathcal J}$, then $\beta(V)$ is tidy for $\beta{\mathcal J}\beta^{-1} = {\mathcal J}$ and ${\mathcal J}(1)$ is normalized by ${\mathcal K}$.

Recall from Corollary~\ref{cor:submultiplicative} that the function $\gamma{\mathcal J}(1) \mapsto \log s(\gamma) + \log s(\gamma^{-1})$ is a norm on ${\mathcal J}/{\mathcal J}(1)$. Hence the set $\left\{ \gamma{\mathcal J}(1)\in {\mathcal J}/{\mathcal J}(1) \mid s(\gamma) \leq M\right\}$ is finite for each constant~$M$. Since, as shown in the previous paragraph, the scale is constant on ${\mathcal K}$-orbits in ${\mathcal J}/{\mathcal J}(1)$, these orbits must be finite.
\end{demo}

Let $\alpha_1\hat{\mathcal H}_l(1)$, \dots, $\alpha_m\hat{\mathcal H}_l(1)$ be a generating set for $\hat{\mathcal H}_l/\hat{\mathcal H}_l(1)$. Then Lemma~\ref{lem:finiteorbits} shows that $\left\{\beta\alpha_i\beta^{-1}\hat{\mathcal H}_l(1) \mid \beta \in N(\hat{\mathcal H}_l)\right\}$ is finite for each $i\in \{1,2,\dots,m\}$. Therefore 
$$
\hat{\mathcal H} := \left\{ \beta\in {\mathcal C} \mid \beta\alpha\beta^{-1}\hat{\mathcal H}_l(1) = \alpha\hat{\mathcal H}_l(1)\right\},
$$
the centralizer in ${\mathcal C}$ of $\hat{\mathcal H_l}$ modulo $\hat{\mathcal H}_l(1)$, has finite index in ${\mathcal C}$ and therefore also in ${\mathcal H}$. Since $\hat{\mathcal H}$ is finitely generated and is 2-step nilpotent modulo the group $\hat{\mathcal H}_l(1)$, Theorem~\ref{thm:nilpotent} shows that $\hat{\mathcal H}$ is flat.

\begin{theorem}
\label{thm:polycyclic}
Let ${\mathcal N}\triangleleft {\mathcal H}\leq \Aut(G)$. Suppose that ${\mathcal N}$ stabilizes some $V$ in ${\mathcal B}(G)$ and that ${\mathcal H}/{\mathcal N}$ is a polycyclic group. Then ${\mathcal H}$ has a finite index subgroup that is flat. 
\endproof
\end{theorem}

\medskip

\noindent{\bf Remark.} Theorem~\ref{thm:polycyclic} cannot be extended to cover solvable groups of automorphisms that are not polycyclic, even if they are finitely generated. The argument breaks down because subgroups of finitely generated solvable groups need not be finitely generated (in fact, the latter property characterizes polycyclic groups among the solvable ones). One example, where $\Gamma$ is the so-called
lamplighter group, is as follows.

Consider the following linear realization of $\Gamma$. 
Let $\field_2$ be the field of two elements and $\Gamma = \ZZ {\ltimes} (\field_2[t,t^{-1}])$, where
$\field_2[t,t^{-1}]$ is the ring of polynomials in $t,t^{-1}$ over $\field_2$, viewed as an abelian group, and $\ZZ$ acts 
through
multiplication by powers of $t$. Then $\Gamma$ embeds naturally 
in the algebraic group $G=SL_2(\field_2((t)))$
(embed $\ZZ$ via powers of the matrix $ A=
\left(\begin{array}{cc} t &1 \\
1 &t^{-1}\end{array}\right)$, and $\field_2[t,t^{-1}]$ in the elementary unipotent 
subgroup).
It is not virtually flat in $G$ as the commutator subgroup 
of any finite 
index subgroup is an unbounded (unipotent) subgroup, while the normalizer of any compact open subgroup of $G$ 
is compact (any open, proper subgroup of 
$SL_n(\field_2((t)))$ is compact). 
In fact, using this idea, $p$-adic  specializations, and the fact that a finitely  generated {\it discrete}  solvable subgroup of $GL_n(\CC$) is necessarily polycyclic, it seems very likely that a finitely generated {\it linear} solvable group is virtually flat (if and) only if it is polycyclic.


\section{The Inner and Outer Commensurator-Normalizer Properties}
\label{sec:Comm_norm}
 
 

The following basic result is also a key tool in the proof of Theorem~\ref{thm:mainresult}.

\begin{proposition}
\label{prop:comm_norm}
If $\Gamma_1$ has the outer commensurator-normalizer property and $\Gamma_2$ is commensurable with $\Gamma_1$ then so does $\Gamma_2$.
\end{proposition}

\begin{demo}
Let $\varphi : \Gamma_2\to \Delta$ be a homomorphism and suppose that $\Lambda\leq \Delta$ is commensurated by $\varphi(\Gamma_2)$. 
Consider first the case when $\Gamma_1\leq \Gamma_2$. Since $\varphi(\Gamma_2)$ commensurates $\Lambda$, so does $\varphi(\Gamma_1)$, and there is $\Lambda'$ commensurable with $\Lambda$ that is normalized by $\varphi(\Gamma_1)$. Then the subgroup $\Lambda'' := \bigcap_{g\in \Gamma_2} \varphi(g) \Lambda' \varphi(g)^{-1}$ is normalized by $\varphi(\Gamma_2)$ and  has finite index in $\Lambda'$, whence it is commensurable with $\Lambda$. 

The proof will be completed by treating the case when $\Gamma_2\leq \Gamma_1$. To this end, define a group $\widetilde{\Delta}$ and homomorphism $\widetilde{\varphi} : \Gamma_1 \to \widetilde{\Delta}$ as follows. Choose a transversal, $X$, for $\Gamma_1/\Gamma_2$ and let $\widetilde{\Delta}$ be the wreath product $S(X) \ltimes \Delta^X$, where the permutation group $S(X)$ acts by composition, that is:
$$
S(X)\times \Delta^X \to \Delta^X : (\sigma,f) \mapsto f\circ \sigma.
$$ 
For $g\in \Gamma_1$, let $\sigma_g \in S(X)$ be the permutation satisfying
$$
gx\Gamma_2 = \sigma_g(x)\Gamma_2, \quad (g\in \Gamma_1,\ x\in X),
$$
so that $g\mapsto \sigma_g: \Gamma_1\to S(X)$ is a homomorphism, and let $\alpha : \Gamma_1\times X \to \Gamma_2$ be the cocycle making the following equality in $\Gamma _1$ hold:
$$
gx =  \sigma_g(x)\alpha(g,x), \quad (g\in\Gamma_1,\,x\in X).
$$
For $g\in \Gamma_1$ let $f_g\in \Delta^X$ be the function $f_g(x) = \varphi(\alpha(g,x))$ and define
$$
\widetilde{\varphi} : \Gamma_1 \to \widetilde{\Delta}:= S(X)\ltimes \Delta^X\ \text{ by }\ \widetilde{\varphi}(g) =(\sigma_g, f_g).
$$ 
Then the cocycle identity for $\alpha$ implies that in the second coordinate one has:
$$
f_{g_1g_2}(x) = \varphi(\alpha(g_1g_2,x)) = \varphi(\alpha(g_1,\sigma_{g_2}(x)))\varphi(\alpha(g_2,x)) = (f_{g_1}\circ\sigma_{g_2})(x)f_{g_2}(x),
$$ 
whence $\widetilde{\varphi}$ is a homomorphism.

It will now be necessary to distinguish between the identity element in $S(X)$, which is denoted by $\iota$, and that in $\Delta^X$, denoted by $e$.  Recall that $\Lambda < \Delta$ 
is the assumed subgroup commensurated by  $\varphi (\Gamma _2)$. The subgroup $\widetilde{\Lambda} = \left\{ (\iota, f)\in  \widetilde{\Delta} \mid f(x)\in \Lambda \text{ for all }x\in X\right\}$ of $ \widetilde{\Delta}$ is commensurated by $\widetilde{\varphi}(\Gamma_1)$. To see this, note that $\widetilde{\varphi}(g)^{-1} =  (\iota, f_g^{-1})(\sigma_{g^{-1}},e)$, hence 
$$
\widetilde{\varphi}(g)\widetilde{\Lambda}\widetilde{\varphi}(g)^{-1} = (\sigma_g,e)\left((\iota, f_g)\widetilde{\Lambda}(\iota, f_g^{-1})\right)(\sigma_{g^{-1}},e).
$$
Since $\Lambda$ is commensurated by $\varphi(\Gamma_2)$, it follows that  $(\iota, f_g)\widetilde{\Lambda}(\iota, f_g^{-1})$ is commensurable with $\widetilde{\Lambda}$. The conjugate by $(\sigma_g,e)$ is still commensurable with $\widetilde{\Lambda}$ because such conjugation leaves $\widetilde{\Lambda}$ invariant.

Since $\widetilde{\Lambda}$ is commensurated by $\widetilde{\varphi}(\Gamma_1)$, there is $\widetilde{\Lambda}' \leq \widetilde{\Delta}$ that is commensurable with $\widetilde{\Lambda}$ and normalized by $\widetilde{\varphi}(\Gamma_1)$. Since $\Delta^X$ is a normal subgroup with finite index in $\widetilde{\Delta}$, it may be supposed that $\widetilde{\Lambda}'\leq \Delta^X$. Assume, as we may, that $x=1$ is the coset representative chosen from the coset~$\Gamma_2$. Then 
$$
\alpha(g,1) = g\text{ and } \sigma_g(1)=1 \text{ for every }g\in \Gamma_2
$$ 
and the projection $\pi : \Delta^X\twoheadrightarrow \Delta$ defined by $\pi(f) := f(1)$ is a homomorphism satisfying $\pi\circ \text{ad}(\widetilde{\varphi}(g)) = \text{ad}(\varphi(g))$, where $\text{ad}(g)$ denotes conjugation by $g$. Therefore $\Lambda' := \pi(\widetilde{\Lambda}')$ is commensurable with $\Lambda$ and normalized by $\varphi(\Gamma_2)$. 
\end{demo}

\medskip
In the other direction, it is useful to keep in mind the following. 
\begin{theorem} 
\label{thm:counter}
The inner  commensurator-normalizer property is not generally invariant 
under passing to finite index subgroups.
\end{theorem}

In the rest of this section we construct the example accounting for
Theorem~\ref{thm:counter}.
Let $\field_2$ be the field of order~$2$. Then the additive groups  $\bigoplus_{l\in {\mathbb Z}} \field_2$ and $\prod_{l\in {\mathbb Z}} \field_2$ are vector spaces over $\field_2$ and are dual to each other through the pairing 
\begin{equation}
\label{eq:duality}
< h,k > = \sum_{l\in{\mathbb Z}} h(-l)k(l), \ \ \left({\textstyle h\in \bigoplus_{l\in {\mathbb Z}} \field_2,\, k\in \prod_{l\in {\mathbb Z}} \field_2}\right).
\end{equation}
Let $\tau$ and $\sigma$ be respectively the automorphisms of $\prod_{l\in {\mathbb Z}} \field_2$ defined by
$$
(\tau k)(l) = k(l+1) \text{ and } (\sigma k)(l) = k(-l), \ \ \left({\textstyle k\in \prod_{l\in {\mathbb Z}} \field_2}\right).
$$
Then $\tau$ and $\sigma$ restrict to be automorphisms of $\bigoplus_{l\in {\mathbb Z}} \field_2$ as well, which will be denoted the same way, and satisfy  the identities
$< \tau h,k>\, =\, < h,\tau k>$ and $<\sigma h,k>\, =\, < h,\sigma k>$. 

The groups $\langle \tau, \sigma\rangle$ and $\langle \tau \rangle$ are isomorphic to the infinite dihedral group and to ${\mathbb Z}$ respectively. Let $\Gamma = \langle \tau , \sigma \rangle \ltimes \bigoplus_{l\in {\mathbb Z}} \field_2$ and  $\Lambda = \langle \tau  \rangle \ltimes \bigoplus_{l\in {\mathbb Z}} \field_2$. Then $\Lambda$  is a subgroup of $\Gamma$ of index~2 and so these two groups are commensurable. 

The group $\Lambda$ does not have the inner commensurator-normalizer property because the subgroup 
$$
\Upsilon := \left\{ (\tau^0,(z_l)) \mid (z_l)\in {\textstyle \bigoplus_{l\in {\mathbb Z}} \field_2} ,\,z_l = 0 \text{ for }l<0\right\}
$$ 
is commensurated by $\Lambda$ and is not commensurable with any normal subgroup of $\Lambda$. However, $\Gamma$ does have the inner commensurator-normalizer property. 
\begin{proposition}
\label{prop:comm-norm}
Let $\Xi$ be a subgroup of $\Gamma$ that is commensurated by $\Gamma$. Then $\Xi$ is commensurable with a normal subgroup of $\Gamma$. 
\end{proposition} 

\begin{demo}  The special case when $\Xi$ is a subgroup of $\bigoplus_{l\in {\mathbb Z}} \field_2$ follows from part~\ref{prop:RC1} of the next lemma. The general case is proved subsequently. 

\begin{lemma}
\label{prop:reduced_case}
\begin{enumerate}
\item Let $\Xi$ be a subgroup of $\bigoplus_{l\in {\mathbb Z}} \field_2$ that is commensurated by $\langle \tau,\sigma\rangle$. Then $\Xi$ is either finite or has finite index in $\bigoplus_{l\in {\mathbb Z}} \field_2$.
\label{prop:RC1}
\item Let $K$ be a closed subgroup of $\prod_{l\in {\mathbb Z}} \field_2$ that is commensurated by $\langle \tau,\sigma\rangle$. Then $K$ is either finite or has finite index in $\prod_{l\in {\mathbb Z}} \field_2$. \label{prop:RC2}
\end{enumerate}
\end{lemma}
\begin{demo}
For $\Xi\leq \bigoplus_{l\in {\mathbb Z}} \field_2$, define 
$$
 \Xi^\perp = \left\{ k\in \prod_{l\in {\mathbb Z}} \field_2 \mid < h,k>\, = 0 \text{ for every }h\in \Xi\right\}.
$$
Then $ \Xi^\perp$ is a closed subgroup of $\prod_{l\in {\mathbb Z}} \field_2$ and 
$\Xi$ is finite (respectively, has finite index) if and only if $\Xi^\perp$ has finite index (resp., is finite).   Similarly, $\Xi$ is commensurated or normalized by $\langle \tau, \sigma \rangle$ if and only if $\Xi^\perp$ is. Hence, once part~\ref{prop:RC2} is proved,  part~\ref{prop:RC1} may be deduced by setting $K = \Xi^\perp$. 

The next two lemmas are needed for the proof of  part~\ref{prop:RC2}. 
\begin{lemma}
\label{lem:stable_finite}
Let $K$ be a closed subgroup of $\prod_{l\in {\mathbb Z}} \field_2$ such that $\tau(K) = K$. Then $K$ is either equal to $\prod_{l\in {\mathbb Z}} \field_2$ or is finite.
\end{lemma}
\begin{demo}
If $K$ is not equal to $\prod_{l\in {\mathbb Z}} \field_2$, then there is a nonzero $h\in \bigoplus_{l\in {\mathbb Z}} \field_2$ that annihilates~$K$. Since $K$ is stable under $\tau$, $< \tau^jh , K >\, = \{0\}$ for every $j\in {\mathbb Z}$. In other words, every element of $K$ satisfies the difference equations 
$$
< \tau^jh,k> \, = \sum_{l\in{\mathbb Z}} h(-l+j)k(l) = 0\ \ \ \  (j\in {\mathbb Z}).
$$
Therefore $K$ is finite. 
\end{demo}

\begin{lemma}
\label{lem:Kexpand}
Let $K$ be a closed subgroup of $\prod_{l\in {\mathbb Z}} \field_2$ such that $\tau(K)\gvertneqq K$. Then there is an integer $J$ such that $\prod_{j\geq J} \field_2$ is contained in $K$.
\end{lemma}
\begin{demo}
Since $\tau(K)$ is strictly larger than $K$, $K$ is a proper subgroup of $\prod_{l\in {\mathbb Z}} \field_2$ and it follows that there is a non-zero $h\in \bigoplus_{l\in {\mathbb Z}} \field_2$ such that  
\begin{equation}
\label{eqn:Kexpand}
< h,k>\, = 0\text{ for every }k\in K.
\end{equation} 
Given any such $h$, every $k\in K$ satisfies the difference equations
\begin{equation}
\label{eq:iterated}
< \tau^{-j}h,k> = \sum_{l\in{\mathbb Z}} h(-l-j)k(l) = 0, \text{ for all }j\geq0
\end{equation}
because $\tau^{-j}(K)\leq K$ for every $j\geq0$. If $< \tau^j h,k>$ were equal to $0$ for every $j\geq 0$ as well, then $K$ would be finite which is impossible under the hypothesis that $\tau(K)\gvertneqq K$. Hence there is a $j\in{\mathbb Z}^+$ such that $< \tau^jh,k>\, =\, < h,\tau^jk>\, \ne0$ for some $k\in K$. 

Choose $h^*$ satisfying (\ref{eqn:Kexpand}) such that the difference between the maximum and minimum integers in the support of $h^*$ is minimized. Denote by $M$ and $m$ respectively these maximum and minimum integers and put $d:=M-m$. By translating $h^*$ if necessary, it may also be supposed that $\langle  h^*, \tau k\rangle \ne0$ for some $k\in K$. We claim:
\begin{equation}
\label{eq:Claim1}
\hbox{\emph{the image of the projection $K\to \prod_{-M}^{-m}\field_2\, :\, k\mapsto k|_{[-M,-m]}$ has order $2^d$}.}
\end{equation}
For, if not, there is $h\ne h^*$ in $\bigoplus_{l\in {\mathbb Z}} \field_2$ with support contained in $[m,M]$ and satisfying~(\ref{eqn:Kexpand}). Then there is $j\in{\mathbb Z}$ such that $h^*-\tau^j(h)$ satisfies~(\ref{eqn:Kexpand}) and has length of support less than $d+1$, contradicting the choice of $h^*$.  

Since $h^*(m) = 1$, it follows from (\ref{eq:Claim1}) that the projection 
$$
\textstyle{K\to \prod_{-M}^{-m-1}\field_2 \,:\, k\mapsto k|_{[-M,-m-1]}} \hbox{ is surjective.}
$$
 Hence, if $k\in K$ satisfies $\langle h^*, \tau k\rangle \ne0$, then $k'\in K$ may be chosen such that $k'|_{[-M,-m-1]} = \tau k|_{[-M,-m-1]}$ and so, putting $k_1 = k-\tau^{-1}k'$, we have $\langle h^*, \tau k_1\rangle \ne0$ and $k_1|_{[-M+1,-m]} = 0$. Then necessarily $k_1(-m+1) = 1$ and, by (\ref{eq:iterated}), $k_1(j) = 0$ for every $j\leq -m$. Since $K$ is invariant under $\tau^{-1}$, it follows that for every $n\geq1$ there is $k_n\in K$ with 
 $$
 k_n(-m+n) = 1\hbox{ and }k_n(j) = 0\hbox{ for every }j< -m+n.
 $$
Therefore, since $K$ is closed, $\prod_{j\geq -m+1} \field_2$ is a subgroup of $K$ and we may take $J=-m+1$.
\end{demo}

\medskip
Returning to the proof of Lemma~\ref{prop:reduced_case}.\ref{prop:RC2}, let $K\leq \prod_{l\in {\mathbb Z}} \field_2$ be commensurated by $\tau$. Then $G : = K +\, \bigoplus_{\mathbb Z} \field_2$ is a subgroup of $\prod_{\mathbb Z} \field_2$ and is equal to the union of all subgroups that are commensurable with $K$. Topologise $G$ by defining ${\mathcal U}\subset G$ to be open if ${\mathcal U}\cap (f+K)$ is open in  $f+K$ for every $f\in G$. Then $G$ is a totally disconnected locally compact group, is stable under $\tau$ and the restriction of $\tau$ to $G$ is an automorphism.  Lemma~1 in \cite{Wi94} shows that there is a natural number $n$ such that $K' := \bigcap_{j=0}^n \tau^j(K)$ satisfies $K' := K'_+K'_-$, where $K'_\pm = \bigcap_{j=0}^\infty \tau^{\pm j}(K')$. It is immediate from the definitions that $K'$ is commensurable with $K$ and that $\tau(K'_+)\geq K'_+$ and $\tau(K'_-)\leq K'_-$. If both $\tau(K'_+) = K'_+$ and $\tau(K'_-) = K'_-$, then both subgroups are finite, by Lemma~\ref{lem:stable_finite}, in which case $K$ is finite. Otherwise, at least one of  $\tau(K'_+) \gvertneqq K'_+$ and $\tau^{-1}(K'_-) \gvertneqq K'_-$ holds. Suppose, without loss of generality, that  $\tau(K'_+) \gvertneqq K'_+$. Then, by Lemma~\ref{lem:Kexpand}, there is an integer $J$ such that $\prod_{j\geq J} \field_2 \leq K'_+ \leq K$. Since $K$ is commensurated by $\sigma$ and $\sigma(\prod_{j\geq J} \field_2) = \prod_{j\leq -J} \field_2$, it follows that a finite index subgroup of $\prod_{j\leq -J} \field_2 \times \prod_{j\geq J} \field_2$ is contained in $K$. Therefore $K$ has finite index in $\prod_{j\in {\mathbb Z}} \field_2$ in this case. 
\end{demo}

\medskip
This completes the proof of Lemma~\ref{prop:reduced_case}, which establishes Proposition~\ref{prop:comm-norm} in the special case when  $\Xi\leq \bigoplus_{l\in {\mathbb Z}} \field_2$. For the general case, let $\Xi\leq \Gamma$ be commensurated by $\langle \tau,\sigma\rangle$. Since $\Lambda$ has index~2 in $\Gamma$, the intersection of $\Xi$ with $\Lambda$ is commensurable with $\Xi$ and so it may be assumed that $\Xi\leq \Lambda$. Then $\Xi =   Z \ltimes(\Xi\cap \bigoplus_{l\in {\mathbb Z}} \field_2)$, where $Z$ is a cyclic subgroup of~$\Lambda$. Since the case when $\Xi\leq \bigoplus_{l\in {\mathbb Z}} \field_2$ has already been treated, it may be supposed that $Z$ is generated by $(\tau^p,h')$, where $p\ne0$. Then $hZh^{-1}\cap \Xi$ has finite index in $hZh^{-1}$ for any $h\in \bigoplus_{l\in {\mathbb Z}} \field_2$, and so, fixing a non-zero $h$, there is a non-zero integer $r$ such that the commutator $[(\tau^p,h')^{kr},h] = \tau^{krp}h - h$ belongs to $\Xi$ for every $k\in {\mathbb Z}$. Hence $\Xi\cap \bigoplus_{l\in {\mathbb Z}} \field_2$ is infinite. Since $\bigoplus_{l\in {\mathbb Z}} \field_2$ is stabilized by $\langle \tau, \sigma \rangle$, $\Xi\cap \bigoplus_{l\in {\mathbb Z}} \field_2$ is also commensurated by this group and so, by Lemma~\ref{prop:reduced_case}, has finite index in $\bigoplus_{l\in {\mathbb Z}} \field_2$. Therefore  $\Xi$ has finite index in $\Gamma$, thus completing the proof.
\end{demo}

\medskip
\noindent {\bf Remark.} It is still true that the inner commensurator-normalizer property passes
from a finite index subgroup to the ambient group.

\medskip

\section{Proof of Theorem~\ref{thm:mainresult}}
\label{sec:Thm1.3Proof}
\subsection{{\bf Proof of part 1}}
In view of Proposition~\ref{prop:comm_norm}, the following result suffices in order to establish part~{\bf 1} of Theorem~\ref{thm:mainresult}.
 \begin{theorem}
 \label{thm:outer}
Let ${\mathbf G}$ be a Chevalley group over a number field $K$. 
If ${\mathbf G} \cong {\mathbf{SL}}_2$ assume that $K$ is not $\QQ$ or an imaginary
quadratic extension of it. Let ${\mathcal O}$ be the ring of integers in $K$. Then ${\mathbf G}({\mathcal O})$ has the outer commensurator-normalizer property.
 \end{theorem}
For the proof, let $\varphi:{\mathbf G}({\mathcal O}) \to \Delta$ be a homomorphism and suppose  that
$\Lambda < \Delta$ is commensurated by $\varphi({\mathbf G}({\mathcal O}))$. Replacing $\Delta$ by  the group generated by $\Lambda$ and $\varphi({\mathbf G}({\mathcal O}))$, it may be supposed that $\Lambda$ is commensurated by~$\Delta$. Put $D =\RCom{\Delta}{\Lambda}$ and let $L$ be the closure of $\Rmap{\Delta}{\Lambda}(\Lambda)$ in $D$. Then the composite of $\varphi$ and the map $\Rmap{\Delta}{\Lambda} : \Delta \to D$ yields a homomorphism $\bar\varphi : {\mathbf G}({\mathcal O}) \to D$. Since 
$\bar\varphi^{-1}(V)$ is commensurable with $\Lambda$ for any compact open subgroup $V$ of $D$, it suffices to show that $\bar\varphi({\mathbf G}({\mathcal O}))$ normalizes a compact, open subgroup of $D$. As a first step, the next proposition will be used to show that each root subgroup of ${\mathbf G}({\mathcal O})$ normalizes a compact, open subgroup (that may depend on the root).

 \begin{proposition}
 \label{prop:roots}
Let ${\mathbf G}$ be a Chevalley group over $K$ and ${\mathcal O}\subset K$ 
its ring of integers. 
\begin{enumerate}
\item Suppose that the rank of ${\mathbf G}$ is at least~2. For every root subgroup, $X_\alpha$, $\alpha\in \Phi$, of ${\mathbf G}({\mathcal O})$ there is a finitely generated nilpotent group, $N_\alpha$, such that $X_\alpha\cap [N_\alpha,N_\alpha]$ has finite index in $X_\alpha$.
\label{prop:roots1}
\item Suppose that ${\mathbf G} = \mathbf{SL_2}$ where $K$ is not ${\mathbb Q}$ or an imaginary quadratic extension of ${\mathbb Q}$. Then corresponding to each root subgroup, $X_i$, $i=1,2$, of ${\mathbf G}({\mathcal O})$ there is a polycyclic group  $P_i$  such that $X_i\cap [P_i',P_i']$ has finite index in $X_i$ for every finite index subgroup $P_i'$ of $P_i$. 
\label{prop:roots2}
\end{enumerate}
 \end{proposition}
 \begin{demo}
 \ref{prop:roots1}. Since the rank of ${\mathbf G}$ is at least~2, each $\alpha\in \Phi$ is the sum, $\alpha = m\beta + n\gamma$, for some $\beta,\gamma\in\Phi$ that are not opposite roots and $m,n\in {\mathbf Z}^+$. Let $N$ be the group generated by $X_\beta$ and $X_\gamma$. Then $N$ is an arithmetic subgroup of the unipotent group, $U$, generated by $X_\beta(K)$ and $X_\gamma(K)$, and is finitely generated and nilpotent. Proposition~7.2.4 in~\cite{Abels} applies to show that $[N,N]$ is an arithmetic subgroup of $[U,U]$. Since $[U,U]$ contains $X_\alpha(K)$, it follows that $[N,N]\cap X_\alpha$ has finite index in $X_\alpha$. 
 
 
 \ref{prop:roots2}. Let $X_i = \left\{ {\sf x}_i(a) \mid a\in {\mathcal O}\right\}$, $i=1,2$ be the root subgroups of $\mathbf {SL_2}({\mathcal O})$ and let $T$ be the subgroup of diagonal elements. Then $X_i$ is isomorphic to the additive group of ${\mathcal O}$ and $T$ to  the group of multiplicative units in ${\mathcal O}$, so that both are finitely generated abelian groups. Denote the elements of $T$ by ${\sf t}(v)$ where $v$ is a unit. Let $P_i$ be the group generated by  $T$ and $X_i$. Then $P_i$ is the semi-direct product of these subgroups and is polycyclic. 
 
The hypotheses on $K$ imply that its group of units  has elements of infinite order, by Dirichlet's Unit Theorem \cite[\S4.5, Corollary~1]{PR}. Hence, for every finite index subgroup $P_i'$ of $P_i$,  the intersection $P_i'\cap T$ contains an element, ${\sf t}(v)$ say, of infinite order. For every ${\sf x}(a)\in  X_i' := X_i\cap P_i'$ we have $[{\sf t}(v),{\sf x}(a)] = {\sf x}_i((v^{\pm2}-1)a)$, where the sign of the exponent of $v$ depends on whether $i=1$ or $2$. Hence the map ${\sf x}(a) \mapsto [{\sf t}(v),{\sf x}(a)]
$ is an endomorphism of $X_i'$ whose range has index equal to the norm of $v^{\pm2}-1$ and the claim holds. 
 \end{demo}
 
\medskip
Suppose now that $X_\alpha$ is a root subgroup of ${\mathbf G}({\mathcal O})$ and  consider first the case when the rank is at least~2. Let $N$ be a finitely generated nilpotent subgroup such that $X_\alpha\cap [N,N]$ has finite index in $X_\alpha$. Then $\bar\varphi(N)$ is a flat subgroup of $D$ by Theorem~\ref{thm:nilpotent} and $\bar\varphi([N,N])$ is uniscalar by Theorem~\ref{thm:flat}. Hence $X_\alpha$ itself is flat and uniscalar and there is a compact, open subgroup of $D$ that is normalized by $\bar\varphi(X_\alpha)$. In the case of $\mathbf{SL_2}({\mathcal O})$, let $P$ be the polycyclic group of the proposition corresponding to $X_\alpha$. Then, by Theorem~\ref{thm:polycyclic}, there is a finite index subgroup $P'$ of $P$ such that $\bar\varphi(P')$ is flat and, by  Theorem~\ref{thm:flat}, $\bar\varphi([P',P'])$ is uniscalar. Hence there is a compact, open subgroup of $D$ that is normalized by $X_\alpha\cap [P',P']$. Since $X_\alpha\cap [P',P']$ has finite index in $X_\alpha$, there is a compact, open subgroup of $D$ that is normalized by $\bar\varphi(X_\alpha)$. The next corollary has thus been established in all cases. 
 
 \begin{corollary}
 \label{cor:rootflat}
Let $\bar\varphi: {\mathbf G}({\mathcal O}) \to D$ be as above. Then for every root subgroup, $X_\alpha$, of ${\mathbf G}({\mathcal O})$ there is a compact open subgroup $V_\alpha \leq D$ that is normalized by $\bar\varphi(X_\alpha)$. 
\endproof
 \end{corollary}

If the subgroup $V_\alpha$ did not depend on $\alpha$, the proof would be  complete because ${\mathbf G}({\mathcal O})$ is  generated by its root subgroups. However ${\mathbf G}({\mathcal O})$ satisfies the stronger property of \emph{bounded generation} relative to its root subgroups, that is: there is an $N$ and a sequence $\{\alpha_n\}_{n=1}^N$ of roots such that ${\mathbf G}({\mathcal O}) = X_{\alpha_1}X_{\alpha_2} \dots X_{\alpha_N}$. It will be seen that this stronger property allows us to find the desired compact, open subgroup of $D$ that is normalized by $\bar\varphi({\mathbf G}({\mathcal O}))$. The following is a compilation of cases from \cite[Theorem~1.2]{Mor} and  \cite[Corollary~1]{Tav}, which extend work in \cite{CarKel}. 

\begin{theorem}
\label{thm:boundedgen}
Let ${\mathbf G}$ be a Chevalley group with irreducible root system $\Phi$ over an algebraic number field $K$ and let ${\mathcal O}$ be the ring of integers in $K$. Suppose that either:
\begin{enumerate}
\item $\hbox{\rm rank}\ \Phi \geq 2$, or 
\item $\hbox{\rm rank}\ \Phi = 1$ and ${\mathcal O}$ has infinitely many units. 
\end{enumerate}
Then ${\mathbf G}({\mathcal O})$ has bounded generation relative to its root subgroups, $X_\alpha$, $\alpha\in \Phi$. 
\endproof
\end{theorem}

Note that by Dirichlet's unit theorem ${\mathcal O}$ has infinitely many units provided that $K$ is not ${\mathbb Q}$ or an imaginary quadratic extension of ${\mathbb Q}$; cf. \cite[\S4.5, Corollary~1]{PR}. 

The metric, $d$, defined on ${\mathcal B}(D)$ (the set of compact open subgroups of $D$) by
$$
d(V_1,V_2) = \log [V_1:V_1\cap V_2][V_2:V_1\cap V_2], \quad V_1,V_2\in {\mathcal B}(D),
$$
will be useful for the next step in the argument (see the proof of 
Proposition~\ref{prop:distor} below). Note that, for any $x\in D$, the conjugation map 
$$
V\mapsto xVx^{-1} \quad (V\in {\mathcal B}(D))
$$ 
is an isometry of $({\mathcal B}(D),d)$.

For each root subgroup, $X_\alpha$, of ${\mathbf G}({\mathcal O})$ choose a compact open subgroup,  $V_\alpha$, of $D$ that is normalized by $\bar\varphi(X_\alpha)$ and set
$$
M = \max \left\{ d(L,V_\alpha) \mid \alpha\in \Phi \right\}.
$$
Then for any $\alpha\in \Phi$ and $x\in X_\alpha$, 
\begin{eqnarray*}
d(\bar\varphi(x)L\bar\varphi(x^{-1}), L) &\leq& d(\bar\varphi(x)L\bar\varphi(x^{-1}), V_\alpha) + d(V_\alpha, L) \\
 &=& d(\bar\varphi(x)L\bar\varphi(x^{-1}),  \bar\varphi(x)V_\alpha\bar\varphi(x^{-1})) + d(V_\alpha, L) \\
 &\leq& 2M.
\end{eqnarray*}
Induction on $j$ shows that 
$$
d(\bar\varphi(x)L\bar\varphi(x^{-1}), L) \leq 2jM \quad (x\in X_{\alpha_1}X_{\alpha_2} \dots X_{\alpha_j})
$$
and so, by Theorem~\ref{thm:boundedgen}, there is an $N$ such that 
$$
d(\bar\varphi(x)L\bar\varphi(x^{-1}), L) \leq 2NM \quad (x\in {\mathbf G}({\mathcal O})).
$$
In other words, the orbit of $L\in {\mathcal B}(D)$ under the conjugation action of $\bar\varphi({\mathbf G}({\mathcal O}))$ is bounded. It follows by \cite[Theorem 6(iii)]{BerLen} (which is quite non-trivial, unlike standard ``circumcentre arguments'') that there is $L'\in {\mathcal B}(D)$ which is stable under the conjugation action of $\bar\varphi({\mathbf G}({\mathcal O}))$, thus completing the proof of the first part of Theorem~\ref{thm:mainresult}. (See also \cite{Sch1} and \cite{Wagner} for work 
related to \cite{BerLen}.)

\medskip

 \begin{remark}
 \label{rem:positivechar}
Let $\field_p$ be a finite field and $\field_p[t]$ denote the ring of 
polynomials over $\field_p$. For $n>2$ the higher rank arithmetic
group in positive characteristic $SL_n(\field_p[t])$ satisfies the Margulis' NST.
Moreover, the algebraic treatment of Venkataramana in \cite{Venkat} 
(already mentioned earlier), shows that every {\it commensurated} subgroup of it is 
finite or of finite index (in particular, this group has the inner commensurator-normalizer 
property). However, unlike its ``cousin'' $SL_{n}(\ZZ)$, 
this group does not satisfy the outer property 
because $SL_n(\field_p[t])$ and $SL_n(\field_p[t^{-1}])$ are both 
(abstractly isomorphic) subgroups of the countable group 
$SL_n(\field_p[t,t^{-1}])$ which commensurate, without almost normalizing,
each other (as can easily be verified). In fact the following can be proved.
  \end{remark}

\begin{theorem}
\label{thm:char}
Let $\Gamma$ be a finitely generated group with the outer commensurator-normalizer property.
Then any linear representation of $\Gamma$ over a field
of positive characteristic has finite image.
\end{theorem}

\begin{demo}
It 
is a rather standard fact, cf.~\cite[Lemma 5.7,5.8]{Lu-Zi}, that
if a finitely generated group $\Gamma$ admits
an infinite representation over such a field, then it admits also some 
unbounded representation 
into $GL_n(k)$ for some local field $k$ of the same (positive) characteristic. 
Hence, by passing to a finite index subgroup and using the implied finite abelianization 
property of $\Gamma$, we may assume the latter representation ranges in  
$SL_n(k)$. However,   
 normalizers
of compact open subgroups of $SL_n(k)$ are (open and different from $SL_n(k)$ hence) 
compact, while the image of $\Gamma$ 
is unbounded, in contradiction to the assumed
outer commensurator-normalizer property.
\end{demo}

\smallskip

\subsection {\bf Proof of part 2 of Theorem~\ref{thm:mainresult}.} 
Let $V^\infty \subseteq S \subseteq V$ be a set of places. Let $\Gamma$ be an $S$-arithmetic subgroup of ${\mathbf G}(K)$  and suppose that  $\Lambda\leq \Gamma$ is commensurated by $\Gamma$. Then $\Lambda\cap {\mathbf G}({\mathcal O}_S)$ is commensurable with $\Lambda$ and it suffices to show that $\Lambda\cap {\mathbf G}({\mathcal O}_S)$ is a standard commensurated subgroup. Hence, replacing $\Gamma$ and $\Lambda$ by $\Gamma\cap {\mathbf G}({\mathcal O}_S)$ and $\Lambda\cap {\mathbf G}({\mathcal O}_S)$ respectively, it may be supposed that $\Gamma$ and $\Lambda$ are subgroups of ${\mathbf G}({\mathcal O}_S)$. Consider $\Gamma\cap {\mathbf G}({\mathcal O})$, which is commensurable with ${\mathbf G}({\mathcal O})$ and commensurates $\Lambda$. Part~1. shows that there is $\Lambda'\leq \Gamma$ that is commensurable with $\Lambda$ and normalized by $\Gamma\cap {\mathbf G}({\mathcal O})$. The proof will be completed by showing that $\Lambda'$ is a standard commensurated subgroup.

Because both $\Gamma\cap {\mathbf G}({\mathcal O})$ and  $\Lambda '$ are commensurated by
$\Gamma$, and the former normalizes the latter, it is easy to verify that the subgroup
they generate $(\Gamma\cap {\mathbf G}({\mathcal O}))\Lambda'$ is also 
commensurated by $\Gamma$ (note that this is not the case in general).
Thus both $\Gamma\cap {\mathbf G}({\mathcal O})$ and $(\Gamma\cap {\mathbf G}({\mathcal O}))\Lambda'$ are commensurated subgroups of $\Gamma$ and 
$$
\Gamma\cap {\mathbf G}({\mathcal O}) \leq (\Gamma\cap {\mathbf G}({\mathcal O}))\Lambda'.
$$ 
By Lemma~\ref{lem:nestedcommensurated}, there is a homomorphism 
\begin{equation*}
\rho : \RCom{\Gamma}{(\Gamma\cap {\mathbf G}({\mathcal O}))} \to \RCom{\Gamma}{(\Gamma\cap {\mathbf G}({\mathcal O}))\Lambda'}
\end{equation*} 
that has dense range. By Proposition~\ref{prop:G(A)}, $\RCom{\Gamma}{(\Gamma\cap {\mathbf G}({\mathcal O}))}$ is isomorphic to ${\mathbf G}(\adele^{S_f})$ 
divided by its center (Note: the notation here and in the sequel 
follows the one preceding Proposition~\ref{prop:G(A)}).
For convenience, let us
establish at this point throughout the rest of the proof the {\it ad hoc\/} notation $G^p$ for the ``projectivization'' 
of a group $G$ (always linear
in our discussion), i.e., $G^p$ is $G$ divided by its center.  
 Thus, using the isomorphism provided by Proposition~\ref{prop:G(A)}  we have a homomorphism
\begin{equation}
\label{eq:HOM}
\rho ' : {\mathbf G}(\adele^{S_f})^p \to \RCom{\Gamma}{(\Gamma\cap {\mathbf G}({\mathcal O}))\Lambda'}
\end{equation} 
that has dense range. That $\rho '$ is in fact surjective follows from the next result which may be of 
independent interest. 

\begin{theorem}
\label{thm:closedrange}
Let $S \subset V$ be any set of places, and 
$\rho : {\mathbf G}(\adele^{S}) \to H$ be a continuous homomorphism, 
where $H$ is locally compact. 
Then $\rho({\mathbf G}(\adele^{S}))$ is closed. 
\end{theorem}

In order not to disrupt the flow of the proof of Theorem~\ref{thm:mainresult}, we postpone the proof of this result
to the end. Continuing, Theorem~\ref{thm:closedrange} implies that the 
homomorphism $\rho '$ in (\ref{eq:HOM}) is surjective, and so 
$\RCom{\Gamma}{(\Gamma\cap {\mathbf G}({\mathcal O}))\Lambda'} \cong 
\rho '({\mathbf G}(\adele^{S_f})^p)$. Moreover, since ${\mathbf G}(\adele^{S_f})^p$ is the 
restricted product of the simple locally compact groups ${\mathbf G}(K_\nu)^p\ (\nu\in S_f)$, there is $S_f'\subset S_f$ (possibly empty) such that $\ker(\rho ') = {\mathbf G}(\adele^{S_f'})^p$ and 
$$
\rho '({\mathbf G}(\adele^{S_f})^p) \cong {\mathbf G}(\adele^{S_f})^p/{\mathbf G}(\adele^{S_f'})^p 
\cong {\mathbf G}(\adele^{S_f\setminus S_f'})^p.
$$ 

Collecting all the homomorphisms obtained so far yields:
$$
\Gamma  \stackrel{{\Rmap{\Gamma}{\Gamma\cap {\mathbf G}({\mathcal O})}}}   \longrightarrow 
\RCom{\Gamma}{(\Gamma\cap {\mathbf G}({\mathcal O}))} \cong  {\mathbf G}(\adele^{S_f})^p   
\stackrel{\rho '}  \longrightarrow 
\RCom{\Gamma}{(\Gamma\cap {\mathbf G}({\mathcal O}))\Lambda'} \cong 
{\mathbf G}(\adele^{S_f\setminus S_f'})^p
$$

The final isomorphism in the sequence is not necessarily unique, and is chosen so that the 
composition of the last two homomorphisms is the 
natural projection. Call the composition of all these maps $f:\Gamma \to {\mathbf G}(\adele^{S_f\setminus S_f'})^p$.
Lemma~\ref{lem:nestedcommensurated} then implies that $(\Gamma\cap {\mathbf G}({\mathcal O}))\Lambda'$ is
 commensurable with $f^{-1}(V)$  for any compact open subgroup, $V$, 
of the last group.
Choosing $V = {\mathbf G}({\Kprod}^{S_f\setminus S_f'})^p$, by tracing back all the arrows
in the definition of $f$ it is easy to see that   
$f^{-1}(V)  = \Gamma \cap ({\mathbf G}(\adele^{S_f'}) \times  
{\mathbf G}(\Kprod^{S_f\setminus S_f'})$). 
Hence $(\Gamma\cap {\mathbf G}({\mathcal O}))\Lambda'$ is commensurable with 
$$
\Gamma \cap ({\mathbf G}(\adele^{S_f'}) \times  {\mathbf G}
(\Kprod^{S_f\setminus S_f'})) <_{com} {\mathbf G}({\mathcal O}_{S'}):= {\mathbf G}({\mathcal O}_{V^{\infty}\cup S_f'}).
$$
We have now reached the reduction step to the Margulis' Normal Subgroup Theorem. 
Since $\Lambda'$ is normalized by $\Gamma\cap {\mathbf G}({\mathcal O})$, it is a normal 
subgroup of the $S'$-arithmetic group $(\Gamma\cap {\mathbf G}({\mathcal O}))\Lambda'$. 
Therefore, by the NST 
\cite[Theorem~VIII.2.6 page 265]{Mar2}, $\Lambda'$ is either finite, or has finite index, 
i.e., it is a standard commensurated subgroup.  We are thus only left with the following:

\medskip

\noindent {\bf Proof of Theorem~\ref{thm:closedrange}.}
The closedness $\rho({\mathbf G}(\adele^{S}))$ is well known to 
experts for finite~$S$. As we could not find a reference the proof is 
sketched here for 
completeness. Suppose that $S= \left\{\nu_1,\dots, \nu_s\right\}$. 
Then ${\mathbf G}(\adele^{S}) = \prod_{j=1}^s {\mathbf G}(K_{\nu_j})$, where each ${\mathbf G}(K_{\nu_j})$ is a projectively simple group over the locally 
compact field $K_\nu$. Let the sequence $\{g_n\}\subset {\mathbf G}(\adele^{S})$ be such that $\rho(g_n)\to x\in H$ as $n\to\infty$. It must be shown that there is $g\in {\mathbf G}(\adele^{S})$ such that $\rho(g) = x$. The Cartan decomposition, see \cite[Theorem~3.14]{Rap} or \cite[Corollary 2.17]{IwMa}, yields that
$$
{\mathbf G}(K_{\nu_j}) = {\mathbf G}(O_{\nu_j})H_j{\mathbf G}(O_{\nu_j}),
$$
where $O_{\nu_j}$ is the ring of integers in $K_{\nu_j}$ and where $H_j$ the Cartan subgroup of ${\mathbf G}(K_{\nu_j})$. Hence, since ${\mathbf G}(O_{\nu_j})$ is compact, it may be supposed that the projection of $g_n$ onto ${\mathbf G}(K_{\nu_j})$ belongs to $H_j$ for each~$j$. The projection of $g_n$ onto ${\mathbf G}(K_{\nu_j})$ will be denoted by $g_{n,j}$. We claim that 
$$
\rho({\mathbf G}(K_{\nu_j})) = \triv\hbox{ for any }j\hbox{ with }\{g_{n,j}\}\hbox{ unbounded}.
$$ 
To see this, assume that $\{g_{n,j}\}$ is unbounded. Then there is there is a root, $\alpha$, such that $\left\{\chi_\alpha(g_{n,j})\right\}$ is unbounded, where $g{\sf x}_\alpha(a)g^{-1} = \chi_\alpha(g){\sf x}_\alpha(a)$ for every $a\in K_\nu$. Hence the identity is an accumulation point of $\left\{ g_{n,j}{\sf x}_{-\alpha}(a)g_{n,j}^{-1}\right\}$ and it follows that $x\rho({X}_{-\alpha}(K_{\nu_j}))x^{-1} = \triv$. Therefore $\rho({X}_{-\alpha}(K_{\nu_j})) = \triv$ and the claim follows because ${\mathbf G}(K_{\nu_j})$ is projectively simple. Since ${\mathbf G}(K_{\nu_j})$ is locally compact for each $j$, the claim implies that there is $g\in {\mathbf G}(\adele^{S})$ such that $\rho(g) = x$. 

Suppose next that $S$ is infinite and let, as before, 
the sequence $\{g_n\}\subset {\mathbf G}(\adele^{S})$ be such that $\rho(g_n)\to x\in H$ as $n\to\infty$. We shall assume hereafter that $H$ is totally disconnected, which is the only case needed for our purposes. See however the remarks proceeding the proof for the general case. Because of this assumption on $H$ we may assume that
$S$ consists of finite places (a more significant use of this assumption is below).
For $\nu\in {S}$ denote the projection ${\mathbf G}(\adele^{S})\to 
{\mathbf G}(K_\nu)$ by $\pi_\nu$. Let ${\mathbf G}(\Kprod^{S})$ denote the product $\prod_{\nu\in S} {\mathbf G}({O}_\nu)$. 
Then $\rho({\mathbf G}(\Kprod^{S}))$ is a compact subgroup of $H$ and so 
there is a compact open subgroup, $V\leq H$, such 
that $\rho({\mathbf G}(\Kprod^{S}))\leq V$.
Since $\rho({\mathbf G}(\Kprod^{S}))$ is compact and $V$ is open, 
$x\rho({\mathbf G}(\Kprod^{S}))x^{-1}\cap V$ has finite index in $V$. Moreover, since $\rho(g_n)\to x$, it may be supposed that 
$$
\rho(g_n{\mathbf G}(\Kprod^{S})g_n^{-1})\cap V = x\rho({\mathbf G}(\Kprod^{S}))x^{-1}\cap V \hbox{ for all }n.
$$
Hence $(g_n{\mathbf G}(\Kprod^{S})g_n^{-1})\cap \rho^{-1}(V)$ 
is independent of $n$. Put 
\begin{equation}
\label{eq:notinV}
P = \left\{ \nu\in S \mid g_n{O}_\nu g_n^{-1} \not\leq  \rho^{-1}(V)\right\}.
\end{equation}
Then $P$ is a finite subset of $S$ because $(g_n{\mathbf G}(\Kprod^{S})g_n^{-1})\cap \rho^{-1}(V)$ has finite index in $g_n{\mathbf G}(\Kprod^{S})g_n^{-1}$. 

Denote by $\pi_P$ and $\pi_{S\setminus P}$ the natural projections 
$$
\pi_P :  {\mathbf G}(\adele^S) \to {\mathbf G}(\adele^P) \hbox{ and } 
\pi_{S\setminus P} : {\mathbf G}(\adele^S) \to {\mathbf G}(\adele^{S\setminus P}).
$$
Then $g_n$ may be factored as $g_n = a_nb_n$ for each $n$, where $a_n\in \ker(\pi_{S\setminus P})$ and $b_n\in \ker(\pi_P)$. Consider $\nu\in S\setminus P$. If $\pi_\nu(g_n)\not\in {\mathbf G}({O}_\nu)$, then $g_n{\mathbf G}({O}_\nu) g_n^{-1} \ne {\mathbf G}({O}_\nu)$ and  $\langle g_n{\mathbf G}({O}_\nu) g_n^{-1}, {\mathbf G}({O}_\nu)\rangle$ has non-compact closure, by \cite[Corollary 2.18]{IwMa}. On the other hand, it follows from (\ref{eq:notinV}) that $\rho(\langle g_n{\mathbf G}({O}_\nu) g_n^{-1}, {\mathbf G}({O}_\nu)\rangle)$ is contained in the compact set $V$, which forces $\rho({\mathbf G}(K_\nu))$ to be trivial. Dividing $S_f\setminus P$ into those places where $\rho$ is trivial and its complement, $b_n$ may be factored as $b_n = c_nd_n$, where $c_n\in  {\mathbf G}( {\Kprod}^{S})$ and $\rho(d_n)$ is trivial. Since 
$ {\mathbf G}({\Kprod}^{S})$ is compact, it may be supposed, by passing to a subsequence if necessary, that $\{c_n\}$ is convergent with limit $c$, say.  Hence $\rho(a_n) = \rho(g_nd_n^{-1}c_n^{-1}) \to x\rho(c^{-1})$ as $n\to\infty$. Since $a_n\in {\mathbf G}(\adele ^P)$ and $P$ is finite, there is $a\in {\mathbf G}(\adele^P)$ such that $\rho(a) = x\rho(c^{-1})$. Therefore $x = \rho(ac)$ as required, completing the proof of 
Theorem~\ref{thm:closedrange}. \endproof

\medskip

\noindent {\bf Remark.} The case of infinite $S$ in Theorem~\ref{thm:closedrange} does not follow directly from the finite case, and special 
properties of the groups ${\mathbf G}(K_\nu)$ are essential for the argument. 
Consider, for example, the restricted direct product of any infinite family of non-trivial 
finite groups w.r.t the trivial subgroup, which is just their direct sum. This group
 embeds densely into the direct product of the groups. 

\medskip

\noindent {\bf Remark.} The assumption made in the proof that $H$ is 
totally disconnected, can indeed be relaxed. The full treatment of the combined
totally disconnected by connected case in this spirit 
seems rather technical (and may depend also on the Montgomery-Zippin work~\cite{MontZip}).
Alternatively, one can give a different, unified, ``high tech'' representation theoretic 
proof of Theorem~\ref{thm:closedrange} along the following lines: 
First, observe that closedness of the image of any 
group
homomorphism $\rho$ with Ker$\rho$ compact is equivalent to
the property that the restriction of the $H$-regular representation $L^2(H)$ to Im$\rho$ is a 
mixing unitary representation. 
This and the well known Howe-Moore Theorem, immediately imply the closedness
of $\rho$ restricted to any simple factor (noticing that an invariant vector is possible
only if the image is bounded, which again contradicts Howe-Moore in view of the Peter-Weyl Theorem).
Hence the matrix coefficients of $L^2(H)$ restricted to each simple factor, unless trivial, 
satisfy the well known
uniform {\it pointwise} exponential decay property, which transfers also 
to (almost) every irreducible  component (appearing in 
a direct integral decomposition of $L^2(H)$). 
Using this and the well 
known (yet key) fact that every irreducible unitary representation of ${\mathbf G}(\adele ^{S})$ is a tensor product
of irreducible (tempered in our case) representations of the local factors, one shows that all
{\it irreducible} components of $L^2(H)$ restricted to  ${\mathbf G}(\adele_{S})$ are mixing, 
which then passes to their direct integral. 

The same argument shows also the following
result about adelic groups, which is well known for (finitely many) simple factors, but requires 
more attention in general. 
\begin{proposition}
\label{prop:adelic-open}
For any ${\mathbf G}$ absolutely simple and simply connected, and any set 
$S \subset V$, any open subgroup
$U<{\mathbf G}(\adele^{S})$ is of the form
 ${\mathbf G}(\adele^{S'}) \times L$ where 
$L<{\mathbf G}(\adele^{S\setminus S'})$ is {\it compact} open. 
\endproof
\end{proposition}
Once again, compare to the situation with a direct sum
of finite groups. Representation theory becomes relevant to this Proposition through the observation that for
an open subgroup $B<A$ the $A$-representation $L^2(A/B)$ is mixing iff $B$ is compact.
 For the proof one looks at the unitary 
representation on $L^2({\mathbf G}(\adele^{S})/U)$, which on every simple factor to which it restricts
 non-trivially
must be mixing (by the individual properties of the simple factors), and now 
continuing as before. One can give here as well a more elementary argument.

\medskip
\noindent {\bf Remark.} One can use Proposition~\ref{prop:adelic-open} in order to prove the following
 result (under the same general conditions on ${\mathbf G}$): 
\begin{proposition}
\label{prop:arithmetic}
If $\Gamma <{\mathbf G}(K)$ contains an infinite $S$-arithmetic subgroup $\Lambda$ (for any $S$
containing the infinite places), then
$\Gamma$ itself is $S'$-arithmetic.  
\endproof
\end{proposition}
The result was proved by 
Lubotzky-Zimmer in \cite[Lemma 2.8]{Lu-Zi} when $\Lambda$ was assumed
 commensurable with ${\mathbf G}(\mathcal O)$ (assumed infinite), and $\Gamma$ is finitely 
generated. 
The first condition is easier to relax. 
 However
to argue similarly to Lubotzky-Zimmer in removing the finite generation condition on $\Gamma$ 
one needs (or at least  can use) 
Proposition~\ref{prop:adelic-open}.

\medskip
\noindent {\bf One Final Remark.} One can use  Proposition~\ref{prop:arithmetic} instead of 
Theorem~\ref{thm:closedrange} in the proof of Theorem~\ref{thm:mainresult}. As suggested by the preceding discussion, they are strongly related. 
Theorem~\ref{thm:closedrange} is also an integral part of Theorem~\ref{thm:adelic}.

\subsection{Alternative approaches, and a question} 
\label{sec:alternative}

We begin by noting that while the approaches are indeed different, 
like us Venkataramana~\cite{Venkat} makes crucial use of the
unipotent elements, leaving the anisotropic cases unapproachable.
Note that part {\bf 2} of our Theorem~\ref{thm:mainresult} and his result~\cite{Venkat} overlap,
but neither one implies the other (while over positive characteristic 
our approach is hopeless). Note also the structural similarity of
the proof of {\bf 1} in Theorem~\ref{thm:mainresult} and that 
in~\cite[Prop 7.14]{PR} and~\cite{Sh3}. We now proceed
to discuss two variations which continue to rely on bounded generation 
as well
as on ``topologization'' of the problem in an essential way (through
Proposition~\ref{prop:recall}). Throughout this subsection $G$ denotes a locally compact
totally disconnected group.

\medskip
\noindent {\bf Distortion.} This approach came up
following a question 
by Pierre-Emmanuel Caprace and Shahar Mozes after the completion of the 
original proof.
Recall that if $\Gamma$ is a finitely generated group then an element
$\alpha \in \Gamma$ is called {\it distorted} if w.r.t some (equivalently any)
word length $\ell$ on $\Gamma$  one has: $\ell(\alpha^n) = o(n)$ 
 as $n\to\infty$. The key observation here is:

\begin{proposition}
\label{prop:distor}
Let $\Gamma < \Aut(G)$ be finitely generated. Then the scale of any distorted element 
$\alpha \in \Gamma$ is $1$.
\end{proposition}

For the proof of this proposition, we explain
a basic remark made in \cite{BaWi:direction} without details, 
claiming that the function 
$d : {\mathcal B}(G)\times {\mathcal B}(G) \to {\mathbb R}^+$ defined by
$$
d(U,W) = \log [U:U\cap W][W:U\cap W]
$$
is a metric on ${\mathcal B}(G)$ (the space of compact open subgroups of $G$).
Indeed it is easy to see that this follows immediately from the asymmetric version of the 
triangle inequality, which is valid for all $U,\,V,\,W\in {\mathcal B}(G)$:
\begin{equation}
\label{eq:triangle}
[U:U\cap W] \leq [U : U\cap V][V:V \cap W].
\end{equation}
This inequality (\ref{eq:triangle}) is established by observing that
\begin{equation*}
\label{eq:firststep}
[U:U\cap W] \leq [U : U\cap V \cap W] =  [U : U\cap V][U\cap V : U\cap V\cap W],
\end{equation*}
where the first relation holds because $U\cap V\cap W$ is a subgroup of $U\cap W$, and the second follows from the chain of inclusions $U\cap V\cap W\leq U\cap V\leq U$. To complete the proof of (\ref{eq:triangle}), observe next that the inclusions  $U\cap V\cap W\leq V\cap W$ and $U\cap V \leq V$ induce a map $(U\cap V)/(U\cap V\cap W) \to V/(V\cap W)$ that is easily seen to be injective, from which it follows that
\begin{equation*}
\label{eq:secondstep}
[U\cap V : U\cap V\cap W] \leq [V:V\cap W]. 
\end{equation*}

Next, returning to Proposition~\ref{prop:distor}, 
applying (\ref{eq:triangle}) to the groups $\alpha\beta(V),\,\alpha(V)$ and $V$, where $\alpha,\beta\in \hbox{Aut}(G)$, yields
$$
[\alpha\beta(V):\alpha\beta(V)\cap V]\leq [\alpha\beta(V):\alpha\beta(V)\cap \alpha(V)][\alpha(V):\alpha(V)\cap V],
$$
whence, since $\alpha$ is an automorphism of $G$,
$$
[\alpha\beta(V):\alpha\beta(V)\cap V]\leq [\beta(V) :\beta(V)\cap V] [\alpha(V):\alpha(V)\cap V].
$$

An induction argument then implies the following. 

\begin{lemma}
\label{lem:estimate}
Let $A$ be a finite subset of $\hbox{Aut}(G)$ and let $V\in {\mathcal B}(G)$. Suppose that $\alpha = \gamma_1\cdots \gamma_l$ is a word of length $\ell$ 
in elements of $A\cup A^{-1}$. Then
\begin{equation}
\label{eq:estimate}
[\alpha(V) : \alpha(V)\cap V] \leq M^l,
\end{equation}
where $M = \max\left\{ [\gamma(V) : \gamma(V)\cap V],\,[\gamma^{-1}(V) : \gamma^{-1}(V)\cap V] \mid \gamma\in A\right\}$.
\end{lemma}

\noindent {\bf Proof of Proposition~\ref{prop:distor}.}
Let $V$ be a minimizing subgroup for $\alpha$. Then, by Proposition~\ref{prop:recall}, 
$$
[\alpha^n(V):\alpha^n(V)\cap V] = s(\alpha^n) = s(\alpha)^n\text{ for every }n\geq0.
$$
On the other hand, Lemma~\ref{lem:estimate} shows that there is a constant $M$ such that one has 
$[\alpha^n(V):\alpha^n(V)\cap V] \leq M^{\ell(\alpha^n)}$ 
 for all $n$, which implies $ n\log s(\alpha) \leq \ell(\alpha^n) \log M = o(n)$
as $n\to\infty$,
hence $s(\alpha)=1$.

Now, arguing exactly as in the proof of Theorem~\ref{thm:mainresult} one 
deduces, using the well known result of Lubotzky-Mozes-Raghunathan~\cite{LMR}
which guarantees distortion of unipotent elements in higher rank non-uniform
lattices:

\begin{theorem}
\label{thm:LMZ}
Let $\Gamma<{\bf G}(K)$ be commensurable with $\bf G (\mathcal O)$, 
where the latter is assumed of rank at least $2$, 
and $\bf G$ is assumed $K$-isotropic. 
If  $\bf G (\mathcal O)$ is boundedly
generated by its unipotent elements then $\Gamma$ has the outer 
commensurator-normalizer property.
\end{theorem}

Of course, with some more structural theoretic arguments (along the
same lines of~\cite{LMR}) 
this result  can also be proved using the approach 
taken in the proof of Theorem~\ref{thm:mainresult}.  
The main issue remains finding (unipotently) boundedly
generated arithmetic groups beyond the split and quasi-split cases 
(and the one other recent known family $SO_{n \ge 5}(f)$ for certain quadratic
forms $f$, due to Erovenko-Rapinchuk~\cite{Er-Ra}). Note also that just
being abstractly boundedly generated (by an arbitrary finite set) is {\it not}
enough for our purposes here (unlike for the CSP, for example).

\medskip
\noindent {\bf Divisibility.} Recall that an element 
$\alpha \in \Gamma$ is {\it infinitely divisible} if for every $n$ there is
$m >n$ and $\gamma \in \Gamma$ with $\gamma ^m = \alpha$. We have, see \cite[Proposition 4]{Wi94}:

\begin{proposition}
\label{prop:divisable}
If $\alpha \in \Aut (G)$ is infinitely divisible then $s(\alpha)=1$.
\end{proposition}

\begin{demo} This follows immediately from the inequality
$s(\gamma ^m)= s(\gamma)^m$ in Proposition~\ref{prop:recall}: 
once  $m> s(\alpha)$ we have  $s(\alpha)^{\frac{1}{m}}<2 $ 
and $\alpha$ cannot be an $m$-power unless its scale is $1$.
\end{demo}

\medskip
In the presence of bounded generation, arguing as in 
the proof of Theorem~\ref{thm:mainresult} yields the following.

\begin{corollary}
\label{cor:div-bdd}
Suppose that $\Gamma < \Aut (G)$ 
is boundedly generated by a finite set of cyclic groups, each generated by an infinitely
divisible element of $\Aut (G)$. Then $\Gamma$ stabilizes a compact open subgroup of $G$.
\endproof
\end{corollary}

Although arithmetic groups themselves do not have infinitely divisible
elements, their unipotent elements have this property inside an ambient
$S$-arithmetic group.
This implies a {\it relative outer commensurator-normalizer property}:
under the same assumptions of Theorem~\ref{thm:LMZ}, if $\Lambda$ is
an $S$-arithmetic group with $S$ {\it containing some finite place},
then for any homomorphism $\varphi: \Lambda \to \Delta$, once 
$\Lambda \cap \bf G (\mathcal O)$ commensurates a subgroup of $\Delta$
then it normalizes a subgroup commensurable to it. While this is
of course weaker than  Theorem~\ref{thm:LMZ}, it does 
offer an approach
to cases like $\bf G= SL_2$, $K=\QQ$, which are not handled by 
Theorem~\ref{thm:mainresult} 
(but are covered in Conjecture~\ref{conj:conj} below). For example, 
if answered positively, 
the following natural question and the preceding discussion, combined with the proof
of part {\bf 2} in Theorem~\ref{thm:mainresult}, would show that any 
commensurated subgroup of $SL_2(\QQ)$ is $S$-arithmetic (which is
only conjectured at this point):

\medskip

\noindent {\bf Question.} Do there exist natural numbers $M,L$ such that
every matrix in $SL_2(\ZZ)$ can be written as a product of at most
$M$ elementary matrices in $SL_2(\QQ)$ whose denominator is bounded (in absolute value) by $L$?

\section{Theorem \ref{thm:adelic} and beyond}
\label{sec:adeles}

\subsection{ Proof of Theorem  \ref{thm:adelic} and the case of totally imaginary $K$.} 
\label{sec:Proof1.4}
Consider first the last statement
of the Theorem: assuming that the normalizer of each compact, open
subgroup of $H$ is amenable, we show possibility {\bf 1} cannot occur. 
Indeed, by part {\bf 1} of 
Theorem \ref{thm:mainresult} (in its topological version proved in Section~\ref{sec:Thm1.3Proof}), 
there exists a compact open subgroup $V < H$ normalized by $\varphi 
(\bf G(\mathcal O))$. By assumption then, $\varphi (\bf G(\mathcal O))$ is contained
in an amenable subgroup of $H$, which is not possible if $\varphi 
(\bf G(\mathcal O))$ is discrete and $\ker\varphi$ is finite. 

Next, notice that the two possibilities in the Theorem are mutually exclusive. 
This follows from Theorem~\ref{thm:closedrange}, the fact that
the closure of $\Gamma$ in ${\bf G}(\mathbb {A} _f)$ is a 
direct product
of a compact subgroup and a restricted direct product over
 ${\bf G}(K_\nu)$ (use Proposition~\ref{prop:adelic-open} above), together with the observation
that the image of $\Gamma$ under any infinite quotient of its closure in ${\bf G}(\mathbb {A} _f)$ is 
non-discrete.

We proceed now to the proof of the main part of the theorem, 
and preface it
with some general relevant facts.
Assume that $\Lambda < \Gamma$ are countable, residually finite groups,
and that $\Lambda$ is commensurated by $\Gamma$. 
Then any (separating) family $\mathfrak {F}$
of finite index normal subgroups of $\Lambda$ which is closed 
under intersection,
defines a system of neighborhoods of the identity in a group 
invariant topology 
on $\Lambda$, making the completion of $\Lambda$ 
a compact (profinite) group. Under the obvious conditions on  $\mathfrak {F}$ w.r.t the conjugation action of $\Gamma$, the same topology
will make (the completion of) $\Gamma$ a locally compact totally disconnected topological group, in which the closure
of $\Lambda$ is a compact open subgroup. This situation has already been encountered previously in
the construction of $\Gamma \parallel\Lambda$ (where $\mathfrak {F}$ is the weakest 
topology with this property), and we shall need two more, 
in the first of which $\mathfrak {F}$ 
consists of {\it all} finite index subgroups of $\Lambda$.

\medskip

\noindent {\bf Notation and convention.} The locally compact group
 obtained as the completion of $\Gamma$ in the latter case 
will be denoted $(\Gamma \parallel\Lambda)$.
Whenever this notation is used, it is implicitly assumed 
that $\Lambda$ is 
residually finite and is commensurated by $\Gamma$.

\medskip

Some basic properties, whose easy verification is left to the reader, 
are given by:

\begin{lemma}
\label{lem:||}
{\bf 1.} If $\Lambda < \Gamma < \Delta$ then the inclusion map 
 $\Gamma \to \Delta \to  (\Delta \parallel \Lambda)$ extends to a topological
isomorphism of the group $(\Gamma \parallel \Lambda)$ with the closure of 
 $\Gamma$ in $(\Delta \parallel \Lambda)$.

\noindent{\bf 2.} If $\Lambda _1, \Lambda _2 < \Gamma$ are commensurable
subgroups, then the identity map on $\Gamma$ extends to a topological group
isomorphism of $(\Gamma \parallel \Lambda_1)$ with  $(\Gamma \parallel \Lambda_2)$.
\end{lemma}

The relevance of this construction to our discussion 
comes from the following universality
property, which is opposite to the one appearing in Corollary~\ref{cor:universality}.

\begin{lemma}
\label{lem:universality}
Assume $\Lambda < \Gamma$ is commensurated, and let 
$\varphi : \Gamma \to H$ be a homomorphism into some locally compact 
totally disconnected group $H$. 
If $\varphi (\Lambda)< H$ is bounded, then
$\varphi$ extends to a continuous homomorphism 
$\tilde \varphi : (\Gamma \parallel \Lambda) \to H$.  
\end{lemma}

\begin{demo}
Let $U < H$ be an open subgroup of $H$.
Since  $\varphi (\Lambda)$ is precompact, $U \cap  \varphi (\Lambda)$
has finite index in  $\varphi (\Lambda)$. Hence
$\varphi^{-1}(U \cap  \varphi (\Lambda)) \subseteq 
\varphi^{-1}(U \cap  \varphi (\Gamma)) \subseteq \Gamma$ has finite index in $\varphi^{-1}
(\varphi (\Lambda))$ and the latter contains $\Lambda$.
 Thus $\varphi^{-1}(U)$ contains a finite index subgroup of $\Lambda$ and is
open in $\Gamma$ in the topology of  $(\Gamma \parallel \Lambda)$, as required.
\end{demo}

\medskip

Aiming at Theorem~\ref{thm:adelic}, let $\Gamma < {\bf G} (K)$ be an $S$-arithmetic 
group, and let $\varphi:\Gamma \to H$. Replacing $H$ by the closure
of $\varphi (\Gamma)$, we may assume the latter is dense.
Let $ U < H$ be a compact
open subgroup, and let $A=\varphi (\Gamma) \cap U < H$. As usual, because
$\varphi (\Gamma)$ commensurates $U$ in $H$, it commensurates $A$, hence
$\Lambda :=\varphi^{-1} (A) <\Gamma$ is commensurated by $\Gamma$.
There are now two cases, in accordance with the alternatives in the statement of  
Theorem  \ref{thm:adelic}. If $\Lambda$ is finite, then so is $A$, and by density of $\varphi (\Gamma)$
it follows that $H$ is discrete and so is $\varphi (\Gamma)$. 
Furthermore, the finiteness of $\Lambda$ implies that $\ker\varphi$ is finite
(in which case, as is well known and easy to see, it is central). This
is exactly possibility {\bf 1} of the Theorem. Otherwise,
$\Lambda$ is infinite, and by Theorem
\ref{thm:mainresult} it is standard, hence in particular it contains a finite
index subgroup $\Lambda _0$ of $\bf G(\mathcal O)$. By
 Lemma 
\ref{lem:universality} it follows that $\varphi$ extends 
continuously to 
$(\Gamma \parallel \Lambda_0)$. As $\Gamma <\Delta:= {\bf G}(K)$ 
we may use 
part {\bf 1} of Lemma~\ref{lem:||} to alternatively 
extend $\varphi$ continuously to the closure of $\Gamma$ in
$({\bf G}(K) \parallel \Lambda _0) \cong ({\bf G}(K) \parallel {\bf G}(\mathcal O))$
(the latter isomorphism coming from part {\bf 2} of  Lemma \ref{lem:||}).
In other words, we have just shown:

\begin{proposition}
\label{prop:extend}
Keep the notations above and let $\Gamma$ be as 
in Theorem  \ref{thm:adelic}, and $H$ be an arbitrary l.c.t.d. group. Any
homomorphism $\varphi : \Gamma \to H$ extends continuously to the closure
of $\Gamma$ in $({\bf G}(K) \parallel {\bf G}(\mathcal O))$. 

\end{proposition}

Thus, it remains to identify  the group 
$({\bf G}(K) \parallel {\bf G}(\mathcal O))$, which turns out to be
the essence 
of the
 Congruence Subgroup Problem. 
As the paper is aimed also at readers who are not familiar with the CSP,
we include some further explanations, along the way using 
the opportunity to define
some of the notions and  notation which have already been used, and will be needed
for the second part of this section.
Readers familiar with the CSP are encouraged to skip at this point to that 
part, as for those the rest of the proof of 
Theorem  \ref{thm:adelic} (as well as its extension to totally imaginary $K$) 
should be quite obvious.

Recall that given any set of valuations $V^\infty \subseteq S \subseteq V $ 
the ring $\mathcal O_S \subseteq K$ of $S$-integers in $K$
is defined by:
$$
\mathcal O _S = \{x \in K  |\quad  \nu (x) \ge 0 \quad \forall \nu \notin S \}
$$

Fix a $K$-embedding ${\bf G} \subseteq {\bf GL_n}$ and define
$\Gamma =  {\bf G}(\mathcal O _S) =
{\bf G}(K) \cap {\bf GL_n} (\mathcal O _S)$. For an ideal $\mathfrak {a}$ of $ \mathcal O _S$ the principal $S$-congruence subgroup of level
$\mathfrak {a}$ is defined by $\Gamma _{\mathfrak {a}}:=\Gamma \cap
{\bf GL_n}(\mathcal O _S,\mathfrak {a})$ where ${\bf GL_n}(\mathcal O _S,\mathfrak {a})$ is the subgroup of ${\bf GL_n} (\mathcal O _S)$
consisting of matrices which are congruent to the identity matrix
mod $\mathfrak {a}$. A subgroup $\Gamma ' < \Gamma$ is said to be
\emph{$S$-congruence} if it contains some principal  $S$-congruence subgroup
(in which case it is obviously of finite index). Of course, when
$S=V^\infty$ one is reduced to the case of ordinary congruence
subgroups of the arithmetic group  ${\bf G}(\mathcal O)$, and in its
classical formulation, the congruence subgroup problem asks whether
any finite index subgroup of it is congruence. The question
makes sense just as well for every $S$-arithmetic group 
$\Gamma <{\bf G}(K)$, working with $S$-congruence subgroups.

A quantitative manner to formulate the CSP is to
notice first that as with all finite index subgroups,
the congruence ones also 
define a (profinite) topology on  $ {\bf G}(\mathcal O)$  
called the congruence topology  
(which does not depend on the $K$-embedding ${\bf G} \subseteq  {\bf GL_n} $).
As with the profinite topology, we may induce the latter topology on 
every subgroup 
$ \Gamma < {\bf G}(K)$, and
 the completion is again 
locally compact and totally disconnected. We shall denote it here by 
(the {\it ad hoc\/} notation) $(\Gamma \parallel^c {\bf G}(\mathcal O))$.
Part {\bf 1} of Lemma \ref{lem:||} applies here just as well, and in particular identifies  $( {\bf G}(\mathcal O)\parallel^c  {\bf G}(\mathcal O))$ 
 with the closure
of  $ {\bf G}(\mathcal O)$ in $({\bf G}(K) \parallel^c {\bf G}(\mathcal O))$.
Returning to the CSP, one obviously has 
a continuous homomorphism: 
$\pi: ({\bf G}(\mathcal O)\parallel  {\bf G}(\mathcal O)) 
\longrightarrow  ({\bf G}(\mathcal O)\parallel^c  {\bf G}(\mathcal O))$
whose kernel $C(\bf {G})$ is called the \emph{Congruence Kernel}. 
The homomorphism $\pi$ is easily seen to be onto 
(the image being dense and compact), and
is a topological group isomorphism  if and only if every finite index subgroup is
congruence (which is the ``classical'' version of the CSP). 
The ``size'' (or more precisely, finiteness property) 
of the congruence kernel, when it is non-trivial,  
may be regarded as measuring the (failure of the) CSP, and following Serre, one
says that the congruence subgroup property holds 
(for ${\bf G}(\mathcal O)$) if  $C({\bf G})$ is finite. A similar
definition and discussion applies to the $S$-arithmetic subgroups, thus defining
the $S$-congruence kernel $C(S,{\bf G})$ using the topology defined by the 
$S$-congruence subgroups.

We can turn now to the proof of 
Theorem \ref{thm:adelic}. Note that by 
Lemma \ref{lem:universality} one has a natural continuous  homomorphism $p$:

\begin{equation}
\label{eq:C}
 1\longrightarrow C   \longrightarrow   ({\bf G}(K)\parallel  {\bf G}(\mathcal O)) 
  \stackrel{p }  \longrightarrow  ({\bf G}(K)\parallel^c  {\bf G}(\mathcal O))
\end{equation}
which is onto by density of the image and $p$ 
being open. 
It is easy to verify directly from the definition
that ker $p= C= C({\bf G})=$ ker $\pi$ (and not just  $ C({\bf G})\subseteq C $),
with $\pi$ defined as above,
which explains
why the computation of the profinite completion 
of ${\bf G}(\mathcal O)$,
and that of the group $({\bf G}(K)\parallel  {\bf G}(\mathcal O))$, 
are essentially the same problem.
To be concrete, it is straightforward to see that the 
congruence completion  $({\bf G}(\mathcal O)\parallel^c  {\bf G}(\mathcal O))$ is
isomorphic to the closure of  ${\bf G}(\mathcal O)$ under its diagonal
embedding in ${\bf G}(\mathbb {A} _f)$, and likewise one may naturally identify 
$({\bf G}(K)\parallel^c  {\bf G}(\mathcal O))$ with the  
closure of ${\bf G}(K)$ in 
${\bf G}(\mathbb {A} _f)$. As ${\bf G}$ is simply connected 
and $K$-isotropic, by strong approximation the  latter
is all of ${\bf G}(\mathbb {A} _f)$, and denoting 
$\hat G = ({\bf G}(K)\parallel  {\bf G}(\mathcal O))$ makes~(\ref{eq:C}) into:
\begin{equation}
\label{eq:CC}
1\longrightarrow C({\bf G}) \longrightarrow 
\hat G  \stackrel{p} \longrightarrow {\bf G}(\mathbb {A} _f)  
\end{equation}

\medskip
Observe that by construction, the extension $p$ splits over 
${\bf G}(K)$, that is, we have a canonical embedding
${\bf G}(K) \to \hat G$ whose composition with $p$ is the diagonal
embedding.
To finish the proof of Theorem~\ref{thm:adelic}, all that remains now is to quote 
the well known solution to the 
CSP in the (split) cases covered by the Theorem (in its final form due to 
Serre~\cite{Se} for ${\bf SL_2}$ and Matsumoto~\cite{Ma} in the higher rank case): 
When $K$ admits a real embedding $ C({\bf G})=1$, i.e. $\hat G \cong {\bf G}(\mathbb {A} _f)$, 
which together with 
Proposition \ref{prop:extend} proves
Theorem \ref{thm:adelic}. 

\medskip
Let us discuss now Theorem~\ref{thm:adelic} when $K$ is totally imaginary. 
In this case the solution to the CSP says that
$C \cong \mu(K)$ is the finite cyclic group of roots of unity in $K$, 
and the extension (\ref{eq:CC}) is central. Thus, arguing as in the 
proof of  Theorem~\ref{thm:adelic}, the CSP in this case 
together with Theorem 
\ref{thm:mainresult} yield that any $\varphi :\Gamma \to H$
as in Theorem  \ref{thm:adelic} extends  continuously to the closure 
of $\Gamma$ in the
locally compact group $\hat G$, in which $\Gamma$ embeds via the 
embedding of all of
${\bf G}(K)$. In other words, a similar theorem holds with ${\bf G}(\mathbb {A} _f)$ replaced 
by its finite extension $\hat G$. Finally, note that upon taking say, $\Gamma ={\bf G}(K)$, 
 the case $H=\hat G$ 
is itself covered by
Theorem~\ref{thm:adelic}, hence the original statement of this Theorem
indeed must be modified when $K$ is totally imaginary, and the
``envelop'' of $\Gamma$ must accommodate this group. 
Hence this modified version of Theorem~\ref{thm:adelic}, when $K$ is totally imaginary, 
is as sharp as it can be, bearing in mind also the following remark.

\begin{remark}
\label{rem:regularity}
The only situations in the split case left out by the above general 
version of  Theorem \ref{thm:adelic} are when ${\bf G} = {\bf SL_2}$ and
$K$ is either $\QQ$ or a quadratic imaginary extension of it. 
In these cases Theorem \ref{thm:adelic} (in its above full version) {\it always fails}. 
(In particular, we witness a rare situation
in rigidity theory with a bold difference between 
$SL_{n\ge 3}(\QQ)$ and $SL_{2}(\QQ)$.) This results from the failure, in these cases, 
of the CSP 
for $SL_2(\mathcal O)$, taking $H= ({\bf G}(K)\parallel  {\bf G}(\mathcal O))$
in (\ref{eq:C}) above. It is a classical observation of Serre that in 
the latter case the congruence
kernel $C$ cannot be topologically finitely generated (for the 
much stronger, so-called ``Rapinchuk Lemma'' see~\cite[p. 138]{Lu-Se}; see 
also the appearance
of this feature for a similar reason in Conjecture \ref{conj:conj} below).
Conjecture  \ref{conj:conj}  entails the expectation 
that a modified version of Theorem~\ref{thm:adelic} (partially motivated by Serre's observation) 
should still hold 
when $\Gamma$ has higher rank. A positive answer to 
the Question posed at the end of Section~\ref{sec:Thm1.3Proof} would enable one to
 approach 
this issue.
\end{remark}

\subsection{ Beyond Theorem \ref{thm:adelic}.} 
\label{sec:beyond}
Although the CSP quoted in the proof of Theorem \ref{thm:adelic} is well known by now, 
it reveals what may 
seem a somewhat curious dependence of the situation on the the global 
field $K$. 
In this section we suggest a conjectural extension of 
Theorem \ref{thm:adelic} which 
accounts for this situation and, more significantly, predicts
sharp forms of different known results in 
the theory of arithmetic groups, placed in a unified manner.

\medskip

\noindent {\bf The setting.} Let $K$ be any number field, $\mathbb {A}$ its ring of adeles,
 and unless assumed specifically otherwise,
 we let now 
${\bf G}$ be \emph{any} absolutely simple, simply connected algebraic group defined over $K$.
Our starting point is a deep work of 
Deligne~\cite{De1}, where using Galois cohomology 
canonical
central extensions of  ${\bf G}$ are constructed, over local and global fields. 
For our purposes, of particular interest is 
the central extension denoted 
${\bf G}(\mathbb {A})^{\widetilde{ }}$  of ${\bf G}(\mathbb {A})$ by
  the finite group $\mu$ of roots of unity in $K$, which splits over ${\bf G}(K)$. In other
words, Deligne constructs the following central extension: 
\begin{diagram}\label{adiagram}
    &        &         &         &                                          &{\bf G}(K)       &                     & \\
    &        &         &         &  \ldTo(1,2)^{s}                        &                & \rdTo(1,2)^{d}    & \\
 1 & \rTo & \mu & \rTo &{\bf G}(\mathbb {A})^{\widetilde{ }} & \rTo^{p} &{\bf G}(\mathbb {A}) & \rTo & 1
\end{diagram}
where $d$ is the standard diagonal embedding of ${\bf G}(K)$ in ${\bf G}(\mathbb {A})$, and $s$ is the splitting,
that is, $p \circ s = d$. The first appearance of a central extension of this type was in the
seminal work of Weil~\cite{We} on the metaplectic kernel, with ${\bf G = Sp_{2n}}$ (where the splitting over 
${\bf G}(K)$, proved analytically by Weil, translates to the quadratic reciprocity law for Hilbert's 
2-symbol). This theme was then studied thoroughly by Moore in his fundamental
paper~\cite{Mo}, at about the same time its relation to the CSP was discovered (see Bass-Milnor-Serre~\cite{BMS}). 
Among other things, Moore showed that if 
$\bf G(K)$ is perfect (now known to be the case
for all $K$-isotropic ${\bf G}$ by the proof of the Kneser-Tits 
conjecture), then there exists a {\it universal} 
central extension $p:{\bf G}(\mathbb {A})^{\widetilde{ }} \to 
{\bf G}(\mathbb {A})$ which splits over ${\bf G}(K)$ as in~(\ref{adiagram})  
above, i.e., any other such extension
is a ``pushout'' (an epimorphic image in the category of central extension) of the universal one.
Using their deep results on the finiteness of the metaplectic kernel (which supply an
``upper bound'') combined with properties of Deligne's construction~(\ref{adiagram}) (supplying a ``lower bound''), Prasad and Rapinchuk have shown that Deligne's construction is indeed Moore's universal 
central extension 
(see Section 8 of their fundamental paper~\cite{PR1}). 
In fact, modulo a mild conjecture (``Conjecture U'' in~\cite{PR1}) needed to 
handle a few exceptional cases, the universality property always holds.

\medskip
\noindent {\bf The Conjecture.} 
To state the conjecture we shall need the
following:

\begin{proposition}
\label{minimal}
Retain the notation from the previous section, including that
in Deligne's central extension~(\ref{adiagram}) above. Let 
$\Gamma < {\bf G}(K)$ be an infinite commensurated (e.g. $S$-arithmetic)
subgroup. Then there exists
a minimal (in the strong sense) {\it open} subgroup of 
${\bf G}(\mathbb {A})^{\widetilde{ }}$ which contains the subgroup
$s(\Gamma)$ (where $s$ is the splitting map of~(\ref{adiagram})). 
\end{proposition}

\begin{notation}\label{ambient}
Given an infinite commensurated subgroup $\Gamma < {\bf G}(K)$, 
the unique minimal open subgroup of ${\bf G}(\mathbb {A})^{\widetilde{ }}$ 
containing $s(\Gamma)$ will be denoted by $\tilde G_\Gamma$.
\end{notation}

In order to proceed to the conjecture, the proof of the Proposition, 
which relies on strong approximation, is deferred to the next subsection.
Abusing notation, by means of the splitting map $s$ in~(\ref{adiagram}) 
we shall freely regard $\Gamma$ as a subgroup
of $\tilde G_{\Gamma}$ as well.

\smallskip

\begin{conjecture}
\label{conj:conj}
Retain the notations above. 
Assume that $\Gamma < {\bf G}(K)$ is a 
commensurated (e.g. $S$-arithmetic) subgroup 
satisfying:  $\, \mathbb {A}$-rank$(\Gamma):= \Sigma K_{\nu}$-rank ${\bf G}(K_{\nu}) 
\ge 2$, where the sum is taken over all those $\nu$ for which
the projection of $\Gamma$ to ${\bf G}(K_\nu)$ is unbounded (note that 
this agrees with the usual notion of $S$-rank when $\Gamma$ is $S$-arithmetic). 
Let $H$ be {\it any} locally compact group
and $\varphi : \Gamma \to H$ a homomorphism. Assume further the following 
{\it regularity} condition:  modulo the 
center of $H$, the subgroup
$Q<H$ generated by all the compact subgroups of $H$ which are normalized 
by $\varphi (\Gamma)$ is {\it tame} (a notion defined below).
Then in each 
one of the following two cases $\varphi$
 extends continuously to a homomorphism 
$\tilde \varphi : \tilde G _{\Gamma} \to H$ with closed image:

\smallskip
\noindent {\bf 1.} The homomorphism $\varphi$ is not proper (for properness see 
{\bf 1} of Theorem  
\ref{thm:adelic}).

\noindent {\bf 2.} Denoting by $H_0$ the connected component of $H$, 
the normalizer of any compact open subgroup of the totally disconnected
group $H/H_0$ is compact.

\end{conjecture}

\medskip

For the ad hoc purpose of this conjecture, a topological group $Q$ is
called {\it tame} if its connected component $Q_0$ is a compact finite 
dimensional (i.e. Lie) group, and the group $Q/Q_0$ is a topologically
finitely generated profinite (compact) group. Note that as is evident from 
Remark~\ref{rem:regularity} above, a 
regularity condition is 
essential for the conjecture in the cases where $\Gamma$ admits an arithmetic subgroup for which 
the CSP fails. Its relevance, already mentioned in that Remark, is related to the 
compactness of the automorphism group of $Q$ (in the connected case it guarantees
finite dimensionality of the whole group by structure theory).
However, we shall not dwell on this issue further.
As the proof of 
Theorem~\ref{thm:adelic} shows, it is to be expected that in the other cases (e.g. when
${\bf G}({\mathcal O})$ has rank at least 2), the regularity condition
in Conjecture \ref{conj:conj} may be relaxed.

\medskip

\noindent {\bf Proof of Proposition~7.5 and preparations towards unfolding Conjecture~\ref{conj:conj}.} 

\smallskip

\noindent{\bf Proof of Proposition~7.5.} Let $C$ (resp. $NC$) be the set of places $\nu \in V$ for which $\Gamma$ projects 
boundedly (resp. unboundedly) to ${\bf G}(K_\nu)$. Let $G^{(c)}_{\Gamma}$ be the
(compact) closure of the projection of $\Gamma$ to $\Pi_{\nu \in C} {\bf G}(K_\nu)$, and 
let $G^{(nc)}_{\Gamma}=\Pi_{\nu \in NC} {\bf G}(K_\nu)$ (restricted product). Define 
$G_{\Gamma}: = G^{(nc)}_{\Gamma} \times G^{(c)}_{\Gamma} <  {\bf G}(\mathbb {A})$.
Any open subgroup $\tilde U < {\bf G}(\mathbb {A})^{\widetilde{ }}$ containing $s(\Gamma)$ 
must project under $p$ of~(\ref{adiagram}) to an open subgroup $p(\tilde U)$ 
of ${\bf G}(\mathbb {A})$ which contains $\Gamma$. 
By Proposition~\ref{prop:adelic-open} (and the fact that non-compact 
${\bf G}(K_{\nu})$ have no proper 
open unbounded subgroups),
this implies  that $p(\tilde U)$ must contain $G_{\Gamma}$.
The key point in the proof is that under our 
assumption on $\Gamma$ the latter is in fact {\it open} in ${\bf G}(\mathbb {A})$. Assuming
this for the moment, we see that $p^{-1}(G_{\Gamma}) < 
{\bf G}(\mathbb {A})^{\widetilde{ }}$  is one open subgroup containing $s(\Gamma)$, hence
by intersecting it we may restrict the minimality claim only for its open subgroups.
However, we now have that any open subgroup $\tilde U < p^{-1}(G_{\Gamma})$ 
containing  $s(\Gamma)$
must project {\it onto} $G_{\Gamma}$ under $p$, hence it is determined by the kernel of its
projection, which is a subgroup of the {\it finite} group $\mu$. Thus, 
clearly there is a 
minimal one.

We are thus left with establishing the openess claim on $G_{\Gamma} < {\bf G}(\mathbb {A})$
which, by definition of $G_{\Gamma}$, is exactly the openess of $G^{(c)}_{\Gamma}$
inside $\Pi_{\nu \in C} {\bf G}(K_\nu)$. In the case where $\Gamma$
is $S$-arithmetic this is an immediate consequence of strong approximation
(as $\Gamma$ is assumed infinite, $NC$ is non-empty). In the general case 
where we only assume
that $\Gamma$ is commensurated by ${\bf G}(K)$, one applies Weisfeiler's strong
approximation theorem in~\cite[Theorem 10.5]{Weis} (noticing
that by (ii) there, the set $V-S$ in (i) must contain all places $\nu $ for which the
projection of $\Gamma$ to ${\bf G}(K_\nu)$ is bounded. The condition that 
the traces of Ad$\{\gamma \}_{\gamma \in \Gamma}$ generate $K$ follows 
essentially from Vinberg's 
theorem -- cf.~\cite[Lemma 2.6]{Pra-Rap}). This completes the proof of Proposition 7.5.

\smallskip

We note that conjecturally, the assumption that $\Gamma$ is commensurated 
by ${\bf G}(K)$  is satisfied exclusively 
by the $S$-arithmetic subgroups; as discussed below, this is actually  
predicted by the
 Conjecture itself (hence making the a priori stronger arithmeticity assumption on $\Gamma$ 
results in an equivalent statement). The advantage of this formulation, however,
lies in its being purely group theoretic.
To appreciate better what this sweeping conjecture entails,
 we shall need to recall first some additional information concerning  
Deligne's central 
extension~(\ref{adiagram}), which was not needed for the proof of 
Proposition~7.5.

Although Deligne's construction is {\it global}, 
for every place $\nu \in V$ the extension  
induces a central extension of the group ${\bf G}(K_{\nu})$ by the 
group $\mu$ of roots of unity in $K$, and
one has
good understanding of its local nature, particularly over the 
$K_{\nu}$-isotropic 
places (which will be more relevant to us).
For every such place $\nu$ which is \emph{non-archimedean},  
the extension $p$ over  ${\bf G}(K_{\nu})$  defines
an element of (maximal) order $\#\mu$ in 
$H^2({\bf G}(K_{\nu}), \mu)  \to H^2({\bf G}(K_{\nu}), \RR/\ZZ)$. 
Thus the group 
$p^{-1}({\bf G}(K_{\nu}))$ is perfect, and like ${\bf G}(K_{\nu})$
it has no proper open unbounded subgroups, and is projectively simple 
with $\mu$ sitting in its center. 
When $K = \CC$ the group ${\bf G}(K_{\nu})= 
{\bf G}(\CC)$ is topologically 
simply connected, hence $p$
necessarily splits over ${\bf G}(K_{\nu})$, whereas 
when $K_{\nu}= \RR$ we have  $\mu = \pm 1$
(as $K$ embeds in $\RR$) and either ${\bf G}(K_{\nu}) $ is topologically 
simply connected,
in which case again $p$ must split, or, when it isn't, $p^{-1}
({\bf G}(K_{\nu}))$ is the unique two sheeted covering of 
 $ {\bf G}(K_{\nu}) $ (note that the latter is always topologically connected,
as {\bf G} is simply connected. 
As in the non-archimedean case, a non trivial cover always occurs 
when $\bf G$ is Chevalley and is
never algebraic.)  
This discussion accounts for part {\bf 1} in the following
result, which collects information of importance when studying
implications of Conjecture~\ref{conj:conj}.

\begin{proposition}
\label{prop:split}
Retain the notation in and before Conjecture \ref{conj:conj}, including
 the one in Deligne's central extension~(\ref{adiagram}), which restricts
to an  extension
$p:\tilde G _{\Gamma} \to G _{\Gamma}$ of the groups defined in Proposition 7.5.
Let $\Gamma < {\bf G}(K)$ be infinite $S$-arithmetic.

\smallskip
\noindent {\bf 1.} Let $\nu \in S$ be such that
$\bf G$ is $K_ {\nu}$-isotropic (equivalently,
the projection of $\Gamma$ to $ G_{\nu}:={\bf G}(K_ {\nu})$ is unbounded). 
Let  $\tilde \varphi$ be a continuous group homomorphism defined on  
$\tilde G _{\Gamma}$, whose
restriction to $p^{-1}( G_{\nu})< \tilde G _{\Gamma} $ 
has {\it infinite kernel}.
If $\nu$ is non-archimedean, 
or $\nu$ is real and 
${\bf G}(K_ {\nu})$ is not simply connected, 
then necessarily $\mu \subset $Ker $\tilde \varphi$. 
Hence, by projective simplicity of $p^{-1}( G_{\nu})$, $\tilde \varphi$ 
must factor through 
a group homomorphism $ \tilde G _{\Gamma} \to   G_{\Gamma} \to G_{\Gamma}^{(\nu)}$
where the latter denotes the canonical complement of $G_{\nu}$ 
in  $G_{\Gamma}$.

\smallskip
\noindent{\bf 2.} Consider the canonical product splitting
$G _{\Gamma}= G^{(sp)}_{\Gamma} \times G^{(nsp)}_{\Gamma}$ where 
$ G^{(sp)}_{\Gamma}$ denotes the (finite, possibly empty) product over all 
archimedean $\nu$ which are either 
complex, or real with ${\bf G}(K_ {\nu})$ topologically simply connected. 
Then the extension~(\ref{adiagram}) splits over $ G^{(sp)}_{\Gamma}$, and hence
(using Lemma 1.4 in~\cite{PR1}), there is a splitting 
$\tilde G _{\Gamma} \cong  G^{(sp)}_{\Gamma} \times \tilde G^{(nsp)}_{\Gamma}$.

\smallskip
\noindent{\bf 3.} After replacing $\Gamma$ with a finite index subgroup 
$\Gamma '$, the extension~(\ref{adiagram}) always splits over the compact 
subgroup $G^{(c)}_{\Gamma '}$ defined in the proof of Proposition 7.5.
Therefore, again  (using~\cite[Lemma 1.4]{PR1}) there is
a splitting $\tilde G _{\Gamma '} \cong  G^{(c)}_{\Gamma '} \times 
 \tilde G^{(nc)}_{\Gamma '}$.

\end{proposition}

The first two
parts are self-explanatory in view of the preceding paragraph. 
To prove part {\bf 3} 
separate first $ G^{(c)}_{\Gamma}$ into the archimedean
and non-archimedean components. 
Recall that when ${\bf G}(\RR)$ is compact it is always simply connected, which follows
from Weyl's celebrated result that $\pi _1$ of simple compact Lie 
groups is always finite, and the real algebraicity of
compact simple Lie groups. Hence we only need to deal with the non-archimedean factor
which is where passing to a finite index subgroup appears.
Indeed, every extension
of a profinite group by a finite group $\mu$ is again profinite, 
hence residually finite, and 
so a finite index subgroup of the extension intersects $\mu$ trivially. 

\medskip
\noindent {\bf Implications of Conjecture~\ref{conj:conj}.} With Proposition~7.8 at hand, 
we may turn to discuss various aspects of
   Conjecture \ref{conj:conj}, some concrete 
predictions it makes, and their compatability with
existing results in the theory of arithmetic groups.

\medskip

\noindent {\bf The CSP.} As we now explain, 
Conjecture~\ref{conj:conj} for $S$-arithmetic $\Gamma$ and finite $H$
 is a quantitative unified
form of the CSP 
(conjectural in general, yet known in many cases): {\it Every homomorphism of
$\Gamma$ to a finite group extends continuously to $\tilde G _{\Gamma}$}.
Note that unlike standard approaches to the CSP which aim
to describe directly (or through computation of the congruence or 
metaplectic kernel) 
the profinite 
completion of $\Gamma$, the group
$\tilde G_\Gamma$ is not compact (and it has non-trivial connected component), 
and $\Gamma < \tilde G_\Gamma$ embeds discretely.

To unfold the above abstract statement, consider first the case $\Gamma = 
{\bf G}(\mathcal O)$.
Any continuous homomorphism
$\tilde \varphi: \tilde G _{\Gamma} \to H$ must vanish on the connected
component $\tilde G _{\Gamma}^0$, and  
 $\Gamma$ projects densely to 
the profinite group $ \tilde G _{\Gamma} / \tilde G _{\Gamma}^0$ 
(by the minimality property of $\tilde G _{\Gamma}$).  
Hence,  Conjecture \ref{conj:conj} predicts
that the latter is the full profinite completion of $\Gamma$. 
As the connected component arises precisely from the archimedean places,  
Proposition~7.8 shows that when there is a real place $\nu$
with ${\bf G}(K_ {\nu})$ not simply connected (case {\bf 1}), 
all finite quotients
of $\Gamma$ must factor through a homomorphism to 
the closure of $\Gamma$ in $ {\bf G}(\mathbb {A}_f)$,
that is, {\it every finite index subgroup is congruence}. 
Otherwise, we are in the setting of  part {\bf 2} of the Proposition, 
where the connected component $\tilde G _{\Gamma}^0$ is 
precisely the group  $G^{(sp)}_{\Gamma}$ there.  Let us see that
in this case the conjecture predicts that the congruence kernel
$C(G)$ is isomorphic to $\mu$ (and is not merely contained in $\mu$, as evidently follows 
from the conjecture). 
We have the following diagram
of exact sequences:

$$\begin{array}{ccccccccc} 1 & \to &  C(G) & \longrightarrow &
  ({\bf G}(K)\parallel  {\bf G}(\mathcal O))  &
\longrightarrow & {\bf G}(\mathbb {A}_f) \cong {\bf G}(\mathbb {A})/\text{connect comp}  & \to & 1 \\ & & \downarrow & & \downarrow & &\parallel & & \\
1 & \to & \mu & \longrightarrow &  {\bf G}(\mathbb {A})^{\widetilde{ }}/\text{connect comp}  &  \longrightarrow & 
{\bf G}(\mathbb {A})/\text{connect comp} \cong {\bf G}(\mathbb {A}_f) & \to
& 1
\end{array}$$

The upper row is the one defining the congruence kernel (see~(\ref{eq:C})
above),
whereas the lower uses the fact that the
connected component is simply connected.
The non-trivial arrow down in the middle 
follows from the universality 
property in Lemma~\ref{lem:universality}
above.
Thus, if $C(G)$ was (mapped to) a proper subgroup of $\mu$, restricting the lower central 
extension
to any of the non-archimedean $\nu$ for which $\bf G$ is $K_\nu$-isotropic would
contradict the known property of (the order of) Deligne's central 
extension as an element of $H^2$ (see the discussion at the beginning of the
paragraph preceding Proposition~7.8). It then follows that $C(G)=\mu$
and that in this case
the central extension obtained 
from~(\ref{adiagram}) 
where both groups are divided by their connected (simply connected) component
is isomorphic to the one 
in~(\ref{eq:CC}) defining
the congruence kernel. 

Assume now that $\Gamma = 
{\bf G}(\mathcal O _S)$ and $S$ contains a 
non-archimedean $\nu$ for which $\bf G$ is $K_\nu$-isotropic. Apply again part {\bf 1} of 
the Proposition 
to this $\nu$ 
and a homomorphism $\varphi$ of $\tilde G _{\Gamma}$ ranging in 
a finite group $H$. It follows that $\mu \subset \text{Ker} \varphi$,  
and as in the first part of the 
analysis above we get that the $S$-congruence kernel is expected to trivialize.

It turns out that the above predictions of Conjecture~\ref{conj:conj} 
are fully compatible with the  deep known results on the CSP, 
which arise as accumulation of work of a long list of distinguished contributors
(cf. the beginning of page 304 in Raghunathan's account~\cite{Rag1} which
contains also very comprehensive bibliography; see also
 the excellent recent survey~\cite{PR2} by Prasad-Rapinchuk,
which is particularly recommended to a non-specialist reader of this 
paper). 
As recalled next, when there is a non-archimedean $K_\nu$-anisotropic place 
$\nu \in S$, Serre's original conjecture for $\Gamma = 
{\bf G}(\mathcal O _S)$ fails. However, Conjecture~\ref{conj:conj} suggests a
natural way to remedy this ``irregularity''.

\smallskip

\noindent{\bf The case of anisotropic $\nu \in S$ and the Margulis-Platonov conjecture.} 
As first observed by Raghunathan,
 when there is a non-archimedean 
$\nu \in S$ for which ${\bf G}$ is $K_{\nu}$-anisotropic the $S$-congruence
kernel must be infinite and non-central.  This is best illustrated when taking $S=V$, 
observing that by definition ${\bf G}(K)$ has {\it no} congruence subgroups
(as $K$ has no non trivial ideals), yet it does contain the finite index 
subgroups induced from its embedding in the profinite group ${\bf G}(K_{\nu})$.
This ``counter-example'' turns out to be the most important one, and is the subject of the  
Margulis-Platonov conjecture ({\bf MP}), which can be stated as the claim 
that
every homomorphism of 
$\Gamma = {\bf G}(K)$ into a finite group extends continuously to  
$G _{\Gamma}$ (cf. the proof of Proposition 7.5), and then, as it must, 
it factors through a homomorphism
$G _{\Gamma} \to \Pi \, {\bf G}(K_{\nu})$, where the product is taken over the
finitely many non-archimedean $\nu$ for which  ${\bf G}$ is 
$K_{\nu}$-anisotropic. Notice that
this same prediction is made by
  Conjecture \ref{conj:conj}, as follows from part {\bf 1} of Proposition \ref{prop:split} 
applied to any one of
the almost all non-archimedean $\nu \in V$ where ${\bf G}$ is 
$K_{\nu}$-isotropic (a similar situation is expected when $S$ is 
co-finite, which is consistent with the general CSP result of 
Prasad-Rapinchuk~\cite[Section 9]{PR1} conditioned on {\bf (MP)}. 
See also Rapinchuk's~\cite{Rapi}
and the references therein for the latest progress on this conjecture.)
 Thus, Conjecture \ref{conj:conj}  
 removes the restriction in Serre's conjecture by 
offering a unified statement. In fact, the same viewpoint 
 suggests a variation on the definition of the congruence kernel, which
agrees with the latter in the isotropic case and might be 
called the {\it adelic kernel}
$\mathbb {A} (\Gamma, {\bf G })$. It makes sense for {\it every} $S$-arithmetic
group $\Gamma$ and is defined as the kernel  of the natural epimorphism from
the full profinite completion of $\Gamma$ to the profinite completion w.r.t
all finite index ``adelic'' normal subgroups, i.e. those   $\Gamma _0 <\Gamma$ 
for which the 
homomorphism $\Gamma \to \Gamma / \Gamma _0$ extends to 
the group $G _\Gamma<{\bf G}(\mathbb {A})$. 
The same arguments as in the foregoing discussion 
around the CSP show that  Conjecture \ref{conj:conj} predicts the 
following uniform statement
{\it for all} (higher rank) $\Gamma$ and $\bf G$: 

\smallskip
{\centerline {$\mathbb {A} (\Gamma, {\bf G })$ 
is central, and is a subgroup of $\mu$.}}

\smallskip
Moreover, $\mathbb {A} (\Gamma, {\bf G })$ should always be trivial when $S$ contains a real place $\nu$ with ${\bf G}(K_\nu)$ 
non simply connected, or a non-archimedean $\nu$ where ${\bf G}$ is $K_\nu$-isotropic.
For $\Gamma = {\bf G}(\mathcal O _S)$ 
it should be either all of $\mu$ or trivial, with the former
occurring precisely when the connected component $G _{\Gamma}^0< G _\Gamma$ is topologically
simply connected, and the totally disconnected quotient 
$ G _{\Gamma} / G _{\Gamma}^0$ is compact. (In this case all subgroups of $\mu$ also occur, for
suitable finite index subgroups of $\Gamma$.)
It seems that these sharp predictions concerning the CSP have not appeared explicitly 
in the literature.

\medskip

\noindent {\bf Margulis' normal subgroup theorem.} The Margulis' NST for
an  $S$-arithmetic group $\Gamma$ is {\it equivalent} to  Conjecture \ref{conj:conj}
in the case where $H$ is a discrete group 
with no non-trivial finite normal 
subgroups and $\varphi$ is surjective. To show NST $\Longrightarrow$ Conjecture,   
given a surjective $\varphi$ we may assume it has infinite kernel 
(otherwise there is nothing to prove). 
The NST implies that  Im$\varphi =H$
is finite, and by the additional condition on $H$ it is actually trivial 
(hence $\varphi$ trivially extends).
In the other direction, aiming at a contradiction, 
let $N<\Gamma$ be an infinite normal subgroup of infinite index. Let 
$L=\Gamma / N$ be the infinite quotient. It is convenient to assume,
as in Margulis' original NST, that $\Gamma$ and hence also $L$,  is finitely 
generated. This is not essential for the argument (see below), 
however it immediately implies that the infinite group $L$ has an {\it infinite} 
quotient $H$ with no non trivial finite normal subgroups (one can use, or
argue similarly to, the result that any infinite f.g. group has an infinite 
{\it just infinite} quotient). Applying case {\bf 1} of the conjecture, as we may for this $H$
and the quotient map $\varphi:\Gamma \to H$, yields a continuous
extension $\tilde \varphi$ defined on $\tilde G _{\Gamma}$ and ranging in a 
discrete group $H$. It is easy to see from the structure of $\tilde G _{\Gamma}$ that any such $\tilde \varphi$
must factor through a compact quotient, hence its image is finite, in 
contradiction to our assumption (on $L$ and $H$). With some 
more care  one can also deal with the case where $S$ is infinite -- see the following
discussion on the Margulis-Zimmer CmSP, which generalizes the current one.

\medskip 

\noindent {\bf The Margulis-Zimmer CmSP.} Extending the previous discussion
we have the following result, which clarifies the
precise relation between the CmSP and Conjecture~\ref{conj:conj}:

\begin{proposition}
\label{prop:CmSP-Conj}
The CmSP is equivalent to Conjecture \ref{conj:conj}
 when $H$ is totally disconnected with no non-trivial 
compact normal subgroups, and
 Im$\varphi$ is dense.
\end{proposition}


\smallskip
Note that Proposition~\ref{prop:CmSP-Conj} shows that the CmSP may be 
viewed as the exact natural complement to the CSP in the description of 
homomorphisms of $S$-arithmetic groups $\Gamma$ into general totally 
disconnected groups $H$. Indeed, if $G _{\Gamma}= G^{(c)}_{\Gamma} \times G^{(nc)}_{\Gamma}$
is the splitting as in the proof of Proposition~7.5, then for the CSP the 
 right factor vanishes under the homomorphisms of interest, while it is the opposite situation
for the CmSP by
Proposition \ref{prop:CmSP-Conj}.
Since we shall not make further use of Proposition \ref{prop:CmSP-Conj} beyond
the preceding discussion, we shall leave out some details in its proof.

\medskip
\begin{demo} 
In one direction assume the CmSP and let 
$\varphi:\Gamma \to H$ have dense image. If Im$\varphi <H$ is discrete  
then so is $H$ and $\varphi$ is surjective, hence, as the CmSP implies NST we 
land in the setting of the previously discussed NST. Thus, we may assume that  
$H$ is not discrete and Im$\varphi$ is dense.
Let $\Lambda = 
\varphi^{-1}(\varphi(\Gamma) \cap V)$ where $V < H$ is compact open (and infinite).
Then $\Lambda <\Gamma$ is an infinite commensurated subgroup, which we may assume to be 
of infinite index, and by 
the CmSP it is $S'$-arithmetic for some $V^{\infty}\subseteq S' \subset S$. 
By strong approximation
the image of $\Gamma$ under its diagonal embedding in
$G={\bf G}(\adele ^{S\setminus S'})$ (notation as before Proposition~\ref{prop:G(A)}) is dense, and 
the closure of $\Lambda$ in $G$ is compact open.  Thus, by Corollary \ref{cor:universality} 
the embedding of $\Gamma$ in $G$ extends
continuously to a homomorphism $\hat \varphi:G \to \RCom{\Gamma}{\Lambda} \cong H$,
where the last isomorphism follows from the fact that $H$ is reduced (no non-trivial
normal compact subgroups) using Lemma \ref{lem:homtotop}.
On the other hand,
clearly 
$G$ itself is a quotient
of ($G_\Gamma$ and)  $\tilde G_\Gamma$, 
 hence composing the two homomorphisms yields
the required extension of $\varphi$ to $\tilde G_\Gamma$.

For the opposite direction let $\Lambda < \Gamma$ be a commensurated subgroup.
Consider the increasing (transfinite) chain of subgroups $\Lambda _{\alpha}$
given by Proposition \ref{prop:chain}, terminating in the commensurated 
subgroup $\Lambda ' < \Gamma$ for which the group $H=  \RCom{\Gamma}{\Lambda}  $
is reduced. We may now apply 
Conjecture \ref{conj:conj} to the natural dense embedding of $\Gamma$ in $H$
(which satisfies the regularity condition as $H$ is reduced).
This readily 
implies that $\Lambda '$ is $S'$-arithmetic. To finish, one only needs to verify 
that a sequence $\Lambda  _{\alpha}$ as in Proposition \ref{prop:chain} 
can terminate
in a $S'$-arithmetic subgroup $\Lambda '$ only if all subgroups in this 
chain are in fact commensurable. 
(Indeed, if $\Lambda '$ is finitely
generated this is a standard noetherianity argument. Otherwise the same argument
 shows that $\Lambda$ virtually contains an infinite finitely 
generated $S''$-arithmetic subgroup. Now one combines Proposition~\ref{prop:arithmetic} 
above, with the fact that
one cannot move upward between two different
$S$-arithmetic subgroups $\Lambda _1 < \Lambda _2 < {\bf G}(K)$ along a chain satisfying the 
conclusion of Proposition \ref{prop:chain}, unless they are commensurable, which can
be proved e.g. by noticing that $\Lambda _1$ is co-amenable in $\Lambda _2$
in the sense of Eymard~\cite{Ey}, hence the semi-simple 
group  $  \RCom{\Lambda _1}{\Lambda _2}   $
is amenable, hence compact, hence finite).
\end{demo}



\smallskip  

\noindent {\bf Margulis' superrigidity.} 
Conjecture~\ref{conj:conj} implies a sharp form of the celebrated superigidity theorem of 
Margulis for $S$-arithmetic groups,
 which is in fact strongly supported 
by Margulis work (including some less familiar 
aspects of it). 
Here one takes
$H$ to be an algebraic group over a local field $F$, noting that
$H<GL_n(F) < SL_{n+1}(F)$ and the latter always satisfies the second assumption
on $H$ of  Conjecture \ref{conj:conj}. It is not difficult to verify that the regularity 
condition in  Conjecture \ref{conj:conj} is not an obstacle in this well understood setting, thus Conjecture~\ref{conj:conj}
predicts that 
{\it all} homomorphisms into algebraic groups over local fields 
should extend continuously to the appropriate envelop. Note that this
is a {\it topological} superrigidity, unlike the so-called ``abstract'' one
which is known as a consequence of Margulis' work 
(see~\cite[Theorem 6, page 5]{Mar2}), and can also
be deduced from the CSP (when known to hold -- see~\cite{BMS}). 
 The knowledgeable 
reader will thus have noticed at this point
a greater level of generality compared to the usual formulations
of Margulis' result, 
which motivates the discussion to follow.

The first compatability issue arises from the standard assumption, made 
in Margulis' theorem, that the
homomorphism $\varphi$ defined on $\Gamma$ has unbounded image. In Margulis' theorem 
this is indeed 
necessary if one aims to extend $\varphi$ continuously to the product of
simple {\it non compact} algebraic groups $G_i$ hosting the lattice  $\Gamma$. 
The treatment of the case where $\varphi (\Gamma)$ is bounded splits naturally into the
 two: either ({\bf i}) the target
group is totally disconnected, in which case $\varphi$ ranges 
in a profinite group,
and the existence of continuous extension is handled by the CSP, or
({\bf ii}) the target group is connected (Lie group). Here the 
Zariski closure
of $\varphi (\Gamma)$ will be a compact real algebraic group, and
as is well 
known among experts and revealed in the {\it proof}
of Margulis' {\it arithmeticity} theorem, upon passing to a finite index subgroup 
$\varphi$ does extend
continuously 
to a real anisotropic place appearing among the factors of  
${\bf G}(\mathbb {A})$ (it is actually
${\bf G}(\mathbb {A})^{\widetilde{ }}$ that should be used 
to account for the finite index ``loss''). This fact can also be deduced from Margulis' Theorem in \cite[Thm 5, page 5]{Mar2}, recalling the key point that compact 
simple real matrix groups are algebraic -- 
see a more comprehensive discussion of this whole aspect 
in~\cite[Prop. 2.3]{La-Lu}.

Another key result of Margulis in our general context
is that {\it every} linear representation $\varphi$ is {\it semisimple} 
(see~\cite[Thm B, page 259]{Mar2}; over positive 
characteristic it is known, as predicted also by
Conjecture\ref{conj:conj}, that the image is finite).
This result of Margulis allows one to omit de facto  
any Zariski density/target
assumption on Im $\varphi$ in Margulis' Theorem. One
last important aspect encountered when comparing  Conjecture 
\ref{conj:conj} and 
Margulis' superrigidity is the precise local nature of the 
target group 
$H$. It is usually 
assumed to be of adjoint type, hence center free, in which case the center
of any Zariski dense homomorphism of a group into it must trivialize. In the setting of  Conjecture \ref{conj:conj} this means that the finite center
$\mu$ (as well as that of ${\bf G)}$ will remain ``invisible''. 
Moreover, when $H$ is
not of adjoint type Margulis gives a counter example 
in~\cite[(5.11) page 231]{Mar2} showing that 
superrigidity
does not necessarily hold.
However, by a careful analysis he demonstrates 
that when the ambient 
algebraic groups hosting $\Gamma$ are {\it simply connected}, 
one can always ``correct'' 
$\varphi$ by multiplying it with a finite homomorphism $f$ of $\Gamma$ to the
center of $H$, resulting in a new, extendable homomorphism. As by the CSP such $f$ should extend  
to $\tilde G _{\Gamma}$, by incorporating this into the extension one concludes that $\varphi$ is 
nevertheless expected to
 extend continuously from $\Gamma$ to $\tilde G _{\Gamma}$, even if it didn't extend to
the ambient group appearing in Margulis' original superrigidity theorem.

\smallskip
 Thus, it seems that all the essential ingredients
are supplied by Margulis' work (with help from the CSP) in order to support the general
version of topological superrigidity into all {\it algebraic} groups over local fields, 
suggested by  Conjecture~\ref{conj:conj}. 
 Note that our discussion here 
did not 
capture the subtle issue of target groups which are non 
algebraic (finite or infinite) covers of such a group. This is discussed next.

\smallskip

\noindent {\bf Central extensions of arithmetic groups.} Consider first
 a motivating question:
Fix $n \ge 3$, 
let ${\bf G=SL_{n}}$, $K=\QQ$, and let $H:= {\widetilde {SL_{n}(\RR)}}$ be the 
universal ($2:1$) cover of 
$SL_{n}(\RR)$. For which $S$ does there exist 
a non-trivial homomorphism $\varphi: \Gamma \to H$ for some $S$-arithmetic subgroup
$\Gamma <SL_n(K)$? 
 This question is included in the setting of  Conjecture \ref{conj:conj}, and is not covered by 
Margulis' theory, as $H$ is not algebraic. 
Of course, one could also take for $H$ central extensions of $SL_n(\QQ _p)$, products of
those, etc. 
When $n=3$ and $S=\{\infty\}$, Millson~\cite{Mi}, via
a rather elaborate argument (motivated by a question of Deligne and Sullivan 
resulting from their own geometric work~\cite{De-Su}), answered
positively this question, leaving open the case $n>3$ for this $S$ (which is 
answered positively below). He also established a {\it negative} answer for arbitrary $n$ and 
some real quadratic extensions $K$ of $\QQ$, when $S$ consists of the two infinite places. 
The latter result is included in
the fundamental Theorem of Deligne~\cite{De2},
to which we shortly return.
Deligne's work~\cite{De2} was later extended by Raghunathan~\cite{Rag3} to some 
anisotropic
groups (using the first proof of CSP in such cases by 
Kneser~\cite{Kn}), and was pivotal in Toledo's famous  
construction of non residually
finite fundamental groups of smooth projective manifolds~\cite{To}. 
Our goal here is to show that  Conjecture \ref{conj:conj} is 
supported
by the above mentioned work of Deligne,
and that together with the existence 
of Deligne's central extension~(\ref{adiagram}), 
it suggests a {\it complete} (generally conjectural)
analysis of this situation.

First let us clarify the relation between the motivating
question and its more 
familiar cousin: when is a lift of an ($S$-)arithmetic 
group $\Gamma < {\bf G}(K)$ 
to $\tilde \Gamma < H$ = central (non linear) extension of ${\bf G}(K)$, 
virtually torsion/center free, or a residually finite group?  
The intimate relation between this question and our approach is perhaps
best illustrated by the following concrete example:

\begin{theorem}
\label{thm:lift}
Fix $n \ge 3$ and let $\pi : \tilde G \to   SL_n(\QQ _p)$ be any non split central
extension of topological groups. Then the lifted group 
$\tilde \Gamma :=\pi ^{-1}(SL_n(\ZZ [{\frac{1}{p}}]))<\tilde G$ is {\it not} 
virtually torsion (or center)
free, nor is it residually finite. 

\end{theorem}
Recall that non-trivial extensions $\tilde G$ as above always exist (and are of number 
theoretic origin), with 
the universal one known to have the same order as that of the finite group of roots 
of unity in  $\QQ _p$ -- cf.~\cite{Pr} and the references therein.

\smallskip
\begin{demo}
We show that every finite index subgroup $\Gamma '< \tilde \Gamma$ has non-trivial (finite) center, 
which clearly implies all three statements  (in fact the argument gives a more precise information). Indeed, otherwise 
$\pi:\Gamma' \to \pi(\Gamma ')< SL_n(\ZZ [{\frac{1}{p}}])$ is a group isomorphism,
and we may apply Theorem \ref{thm:adelic} with $\Gamma =  \pi(\Gamma ')$,
$H=\tilde G$ and $\varphi = \pi^{-1}|_\Gamma '$. The closure of 
$\Gamma $ in $SL_n(\adele _f)$ is isomorphic to
$ SL_n(\QQ _p) \times L$ for some compact group $L$, and by non-splitness
of the original extension, the restriction of any continuous extension $\tilde \varphi$ to 
$SL_n(\QQ _p)$ must be trivial. Hence $\tilde \varphi$ must have bounded image in $H$ -- 
contradiction.
\end{demo}

\medskip

While Theorem  \ref{thm:adelic} is based on additional
ingredients, it is actually only the CSP
which is needed for the proof of Theorem \ref{thm:lift}. In a similar vein, 
taking this time $H= {\widetilde {SL_{n}(\RR)}}$ and $\pi:H \to SL_{n}(\RR)$,
Conjecture~\ref{conj:conj} also implies
that for the same group $\Gamma = SL_n(\ZZ[{\frac{1}{p}}])$
no finite index subgroup of $\pi^{-1}(\Gamma)=\tilde \Gamma < H$ 
 is center free or residually finite. As showed by Deligne~\cite{De2}, 
while this fact does not follow directly from the CSP, its proof does involve
similar ingredients to those showing in the CSP 
(which is not a surprise in light of the proof of 
Theorem \ref {thm:lift} above).

\smallskip

On the positive side, answering positively Millson's question from~\cite{Mi} mentioned earlier, 
observe that in the same last example there actually
{\it is} a splitting of the central extension 
$\pi : {\tilde SL}_{n}(\RR) \to SL_{n}(\RR)$ over a finite index subgroup  $\Gamma '< SL_n(\ZZ)$ (hence $\pi ^{-1}(SL_n(\ZZ))$ {\it is} 
virtually center free and residually finite). To see this apply 
part {\bf 3} of Proposition \ref{prop:split} to get, for a finite index  subgroup $\Gamma '$, 
$\tilde G _{\Gamma'} =  {\widetilde {SL_{n}(\RR)}} \times  
C$ where $C < {\bf G}(\mathbb {A}_f)$ is compact open. Now project 
in Deligne's extension~(\ref{adiagram}) the group
$\Gamma '\cong s(\Gamma ')< \tilde G _{\Gamma'}$ into the left factor
$ {\widetilde {SL_{n}(\RR)}} $, to get the required splitting.
(note that this suggests that the 
splitting of the above covering $\pi$
over a fixed lattice in $SL_{n}(\RR)$, 
{\it without} passing to a finite index subgroup, may be determined by a
 purely $p$-adic data -- see below.)
Similarly, one shows that over a finite index subgroup 
of $SL_n(\ZZ[{\frac{1}{p}}])$, the $2:1$ extension 
$$
p: (SL_n(\RR) \times SL_n(\QQ _p))^{\widetilde{ }}
 \to  SL_n(\RR) \times SL_n(\QQ _p)
$$
obtained by restriction of Deligne's extension~(\ref{adiagram}) to 
places $(p, \infty)$,
does split. 
One may  replace here $\ZZ[{\frac{1}{p}}]$ by $\ZZ[{\sqrt 2}]$ and $\QQ _p$ 
by $\RR$, with a similar outcome, but replacing $(SL_n(\RR) \times SL_n(\QQ _p)^{\widetilde{ }}$
with ${\widetilde SL_n(\RR)} \times {\widetilde SL_n(\QQ _p)}$ yields 
a $4:1$ extension for which, by Conjecture \ref{conj:conj},
the
lifted arithmetic group is {\it not} virtually torsion free. 

\smallskip

The general existence question of homomorphisms of $S$-arithmetic groups 
$\Gamma < {\bf G}(K)$
into groups $H$ which are central extensions of (products of) algebraic groups over a
local fields proceeds in a similar spirit, once ``projectivization'' and Margulis' superrigidity 
can be used to reduce the problem to 
the case where
$H$ is connected in the Schur theoretic sense, i.e., it is a perfect group, and 
is locally isomorphic to $G_S:=\Pi_{v \in S} {\bf G}(K_v)$ (one may, as we shall,
 omit here the non-archimedean
anisotropic places in $S$).
In the latter case,
one aims to obtain either an obstruction to a homomorphism into $H$ using
 Conjecture \ref{conj:conj} (a la part {\bf 1} of Proposition \ref{prop:split}), or 
a construction as illustrated with Millson's question for $SL_n(\ZZ)$, using Deligne's 
extension~(\ref{adiagram}) and the splitting in {\bf 3} of Proposition \ref{prop:split}.
It is quite remarkable that (at least in the $K$-isotropic case) with the help of Margulis' 
supperrigidity,
Deligne's main result in~\cite{De2} can be interpreted - even if it is not stated
in such terms -
as the assertion that these two approaches precisely complement each other.
In other words, for $K$-isotropic $\bf G$ (as always, absolutely simple and simply connected), 
the following holds: 

\smallskip
{\it Fix some $V^\infty \subseteq S$ with the usual higher rank assumption, and let 
$H$ be a (Schur theoretically-) connected group, locally (i.e. projectively) 
isomorphic to $G_S$ as above. Then
a non-trivial 
homomorphism of some $S$-arithmetic subgroup $\Gamma < {\bf G}(K)$ to the group $H$ exists,
if and only if $H$ is
a quotient of the group $ \tilde G^{(nc)}_{\Gamma } $  appearing in part {\bf 3} of
Proposition~\ref{prop:split}.}

\smallskip
Evidently, this is compatible with Conjecture~\ref{conj:conj}. Note however, that in the 
positive case the conjecture implies a stronger result, predicting when a given $S$-arithmetic 
group 
should, or should not, admit a non-trivial homomorphism into $H$, without allowing passage to a finite index 
subgroup. Clarifying this ``local-global'' picture would be of interest.

\smallskip
A different angle at 
Deligne's concise paper  can be found in section 1 of Raghunathan's~\cite{Rag3},
where he relaxes Deligne's original 
condition that ${\bf G}$ be quasi-split
to assuming that ${\bf G}(K)$ is perfect (now known to hold when $\bf G$ is $K$-isotropic), 
together with the centrality of the
$S$-congruence kernel for  ${\bf G}$ (used also by Deligne, and again known today at least in
the $K$-isotropic case). Of course,
many more details are skipped in the discussion above and
the general case remains conjectural.




\begin{thebibliography}{MM}
%
%
\bibitem{Abels} H.~Abels,
\newblock {\it Finite presentability of S-arithmetic groups, compact presentability of solvable groups}, \newblock Springer Lecture Notes\,{\bf
1261}, 1987.
%
%
\bibitem{BLS} H.~Bass, M.~Lazard and J-P.~Serre, {\it Sous-groupes
d'indices finis dans $SL(n , Z)$,} Bull. Amer. Math. Soc.\,{\bf
70}(1964), 385-392.
%
%
\bibitem{BMS} H.~Bass, J.~Milnor and J-P.~Serre, {\em Solution of
the congruence subgroup problem for $SL_n$ $(n \geqslant 3)$ and
$Sp_{2n}$ $(n \geqslant 2),$} Publ. Math. IHES\:{\bf 33}(1967),
59-137.

\bibitem{BaW} U. Baumgartner and G.\,A. Willis,
\emph{Contraction groups and scales of
automorphisms of totally disconnected locally compact
groups},
Israel J.\ Math.\ {\bf 142}\,(2004), 221--248.
%
%
\bibitem{BaWi:direction}
U. Baumgartner and G.A. Willis,
\newblock The direction of an automorphism of a totally disconnected locally compact group,
\newblock {\sl Math. Z.\/} {\bf 252} (2006), 393--428.
%
%
\bibitem{BerLen} G.\,M.\,Bergman and H.\,W.\,Lenstra,
\newblock Subgroups close to normal subgroups,
\newblock J. Algebra {\bf 127} (1989), 80--97. 





%
%
\bibitem{CarKel} D. Carter and G. Keller,
\newblock Bounded elementary generation of $SL_n({\mathcal O})$,
\newblock Amer. J. of Math. {\bf 105} (1983), 673--687.
%
%

\bibitem{Ch-Ve}
I.~Chatterji, T.~N.~Venkataramana,
{\it Discrete linear groups containing arithmetic groups},
preprint.


\bibitem{De2} P.~Deligne {\it Extensions centrales non residullement finies de grupes arithmetiques}  C.R Acad. Sci. Paris Ser A-B\:{\bf 287}(1978) no. 4,
203-208.

\bibitem{De1} P.~Deligne {\it Extensions Centrales de grupes alg\'ebriques simplement connexes et
cohomologie galoisienne} Publ. Math. IHES\:{\bf 84}(1996), 35-89.

\bibitem{De-Su} P.~Deligne and D.~Sullivan, 
{\it Fibres vectoriels complexes a groupe structural discret},
C. R. Acad. Sci. Paris Sér. A-B {\bf 281} (1975), no. 24, 81--83.


\bibitem{Er-Ra} I.~Erovenko and A.~Rapinchuk, 
{\it Bounded generation of $S$-arithmetic subgroups of isotropic
orthogonal groups over number fields},
J. Number Theory {\bf 119} (2006), no. 1, 28--48.

\bibitem{Ey}
P.~Eymard,
{\it Moyennes invariantes et representations unitaires}, 
Lect. Notes in Math. {\bf 300}, Springer, Berlin-New York, 1972.





\bibitem{Gl:skew} H. Gl\"ockner, Scale functions on linear groups over local
skew fields, J. Algebra {\bf 205} (1998), 525--541.
%
%
\bibitem{Gl:padic} H. Gl\"ockner, Scale functions on $p$-adic Lie groups,
Manuscripta Math. {\bf 97} (1998), 205--215.
%
%
\bibitem{GlWi} H.\,Gl\"ockner and G.\,A.\,Willis,
\newblock Topologization of Hecke pairs and Hecke C*-algebras,
\newblock Topology Proceedings {\bf 26}, (2001--2002), 565--591.
%
%
\bibitem{HLM}
K.\,H.\,Hofmann,  J.\,R.\,Liukkonen, M.\,W.\,Mislove,
\newblock \emph{Compact extensions of compactly generated nilpotent groups are pro-Lie},
\newblock Proc. Amer. Math. Soc.\ {\bf 84}\,(1982),  443--448.
%
%
\bibitem{HM}
K.\,H.\,Hofmann,  A.\,Mukherjea,
\newblock \emph{Concentration functions and a class of non-compact groups},
\newblock Math. Ann. {\bf 256} (1981), 535--548.
%
%
\bibitem{Hofmann} Karl H. Hofmann,
\newblock Characteristic subgroups in locally compact totally disconnected groups and their applications to a problem on random walks on locally compact groups,
\newblock {\sl Technische Hochschule Darmstadt}, Preprint No. 606, {\bf 1981}.
%
%
\bibitem{IwMa}
N.\,Iwahori, H.\,Matsumoto, 
\newblock \emph{On some Bruhat decompositions and the structure of the Hecke rings of $p$-adic Chevalley groups},
\newblock Publ. Math. I.H.E.S., {\bf 25} (1965), 5--48.
%
%
\bibitem{JRW}
W. Jaworski, J. M. Rosenblatt and G.A. Willis,
\newblock Concentration functions in locally compact groups,
\newblock {\sl Math. Ann.} {\bf 305}(1996), 673--691.
%
%
\bibitem{Kn} M.~Kneser, 
\newblock \emph{Normalteiler ganzzahliger
Spingruppen}, 
\newblock Crelle Journal\:{\bf 311/312}(1979),
191-214.

\bibitem{La-Lu} M.~Larsen and A.~Lubotzky, 
{\it Representation growth of linear groups},
J. Eur. Math. Soc. {\bf 10} (2008), no. 2, 351--390.

\bibitem{LMR} A.~Lubotzky, S.~Mozes, and M.~S.~Raghunathan, 
{\it The word and Riemannian metrics on lattices of semisimple groups,}
Inst. Hautes Études Sci. Publ. Math. {\bf 91} (2001), 5--53.

\bibitem{Lu-Se}
A.~Lubotzky, D.~Segal, 
{\it Subgroup Growth},
Progress in Mathematics, {\bf 212} Birkhäuser, Basel, 2003.

\bibitem{Lu-Zi} A.~Lubotzky and R.~Zimmer, 
{\it Arithmetic structure of fundamental groups and actions of semisimple
Lie groups,}
Topology {\bf 40} (2001), no. 4, 851--869.

\bibitem{Mar1} G.A.~Margulis, {\it Finiteness of quotients of
discrete groups,} Funct. Anal. Appl.\:{\bf 13}(1979), 28-39.

\bibitem{Mar2} G.A.~Margulis, {\it Discrete Subgroups of
Semisimple Lie Groups,} Springer-Verlag, 1991.

%
%
\bibitem{Ma} H.~Matsumoto, {\it Sur les sous-groupes
arithm\'etiques des groupes semi-simples d\'eploy\'es,} Ann. Sci.
Ecole Norm. Sup.\,(4) {\bf 2}(1969), 1-62.

\bibitem{Me} J.~Mennicke, {\it Finite factor groups of the
unimodular group,} Ann.\:Math.\,{\bf 81}(1965), 31-37.

\bibitem{Mi} J.~Millson {\it Real vector bundles with discrete structure group} Topology\:{\bf 18}(1979) no.1
83-89.

\bibitem{Moeller} R. M\"oller,
\newblock Structure theory of totally disconnected locally compact groups via graphs and permutations,
\newblock Canad. J. Math. {\bf 54} (2002), 795--827.
%
%

\bibitem{MontZip}
D. Montgomery and L. Zippin, Topological Transformation
Groups, Interscience Publishers, New York--London 1955, reprinted 1974.
%
%
\bibitem{Mo} C.~Moore, {\it Group extensions of $p$-adic and
adelic groups,} Publ. Math. IHES\:{\bf 35}(1968), 5-70.

\bibitem{Mor} D.~Morris-Witte, {\it Bounded generation of ${\rm SL}(n,A)$ (after
D. Carter, G. Keller, and E. Paige)}, New York J. Math.\,{\bf
13}(2007), 383-421.

\bibitem{PR} V.P.~Platonov, A.S.~Rapinchuk, {Algebraic Groups
and Number Theory,} Acad. Press, 1991.

\bibitem{Pr} G.~Prasad,
{\it Deligne's topological central extension is universal,} 
Adv. Math. {\bf 181} (2004), no. 1, 160--164.


\bibitem{PRag1} G.~Prasad and  M.S.~Raghunathan, {\it On the
congruence subgroup problem: Detrmination of the metaplectic
kernel,} Invent.\:Math.\,{\bf 71}(1983), 21-42.

\bibitem{PR1} G.~Prasad and A.S.~Rapinchuk, {\it Computation of the
metaplectic kernel,} Publ. Math. IHES\:{\bf 84}(1996), 91-187.

\bibitem{PR2} G.~Prasad and A.S.~Rapinchuk, {\it Developments on the
congruence subgroup problem after the work of Bass, Milnor and Serre},
Preprint.

\bibitem{Pra-Rap} G.~Prasad and A.S.~Rapinchuk,
{\it Weakly commensurable arithmetic groups and isospectral locally symmetric spaces},
Publ. Math. IHES {\bf 109} (2009), 113--184.


\bibitem{Rag2} M.S.\:\:Raghunathan, {\it On the congruence subgroup
problem,} Publ. Math. IHES {\bf 46}(1976), 107-161.

\bibitem{Rag3} M.S.~Raghunathan, {\it Torsion in cocompact
lattices in coverings of $\mathrm{Spin}(2,n),$} Math. Ann.\:{\bf
266}(1984), 403-419 (Corrigendum: ibid., {\bf 303}(1995), 575-578).

\bibitem{Rag1} M.S.~Raghunathan, {\it The congruence subgroup problem,}
Proc. Indian Acad. Sci. (Math. Sci.)\,{\bf 114}(2004), no. 4,
299-308.

\bibitem{Rap} A.S.~Rapinchuk, Congruence subgroup problem for algebraic
groups: old and new. Journ\'ees Arithm\'etiques, 1991 (Geneva).
Ast\'erisque \# {\bf 209}(1992), 73--84.

\bibitem{Rapi} A.S.~Rapinchuk,
{\em The Margulis-Platonov conjecture for $SL_{1,D}$ and 2-generation of finite simple groups},
Math. Z. {\bf 252} (2006), no. 2, 295--313.   



\bibitem{Sch1} G.\,Schlichting,
\newblock Operationen mit periodischen Stabilisatoren,
\newblock {\em Archiv der Math.\/}, {\bf 34} (1980), 97--99.

\bibitem{Se} J-P.~Serre, {\it Le probl\`eme des groupes de congruence
pour $SL_2,$} Ann.\:Math.\,{\bf 92}(1970), 489-527.

\bibitem{Sh1} Y.~Shalom,
{\it Rigidity of commensurators and irreducible lattices},
 Invent. Math., {\bf 141}(1):(2000) 1--54.
 
\bibitem{Sh3} Y.~Shalom,
{\it Bounded generation and Kazhdan's property (T)},
IHES Publ. {\bf 90} (2001) 145--168.

 \bibitem{Tav} O. I. Tavgen
 \newblock{\em Bounded generation of Chevalley groups over rings of algebraic $S$-integers},
 \newblock Math. USSR Izvestiya {\bf 36}(1991), 101--128. (Tom 54 (1990) 97--122.)
 
 \bibitem{To} D.~Toledo,
{\it Projective varieties with non-residually finite fundamental group,}
Inst. Hautes Études Sci. Publ. Math. {\bf 77} (1993), 103--119.

%
 %
 \bibitem{Tz} K.\,Tzanev,
 \newblock C*-alg\`ebres de Hecke at K-theorie,
 \newblock Doctoral Thesis, Universit\'e Paris~7 -- Denis Diderot (2000).
 %
 %
\bibitem{Tz1} K.\,Tzanev,
 \newblock Hecke C*-alg\`ebres and amenability,
 \newblock J. Operator Theory {\bf 50} (2003), 169--178.
 %
 %
 \bibitem{Venkat} T. N. Venkataramana,
 \newblock {\em Zariski dense subgroups of arithmetic groups},
\newblock J. Algebra {\bf 108} (1987), 325--339.
%
%

\bibitem{Wagner} F.\,O.\,Wagner,
\newblock Almost invariant families,
\newblock Bull. London Math. Soc., {\bf 30(3)} (1998), 235--240.
%
%
\bibitem{Wang} S.\,P.\, Wang,
\newblock Compactness properties of topological groups III,
\newblock Tran. Amer. Math. Soc. {\bf 209} (1978), 399--418. 

\bibitem{We} A.~Weil, {\it Sur certains grupes d'operateurs unitairs}, Acta Math. {\bf 111}(1964)
143-211.

\bibitem{Weis} B.~Weisfeiler,
{\em Strong approximation for Zariski-dense subgroups of semisimple algebraic groups},
Ann. of Math. (2) {\bf 120} (1984), no. 2, 271--315.


\bibitem{Wi94}
G.\,A. Willis,  {\em The structure of totally disconnected,
locally compact groups},
Math.\ Ann.\ {\bf 300}\,(1994), 341--363.
%
%
\bibitem{Wi01} G.\,A. Willis, {\em Further properties of the scale function on a totally disconnected group}, J.~Algebra, {\bf 237}\,(2001), 142--164.
%
%
\bibitem{Wi04}
G.\,A. Willis,  {\em Tidy subgroups for commuting automorphisms of
totally disconnected groups: An analogue of
simultaneous triangularisation of matrices},
New\ York\ J.\ Math. {\bf 10}\, (2004), 1--35.
%
%
\end{thebibliography}
\end{document}